\documentclass[manuscript]{acmart}

\usepackage{etoolbox} 
\usepackage{algorithm} 
\usepackage{algpseudocode} 
\usepackage{listings} 
\usepackage[nameinlink,capitalize]{cleveref} 
\usepackage{multirow} 
\usepackage{enumitem} 
\usepackage{underoverlap} 

\newcommand{\figdir}{.//}

\definecolor{contrastingred}{RGB}{180,0,0} 

\lstdefinelanguage[ppmd]{python}[]{python}{%
  emph={}
}
\definecolor{DarkBlue}{rgb}{0.00,0.00,0.55}
\definecolor{DarkRed}{rgb}{0.55,0.00,0.00}
\definecolor{DarkGreen}{rgb}{0.00,0.55,0.00}
\definecolor{Gray}{rgb}{0.95,0.95,0.95}
\definecolor{Purple}{rgb}{0.5,0.0,0.5}
\definecolor{Bittersweet}{rgb}{1.0,0.44,0.37}
\renewcommand{\vec}[1]{\boldsymbol{#1}}
\newcommand{\VDG}[2]{\ensuremath{\mathbb{V}_{#1,#2}^{(\text{DG})}}}
\newcommand{\VCG}[2]{\ensuremath{\mathbb{V}_{#1,#2}^{(\text{CG})}}}
\newcommand{\FDG}[2]{\ensuremath{\mathbb{F}_{#1,#2}^{(\text{DG})}}}

\newcommand{\jump}[1]{\ensuremath{[\![#1]\!]}}
\newcommand{\avg}[1]{\ensuremath{\{\!\{#1\}\!\}}}

\newtheorem{definition}{\protect\definitionname}[section]
\providecommand{\definitionname}{Definition}

\newtheorem{observation}{\protect\observationname}[section]
\providecommand{\observationname}{Observation}

\newtheorem{technique}{\protect\techniquename}[section]
\providecommand{\techniquename}{Technique}

\begin{document}

\acmJournal{TOMS}
\citestyle{acmauthoryear}
\title{Implementation techniques for multigrid solvers for high-order Discontinuous Galerkin methods}
\author{Sean~Baccas}
\affiliation{%
    \institution{Durham University}
    \department{Advanced Research Computing}
    \city{DH1 3LE Durham}
    \country{United Kingdom}
}
\authornote{%
    Current affiliation: United Kingdom Atomic Energy Authority, %
    OX14 3DB Abingdon, %
    United Kingdom
}
\orcid{0009-0008-5263-1230}
\author{Alexander~A.~Belozerov}
\affiliation{
    \institution{University of Bath}
    \department{Department of Mathematical Sciences and Institute for Mathematical Innovation (IMI)}
    \city{BA2 7AY Bath}
    \country{United Kingdom}
}
\orcid{0000-0002-9292-9540}
\author{Eike~H.~M\"{u}ller}
\affiliation{
    \institution{University of Bath}
    \department{Department of Mathematical Sciences}
    \city{BA2 7AY Bath}
    \country{United Kingdom}
}
\orcid{0000-0003-3006-3347}
\email{e.mueller@bath.ac.uk}
\author[5]{Tobias~Weinzierl}
\affiliation{
    \institution{Durham University}
    \department{Department of Computer Science}
    \city{DH1 3LE Durham}
    \country{United Kingdom}
}
\orcid{0000-0002-6208-1841}

\begin{abstract}
  Matrix-free geometric multigrid solvers for elliptic PDEs that have been discretised with Higher-order Discontinuous Galerkin (DG) methods are ideally suited to exploit state-of-the-art computer architectures. 
Higher polynomial degrees offer exponential convergence, while the workload fits to vector units, is straightforward to parallelise, and exhibits high arithmetic intensity. 
Yet, DG methods such as the interior penalty DG discreisation do not magically guarantee high performance: they require non-local memory access due to coupling between neighbouring cells and break down into compute steps of widely varying costs and compute character. 
We address these limitations by developing efficient execution strategies for $hp$-multigrid. Separating cell- and facet-operations by introducing auxiliary facet variables localizes data access, reduces the need for frequent synchronization, and enables overlap of computation and communication. Loop fusion results in a single-touch scheme which reads (cell) data only once per smoothing step. 
We interpret the resulting execution strategies in the context of a task formalism, which exposes additional concurreny. The target audience of this paper are practitioners in Scientific Computing who are not necessarily experts on multigrid or familiar with sophisticated discretisation techniques. By discussing implementation techniques for a powerful solver algorithm we aim to make it accessible to the wider community.

\end{abstract}



\maketitle

\section{Introduction}

%
%
Elliptic partial differential equations (PDEs) are ubiquitous in computational science and engineering. 
They describe stationary solutions, arise from (semi-)implicit time discretisations of evolutionary PDEs or express non-local global constraints. 
Efficient numerical algorithms are required to solve this prominent class of PDEs. Multigrid methods \cite{mccormick1987multigrid,trottenberg2001multigrid} are a family of hierarchical iterative solvers that are known to be algorithmically optimal:
For certain elliptic PDEs it can be 
proven that the computational cost (quantified by counting floating point operations) and storage requirements grow linearly in proportion to the number of unknowns in the discretised system (see e.g. \cite{hackbusch2013multi,reusken2008introduction}). For many other, more complicated problems this linear growth is observed empirically.
Although multigrid has been attested excellent potential for several exascale applications \cite{anzt2020preparing,baker2012scaling,kohl2022textbook,ibeid2020fft}, the implementation of the algorithm in practice remains non-trival. Multigrid often struggles to fit to modern architectures as these favour localized and continuous data access, perform best for small, dense matrix-vector multiplications (mat-vecs), and suffer from computational steps with reduced concurrency due to coarse grid solves.  
As a consequence, it is not trivial to translate the algorithmic optimality of multigrid directly into fast code on exascale machines. 
To make multigrid software fit for practical applications, it needs to satisfy a range of requirements:

\begin{enumerate}[label=(\alph*)]
  \item It needs to exploit vector processing capabilities e.g.~through AVX (x86) or SVE (ARM). 
  \item It has to have low matrix assembly cost, as the discretisation might change in each and every solution step e.g.~due to non-linearities, dynamic adaptive mesh refinement, or incremental numerical integration \cite{murray2021assembly}.
  \item It has to be memory-modest. For large scale applications---notably where multigrid is used as one building block only---memory quickly becomes a constraining factor. 
  \item It needs to scale, notably exploiting shared-memory parallelism. This is because on many modern machines the number of cores per node continues to increase while the node count stagnates.
\end{enumerate}

\noindent
This list is not comprehensive. It does, for example, not include issues that arise from heterogeneous architectures.

In this paper, we focus on a particular cocktail of multigrid ingredients to address these problems: 
High-order Discontinuous Galerkin (DG) methods \cite{reed1973triangular,johnson1986analysis} form the backbone of many simulation codes today, since they naturally result in systems composed of small, dense matrices that are well-suited for vectorisation. 
Their matrix-free implementations eliminate storage constraints and naturally adapt to changes in the matrix during the iterative solution of the problem. 
Exploiting both $h$-refinement in the mesh resolution and $p$-refinement in the polynomial degree of the local DG basis functions results in $hp$-multigrid, which combines the computational efficiency on the finest, most costly mesh, with raid convergence of traditional multigrid methods.
Finally, combinations of data and task parallelism offer the ``concurrency freedom'' to exploit clusters with nodes hosting hundreds of hardware threads.
The design space for combining these building blocks is vast. 
They need to be adapted and tailored to each other.

While sophisticated $hp$-multigrid is very powerful, its efficient implementation can be daunting, especially for those who are not familiar with advanced finite element discretisations. In this paper we aim to provide a collection of practically relevant techniques that will be useful for HPC practitioners without detailed knowledge of high order DG methods. To make our work accessible to a wider audience, we start by reviewing the discretisation and solver algorithms before explaining how they can be translated into efficient computer code.

%
%
\paragraph{Model problem and multigrid components}
In this paper, we concentrate on the solution of the Poisson equation, i.e.~the elliptic PDE

\begin{equation}
  \begin{aligned}
    -\Delta u(x) = -\sum_{j=1}^{d} \frac{\partial^2 u}{\partial x_j^2} & = f(x) \qquad\text{for all $x\in\Omega=[0,1]^d$ with Dirichlet BCs}\\
    u(x)         & = g(x) \qquad\text{for all $x\in \Gamma=\partial \Omega$}.
  \end{aligned}
  \label{eqn:introduction:poisson}
\end{equation}

\noindent
Here $\Delta$ is the $d$-dimensional Laplace operator and $f(x)$, $g(x)$ are given functions. 
Despite its simplicity, Equation \eqref{eqn:introduction:poisson} is a commonly used benchmark since its solution requires addressing the characteristic implementation challenges outlined above.

Our implementation uses the Interior Penalty Discontinuous Galerkin (DG) method \cite{wheeler1978elliptic,arnold1982interior,oden1998discontinuoushpfinite,baumann1999discontinuous,riviere1999improved,cockburn2009unified} on a (potentially dynamically adaptive) Cartesian mesh constructed through spacetrees \cite{weinzierl2011peano,weinzierl2019peano}. This results in a block-sparse linear system that we solve through a geometric multigrid solver with $hp$-coarsening \cite{siefert2014algebraic,kronbichler2018performance,bastian2019matrix}: 
Similar to \cite{bastian2019matrix}, we employ a single $p$-coarsening step to restrict the residual on the finest mesh to a low-order continuous Galerkin (CG) function space on the same mesh. The resulting error correction equation in CG space is solved (approximately) with a classic $h$-multigrid cycle. The computational bottleneck of the resulting scheme is the block-Jacobi smoother in the high-order DG space on the finest level.

As it is common practice for many applications in supercomputing, we employ a matrix-free approach where the global stiffness matrix is never assembled: 
instead, we exploit the homogeneity and isotropy of \eqref{eqn:introduction:poisson} to pre-compute a fixed number of small, representative matrices at compile time. 
From these, the cell- and face-local matrices are constructed and applied at each smoother iteration. 
All steps of the multigrid algorithm are expressed as mesh traversals which apply local matrix-blocks on-the-fly. 

The axis-alignment of the spacetree, a generalisation of the octree concept, would allow us to use factorisation techniques \cite{kronbichler2018performance} for an even more efficient on-the-fly assembly.
Since at higher discretisation orders the most significant storage costs arise from the DG matrix on the finest level, our pure matrix-free, rediscretisation-based approach can also be used within a geometric-algebraic combination of multigrid ingredients \cite{weinzierl2018quasi} or in combination with explicit assembly on coarser levels as in \cite{bastian2019matrix}; this can be advantageous for some applications.
\paragraph{Novelty and main achievements} 

The matrix-free implementation of the DG block-Jacobi smoother does not automatically guarantee high performance: the mesh-traversal requires non-local memory accesses since the unknowns in neighbouring cells are coupled. The naive implementation of the smoother requires more than one mesh traversal, which results in repeated access to the same memory. If we have complex dependencies on neighbouring cells this also results in potentially scattered memory accesses. 
On modern chip architectures, reading data from memory is about an order of magnitude more expensive than performing floating point operations. As a consequence, the advantages of the matrix-free approach are lost if we do not take care to avoid repeated and unstructured memory access.
A naive domain decomposition of the block-Jacobi iteration also might result in insufficient concurrency for modern manycore chips.
We tackle these challenges as follows:

\begin{description}
  \item[First,] we introduce a set of implementation techniques to translate DG (and CG) iterations into equivalent single-touch algorithms. As a consequence, the (dominant) data is read and written only once per mesh sweep. 
  \item[Second,] we introduce a task formalism to break down the components of interior penalty DG into small tasks. Scheduling some tasks for asynchronous execution with the runtime system exposes additional concurrency, which complements the  parallelism achieved with domain decomposition.
  \item[Third,] we combine the multigrid building blocks into state-of-the-art multigrid solvers. Convergence studies confirm the algorithmic efficiency for a range of mesh resolutions and polynomial degrees.
\end{description}

\noindent

%
%
\paragraph{Structure}

The remainder of the paper is organised as follows: we start by reviewing the mathematical background required for the formulation of the interior penalty DG discretisation in \cref{section:discretisation}. Our notion of efficiency, which will guide the implementation strategies in later sections, is defined in \cref{section:efficiency}. 
The methodological core contribution can be found in \cref{section:single-level}. There we introduce implementation techniques for the DG block-Jacobi iteration, which forms the computational bottleneck of the $hp$-multigrid algorithm. These methods are interpreted in the context of a task-language, which allows for further optimisations. We demonstrate how the block-Jacobi smoother can be integrated into an $hp$-multigrid algorithm in \cref{section:multigrid}. 
Convergence studies (\cref{section:results}) confirm that the algorithm is robust with respect to increases in mesh resolution and polynomial degree. This is complemented by an evaluation of the computational performance of our techniques (\cref{section:runtime}). A brief discussion and outlook in \cref{section:conclusion} summarise the main achievements and highlight directions for future work.

\section{Mesh geometry and function space discretisation}
\label{section:discretisation}
Without loss of generality we assume homogeneous Dirichlet boundary conditions, $g(x)=0$ for the solution of the Poisson problem in \eqref{eqn:introduction:poisson} on the $d$-dimensional unit box $\Omega=[0,1]^d$.
More realistic geometries require proper scaling, mesh distortion, immersed boundary or marker-and-cell techniques to be mapped onto such a mesh. They are beyond the scope of the present paper.

\subsection{Hierarchical Cartesian meshes}
Our code employs a spacetree data structure  with a subdivision factor of
three \cite{weinzierl2011peano,weinzierl2019peano} to divide the domain $\Omega$ into a computational fine level mesh $\Omega_h$, embedded into a hierarchy of coarser meshes: the hypercube $[0,1]^d$ that defines the domain can be considered as a trivial mesh consisting of a single cell with $2d$ facets and $2^d$ vertices.
Next, we split up the hypercube equidistantly into three parts along each
coordinate axis, and therefore obtain a regular, axis-aligned Cartesian mesh with $3^d$ cells, $4d \cdot 3^{d-1}$ facets and $4^d$
vertices.

We repeat the splitting recursively, yet decide independently for each of these hypercubes whether to refine or not. This results in a hierarchy of nested meshes which are embedded into each other.
The cubic cells naturally form a rooted tree with the original hypercube $[0,1]^d$ at the root.

\begin{figure}[htb]
  \begin{center}
    \includegraphics[width=0.35\textwidth]{\figdir/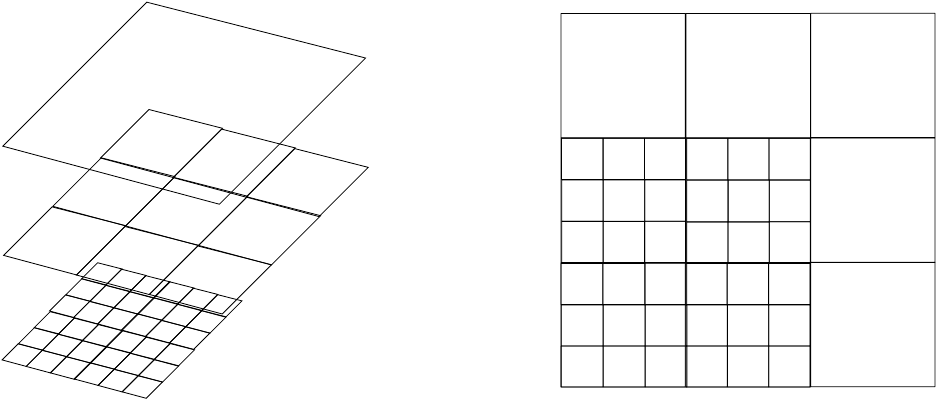}
    \hspace{1.2cm}
    \includegraphics[width=0.15\linewidth]{\figdir/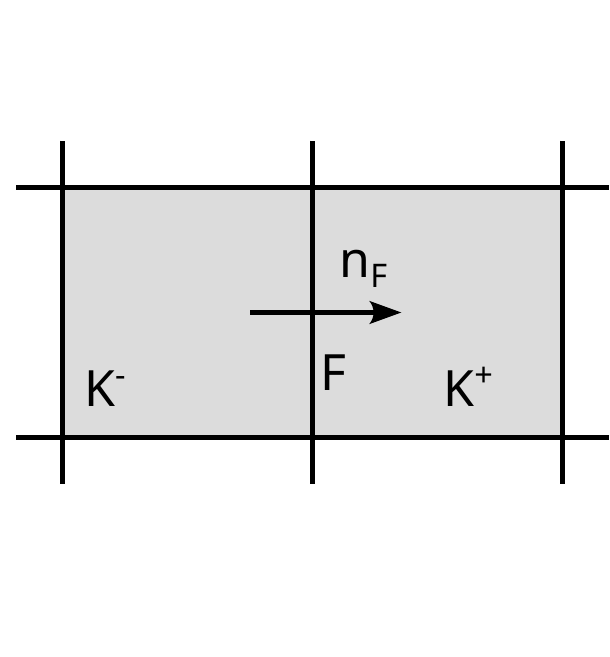}
  \end{center}
  \caption{
    Left: 
    We construct the spatial discretisation through a sequence of hypercubes
    (squares as we work in a $d=2$-dimensional setting in the illustration)
    embedded into each other.
    They form a spacetree.
    Middle: 
    The union of the unrefined leaf nodes of the tree forms an adaptive
    Cartesian mesh.
    Right:
    Any two neighbouring cells $K^-$, $K^+$ share one facet $F$ with associated facet normal $n_F$.
    \label{figure:discretisation:spacetree}
  }
\end{figure}

The union $\Omega_h = \bigcup_i K_i$ of all unrefined leaf-cells $K_i$, i.e.~cells that are not subdivided further, forms an adaptive Cartesian fine level mesh $\Omega _h$. 
Each $K_i$ is an open hypercube with the corresponding closure $\overline{K}_i$ where by construction $\bigcup_i \overline{K}_i = \Omega$ and $K_i\bigcap K_j=\emptyset$ for $i\neq j$. 
We define the skeleton $\mathcal{E}_h$ as the set of facets $F$ of the mesh:

\[
  \mathcal{E}_h = \{F\} = \{\overline{K}_i\cap
   \overline{K}_j, \quad \forall K_i, K_j\in \Omega_h, i\neq j \}.
\]

\noindent
The skeleton $\mathcal{E}_h$ can we divided into two disjoint sets:
boundary facets $\mathcal{E}_h^\partial$ and interior facets $\mathcal{E}_h^\text{i}$.
%
%
The set of vertices $\mathcal{V}_h$ of the mesh can be separated into interior and boundary vertices in a similar way.

For each facet $F\in\mathcal{E}_h$ we choose a unique unit normal vector $n_F$ which defines the orientation of the facet. 
Since the mesh is rectangular and axis-aligned, the $d$-dimensional vector $n_F$ has exactly one non-zero component and interior facets can be normalised to point in the direction of the coordinate axes. 
For boundary facets, we choose $n_F$ to coincide with the outward normal of the domain $\Omega$ (cmp.~convention in \cite{bastian2012algebraic}).

Each interior facet $F$ forms the boundary of exactly two cells, which we denote by $K^+$ and $K^-$ such that $n_F$ points from $K^-$ to $K^+$ (\cref{figure:discretisation:spacetree}), which we label as ``from left to right''. 
We will later associate unknowns with these left ($^-$) and right ($^+$) adjacent cells.
Let furthermore $\mathcal{N}(K)\subset \Omega_h$ be the set of facet-connected neighbours of a cell, i.e.~all cells $K' \not= K$ which share a facet with $K$. 
Analogously, the facets of a cell are denoted by $\mathcal{F}(K)\subset \mathcal{E}_h$, and we write $\mathcal{V}(K)$ for the set of its vertices.

Coarse level meshes, which are required for the $h$-multigrid algorithm, can be constructed in an analogous way by using cells on the coarser levels of the spacetree hierarchy.
\subsection{Function spaces}

\begin{figure}[htb]
  \begin{center}
    \includegraphics[width=0.95\textwidth]{\figdir/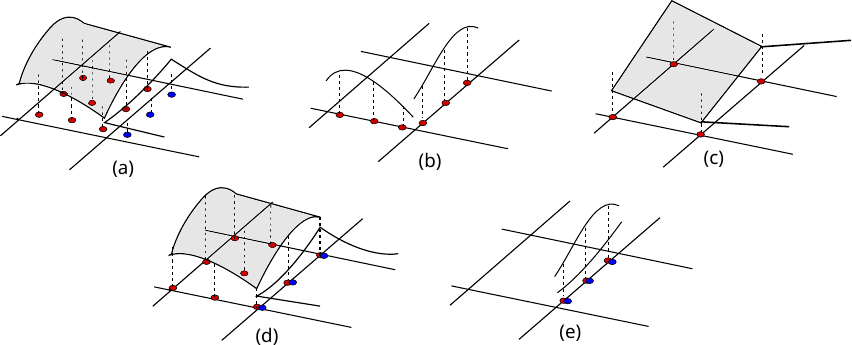}
  \end{center}
  \caption{
    Illustration of the fundamental function spaces used in this work.  Discontinuous Galerkin (DG) space $\VDG{h}{p=2}$ with Gauss-Legendre nodes (a), DG facet space $\FDG{h}{p=2}$ with Gauss-Legendre nodes  (b), lowest order piecewise linear Continuous Galerkin (CG) space $\VCG{h}{p=1}$ (c), DG space $\VDG{h}{p=2}$ with Gauss-Lobatto nodes (d) and composite DG facet space $\FDG{h}{p=2}\times \FDG{h}{p=2}$ with Gauss-Legendre nodes (e).
    \label{figure:discretisation:shape-functions}
  }
\end{figure}

We require three different families of function spaces: 
volumetric Discontinuous Galerkin (DG) spaces $\VDG{h}{p}$, DG facet spaces $\FDG{h}{p}$ and lowest order volumetric Continuous Galerkin (CG) spaces $\VCG{h}{p=1}$ (\cref{figure:discretisation:shape-functions}).

\subsubsection{Volumetric DG spaces}\label{sec:discretisation:volumetric_dg_spaces}
A discretisation of~\eqref{eqn:introduction:poisson} in the Discontinuous Galerkin framework aims to find solutions of the underlying PDE in the weak form in the function space $\VDG{h}{p}$ defined by
\begin{equation}
  \VDG{h}{p} = \{ u\in L_2(\Omega) : u|_K \in \mathbb{P}_p(K)\;\text{for all cells $K\in \Omega_h$}\}.
  \label{eqn:DG_functionspace}
\end{equation}

\noindent
Here and in the following $\mathbb{P}_p(S)$ is the space of multivariate polynomials of degree $p$ on the subdomain $S\subset \Omega$. 
A natural basis on each reference element consists of tensor-products of Lagrange polynomials whose nodes are Gauss-Legendre or Gauss-Lobatto points (\cref{figure:discretisation:shape-functions}~a,d). While the resulting functions in $\VDG{h}{p}$ are square integrable over $\Omega$, they are generally not continuous across interior facets where they might exhibit jumps. To represent functions from the DG space $\VDG{h}{p}$, we need to store $(p+1)^d$ values per fine grid cell $K \in \Omega _h$. Depending on the choice of basis, these degrees of freedom correspond to the evaluation of the function at the $(p+1)^d$ Gauss-Legendre or Gauss-Lobatto nodes of $K$.

\subsubsection{DG facet spaces}\label{sec:discretisation:dg_facet_spaces}
\label{subsection:discretisation:dg-facet-spaces}
Our implementation also requires the storage of functions on facets of the mesh. For this we define the DG facet space (\cref{figure:discretisation:shape-functions}~b)
\begin{equation}
  \FDG{h}{p} := \{\lambda\in L_2(\mathcal{E}_h):\lambda|_F \in \mathbb{P}_p(F)\;\text{for all facets $F\in\mathcal{E}_h$}\}.
  \label{eqn:discretisation:fdg_space}
\end{equation}

\noindent
%
Since each facet can be interpreted as a $d-1$-dimensional cell, functions in $\FDG{h}{p}$ can be stored similarly to the volumetric DG functions described in \cref{sec:discretisation:volumetric_dg_spaces}: On each facet the function is represented as the tensor-product of $d-1$ Lagrange polynomials, the nodes of which are the $p+1$ Gauss-Legendre and Gauss-Lobatto points in one dimension. This results in $(p+1)^{d-1}$ degrees of freedom per facet $F\in\mathcal{E}_h$.

Since functions in $\VDG{h}{p}$ can have jumps across facets, the projection of a function $u\in\VDG{h}{p}$ onto the facet space $\FDG{h}{p}$ hence is not well defined. To address this issue, we can project the solution separately for the left and right cell that touches the facet. The combination of these two projections can be stored in the product space $\FDG{h}{p} \times \FDG{h}{p}$ (\cref{figure:discretisation:shape-functions}~e); each nodal point of this space stores two values which correspond to the left and right projection of the volumetric DG function $u$.
 
\subsubsection{Volumetric CG spaces}\label{sec:discretisation:volumetric_cg_spaces}
The continuous Galerkin spaces $\VCG{h}{p} \subset \VDG{h}{p}$ contain those functions that can be represented by a $p$-dimensional polynomial in each cell while being continuous across facets. 
There are two ways of constructing such a continuous function space. 
We can either build the continuity directly into the construction of the (global) basis functions or use the corresponding volumetric DG space and enforce the continuity by synchronising the DG unknowns that correspond to a single CG unknown.
%

In this work, we use only the special case of $\VCG{h}{1}$, which consists of piecewise linear continuous functions over $\Omega_h$ (\cref{figure:discretisation:shape-functions}~c). 
In this case it makes sense to associate the unknowns directly within the vertices of the mesh, which encodes the continuity directly in the function space data structure.

\subsection{Weak formulation in discontinuous Galerkin space}
\label{subsection:dg-discretisations}
We now write down the discretiation of~\eqref{eqn:introduction:poisson} in the DG function space $\VDG{h}{p}$ defined in~\eqref{eqn:DG_functionspace}.
To simplify notation, we write volume- and surface-integrals as

\[
  \begin{aligned}
    \left( f \right)_{K}  & := \int_K f(x)\;dx\quad\text{for $d$-dimensional cells $K$, and }    \\
    \langle f \rangle_{F} & := \int_F f(x)\;ds\quad\text{for $d-1$-dimensional facets $F$}.
  \end{aligned}
\]

\noindent
Multiplication of the left-hand side of~\eqref{eqn:introduction:poisson} by a test-function $v\in \VDG{h}{p}$, integration over the domain $\Omega$ and integration by parts lead to

\begin{equation}
  \begin{aligned}
    \int_\Omega v (-\Delta u)\;dx = \sum_{K\in \Omega_h} \left(v (-\Delta u) \right)_{K}
    \mapsto \sum_{K\in \Omega_h} \left(\nabla v \cdot \nabla u \right)_{K} -
    \sum_{K\in \Omega_h} \langle v (n\cdot\widehat{\nabla u})\rangle_{\partial K}=:a(u,v),
  \end{aligned}
  \label{eqn:weak_form_naive}
\end{equation}

\noindent
where $n$ is the outward unit normal of cell $K$. 
The integration by parts requires the evaluation of the gradient $\nabla u$ on the surface of the cell. 
This term is not well defined since the solution is not continuous (let alone differentiable) across the facets. As a consequence we have to approximate the gradient with a \emph{numerical flux} $\widehat{\nabla u} \approx \nabla u$.
\subsubsection{Interior Penalty method}\label{sec:ip_method}
For DG methods this numerical flux is constructed by using the function representations in the neighbouring cells.
We adopt common notation from the DG literature to define the extrapolation of cell-based values to facets:
For some function $v\in \VDG{h}{p}$ which is continuous within each cell (but not across facets) we write

\begin{equation}
  v^\pm (x) := \lim_{\varepsilon\rightarrow 0^+} v(x\pm \varepsilon n_F)\qquad\text{for $x\in F$}\label{eqn:v_limits}
\end{equation}

\noindent
and define the jump and average as

\[
  \jump{v} = v^- - v^+ \quad \text{and} \quad   
  \avg{v} = \frac{1}{2}\left(v^-+v^+\right) \quad  \text{respectively.}
\]

\noindent
For boundary facets $\jump{v}:=v^-$, and we often simply write $v$,
since there is only one way of taking the limit in \eqref{eqn:v_limits}. 
This allows re-writing the second term in \eqref{eqn:weak_form_naive} as a sum over facets such that the weak form becomes

\begin{equation}
  a(u,v) = \sum_{K\in \Omega_h} \left(\nabla v \cdot \nabla u \right)_{K} -
  \sum_{F\in \mathcal{E}_h} \langle \jump{v} (n_F\cdot\widehat{\nabla u})\rangle_{F}.
  \label{eqn:weak_form_facetwise}
\end{equation}

\noindent
A naive choice for the numerical flux would be to take the average of the gradients from the adjacent cells, i.e.~$\widehat{\nabla u}=\avg{\nabla u}$. 
Unfortunately, this leads to an indefinite problem, as $a(u+w,v)=a(u,v)\Leftrightarrow a(w,v)=0$ for all $v\in\VDG{h}{p}$ where $w\in\VDG{h}{0}$ is an arbitrary function which is piecewise constant on the cells. 
To address this issue, additional stabilisation or regularisation terms have to be added to~\eqref{eqn:weak_form_facetwise}. 
The interior penalty formulation (see e.g.~\cite{bastian2012algebraic}) uses $\widehat{\nabla u}=\avg{\nabla u}$ and therefore extends~\eqref{eqn:weak_form_facetwise} into

\begin{equation}
  \begin{aligned}
    a(u,v) = & \sum_{K\in\Omega_h} \left(\nabla u\cdot\nabla v\right)_{K}
     + \sum_{F\in\mathcal{E}_h^\text{i}} \left(
    -\langle\jump{v}\avg{n_F\cdot \nabla u}\rangle_F
    +\theta \langle\jump{u}\avg{n_F\cdot \nabla v}\rangle_F
    \right.
    \\
     &
     \left.
    + \gamma_F \langle \jump{u}\jump{v}\rangle_F
    \right) 
     + \sum_{F\in\mathcal{E}_h^\partial} \left(
    -\langle v (n_F\cdot \nabla u)\rangle_F
    +\theta \langle u (n_F\cdot \nabla v)\rangle_F
    + \gamma_F \langle uv \rangle_F
    \right),
  \end{aligned}\label{eqn:weak_form_IP}
\end{equation}

\noindent
where $\theta >0$ and $\gamma_F >0$ are positive parameters which are usually chosen to be constant.
The boundary terms over $\mathcal{E}_h^\partial$ in \eqref{eqn:weak_form_IP} enforce the homogeneous Dirichlet boundary conditions weakly. 
As the solution converges to the smooth, true solution with finer and finer meshes, the magnitude of the jumps $\jump{u}$ and the value of $u$ on the boundary decrease. 
Hence, the formulation in \eqref{eqn:weak_form_IP} is consistent since the penalty terms vanish in the limit $h\rightarrow 0$.

The DG solution $u\in\VDG{h}{p}$ to \eqref{eqn:introduction:poisson} is then obtained by solving
\[
  a(u,v) = b(v) := \left(vf\right)_{\Omega_h}\qquad\text{for all $v\in\VDG{h}{p}$}.
\]

\subsubsection{Matrix representation}
As the DG method employs basis functions with cell-local support, the matrix arising from \eqref{eqn:weak_form_IP} decomposes over the mesh: the global matrix which couples all unknowns is a composite of many small local matrices. Each of these local matrices describes the coupling of unknowns in a single cell or between a pair of adjacent cells, respectively. 
The weak problem in \eqref{eqn:weak_form_IP} leads to the following block-sparse system of linear equations:

\begin{equation}
   \underbrace{
      \begin{pmatrix}
        A_{K_1 \gets K_1} & A_{K_1 \gets K_2} & A_{K_1 \gets K_3} & A_{K_1 \gets K_4} & \ldots \\
        A_{K_2 \gets K_1} & A_{K_2 \gets K_2} & A_{K_2 \gets K_3} & A_{K_2 \gets K_4} & \ldots \\
        A_{K_3 \gets K_1} & A_{K_3 \gets K_2} & A_{K_3 \gets K_3} & A_{K_3 \gets K_4} & \ldots \\
        A_{K_4 \gets K_1} & A_{K_4 \gets K_2} & A_{K_4 \gets K_3} & A_{K_4 \gets K_4} & \ldots \\
        \vdots            & \vdots            & \vdots            & \vdots            & \ddots \\
      \end{pmatrix}
    }_{=:A}
    \underbrace{
      \begin{pmatrix}
        \vec{u}^{(c)}|_{K_1} \\
        \vec{u}^{(c)}|_{K_2} \\
        \vec{u}^{(c)}|_{K_3} \\
        \vec{u}^{(c)}|_{K_4} \\
        \ldots               \\
      \end{pmatrix}
    }_{=:\vec{u}^{(c)}}
    =
    \underbrace{
      \begin{pmatrix}
        \vec{b}^{(c)}|_{K_1} \\
        \vec{b}^{(c)}|_{K_2} \\
        \vec{b}^{(c)}|_{K_3} \\
        \vec{b}^{(c)}|_{K_4} \\
        \ldots               \\
      \end{pmatrix}
    }_{=:\vec{b}^{(c)}}.
  \label{eqn:matrix_form_IP}
\end{equation}

\noindent
The global vector $\vec{u}^{(c)}= (\vec{u}^{(c)}|_{K_1}, \vec{u}^{(c)}|_{K_2}, \ldots)$ collects all DG unknowns such that the local vector $\vec{u}^{(c)}|_{K}$ contains the unknows within a particular cell $K$.
For each cell $K$ the small dense block-matrix $A_{K \gets K} $ describes the coupling of the unknowns in the cell with themselves. The couplings of unknowns between two different cells $K_i\neq K_j$ are described by the other matrix blocks $A_{K_i \gets K_j}$.
Since the matrices $A_{K_i \gets K_j}$ for $i\neq j$ arise from the facet integrals that describe the numerical flux and penalty terms in \eqref{eqn:weak_form_IP}, they are non-zero only if $K_i$ and $K_j$ are direct neighbours, i.e. if they share a common facet:

\[
A_{K_i \gets K_j} = 0 \quad \text{if } K_j \not \in \mathcal{N}(K_i) \text{ and } i\not =j.
\]

\noindent
The stiffness matrix $A$ defining the linear system $A \vec{u}^{(c)} = \vec{b}^{(c)} $ consequently is block-sparse and the \textit{global} system of equations \eqref{eqn:matrix_form_IP} can be written as a sequence of \textit{local} systems
\begin{equation}
  \forall\;\text{cells}\;K: \qquad
  A_{K\gets K'} \vec{u}^{(c)}|_{K} + \sum_{K'\in \mathcal{N}(K) }A_{K\gets K'}\vec{u}^{(c)}|_{K'} = \vec{b}^{(c)}|_{K}.
  \label{eqn:local_matrix_cells}
\end{equation}

\subsubsection{Block-Jacobi updates}

A simple solver for the equation system \eqref{eqn:matrix_form_IP} is the block Jacobi iteration. 
It updates the entries of $\vec{u}^{(c)} = (\vec{u}^{(c)}|_{K_1}, \vec{u}^{(c)}|_{K_2}, \ldots)$ per cell concurrently.
Let the superscript ``$^k$'' denote the iterate.
Given the current iterate $(\vec{u}^{(c)})^k =\left( (\vec{u}^{(c)})^k|_{K_1}, (\vec{u}^{(c)})^k|_{K_2}, \ldots \right)$, the values $(\vec{u}^{(c)})^k|_{K}$ in cell $K$ are updated to obtain $(\vec{u}^{(c)})^{k+1}|_{K}$ at the next iteration according to the rule

\begin{equation}
  \begin{aligned}
    \forall\;\text{cells}\;K: \qquad (\vec{u}^{(c)})^{k+1}|_K & = (\vec{u}^{(c)})^{k}|_K + \omega A_{K\leftarrow K}^{-1} \left(\vec{b}^{(c)} - A(\vec{u}^{(c)})^k \right) |_K \\
                                                              & =: (\vec{u}^{(c)})^{k}|_K + \omega A_{K\leftarrow K}^{-1} (\vec{r}^{(c)})_K,
  \end{aligned}
  \label{eqn:block-Jacobi-IP}
\end{equation}

\noindent
with the residual $\vec{r}^{(c)}$ given by

\[
  (\vec{r}^{(c)})|_K = \vec{b}^{(c)}|_K - A_{K\gets K} (\vec{u}^{(c)})^k|_K - \sum_{K'\in \mathcal{N}(K)}A_{K\gets K'}(\vec{u}^{(c)})^k|_{K'}.
\]

\noindent
The 
relaxation parameter $0 < \omega \leq 1$ in \eqref{eqn:block-Jacobi-IP} is chosen to guarantee and accelerate convergence. 
If the block-Jacobi update is used as a smoother in a multigrid method, 
it is necessary to pick $\omega$ such that errors that fluctuate rapidly across the domain are damped efficiently.

\subsection{Weak formulation in piecewise linear function space}
\label{sec:weak_form_piecewise_linear}
On the coarse levels of the multigrid hierarchy we also need to discretise  \eqref{eqn:introduction:poisson} in the piecewise linear function space $\VCG{h}{1}$. In this case the surface integrals in \eqref{eqn:weak_form_naive} cancel between neighbouring cells as their outward normals point in opposite directions. The weak form reduces to
\begin{equation}  
  a(u,v) = \sum_{K\in \Omega_h} \left(\nabla v \cdot \nabla u \right)_{K}\qquad\text{for $u,v\in\VCG{h}{1}$}.
  \label{eqn:weak_form_cg}
\end{equation}

\noindent
As discussed in \cite[Appendix A.2]{bastian2019matrix}, the same result is obtained by restricting the weak form of the interial penalty DG method in \eqref{eqn:weak_form_IP} to functions in the subspace $\VCG{h}{1}\subset \VDG{h}{p}$. 
Instead of the block-smoother \eqref{eqn:block-Jacobi-IP} we employ a point-Jacobi iteration in this case. 

\[
(\vec{u}^{(c)})^{k+1}|_V  = (\vec{u}^{(c)})^{k}|_V + \omega D^{-1} (\vec{r}^{(c)})_V.
\]

\noindent
The degrees of freedom of $\vec{u}^{(c)}$ are now given by the values at the vertices of the mesh.
The diagonal matrix $D$ is constructed by accumulating the diagonal values of the cell-local matrices for all cells that touch a particular vertex $V$,
while the residual $(\vec{r}^{(c)})_V$ results from the matrix-vector products of these matrices with the current solution $(\vec{u}^{(c)})^{k}$.

\section{Criteria for an efficient implementation}
\label{section:efficiency}

\begin{algorithm}
  \caption{Matrix-free block-Jacobi iteration. Input: initial solution $\vec{u}^{(c)}_0$, right hand side $\vec{b}^{(c)}$, relaxation parameter $\omega$, number of iterations $n_{\text{iter}}$, tolerance $\epsilon$. Output: solution $\vec{u}^{(c)}$ after $n_{\text{iter}}$ iterations.}
  \label{alg:vanilla-interior-penalty-single-level}
  \begin{algorithmic}[1]
    \State{Set $\vec{u}^{(c)}\gets \vec{u}^{(c)}_0$}
    \For{$k=1,2,\dots,n_{\text{iter}}$}
    \State{Set $\vec{u}^{(c)}_{\text{old}}\gets\vec{u}^{(c)}$}
      \label{alg:vanilla-interior-penalty-single-level:back-up-solution}
    \For{every cell $K \in \Omega_h$}
    \Comment{on-the-fly assembly to stay matrix-free}\label{alg:vanilla-interior-penalty-single-level:iteration:start}
    \State {$\vec{r}^{(c)}|_K \gets \vec{b}^{(c)}|_K -  A_{K\gets K} \vec{u}^{(c)}_{\text{old}}|_K$}\label{alg:vanilla-interior-penalty-single-level:residual}
    \Comment {cell-local contribution to residual}
    \For{every neighbouring cell $K'\in\mathcal{N}(K)$}
    \State {$\vec{r}^{(c)}|_K \gets \vec{r}^{(c)}|_K - A_{K\gets K'} \vec{u}^{(c)}_{\text{old}}|_{K'}$}
      \Comment {contribution from neighbouring}
      \label{alg:vanilla-interior-penalty-single-level:loop-over-neighbours}
    \EndFor
      \Comment {cells to residual}
    \State {$\vec{u}^{(c)}|_K \gets \vec{u}^{(c)}_{\text{old}}|_K + \omega  A_{K\gets K}^{-1} \vec{r}^{(c)}|_K$}\label{alg:vanilla-interior-penalty-single-level:solution_update}
    \Comment{update state}
    \EndFor\label{alg:vanilla-interior-penalty-single-level:iteration:end}
    \EndFor
    \State{\Return $\vec{u}^{(c)}$}
  \end{algorithmic}
\end{algorithm}

We employ the block-Jacobi iteration \eqref{eqn:block-Jacobi-IP} as smoother of our $hp$-multigrid algorithm. A vanilla implementation for the iterative solution of \eqref{eqn:matrix_form_IP} is shown in \cref{alg:vanilla-interior-penalty-single-level}.
If it is applied as solver rather than a smoother with a fixed iteration count, it might be supplemented with an early termination criterion. The can for example measure the relative reduction of some norm of the residual or the change in the solution from one iteration to the next.

An efficient implementation of the block-Jacobi iteration has to map well to the underlying hardware.
Since the majority of the runtime of the multigrid algorithm will be spent on the finest level, it makes sense to concentrate on formulating criteria that guarantee the efficiency of the DG smoother, bearing mind that these requirements apply to the coarse grid corrections, too.

\subsection{Efficiency criteria}

\paragraph*{Vector efficiency}

The vanilla implementation relies on the multiplication with small dense local matrices.

\begin{definition}\label{def:vectorisation}
 A \emph{core-efficient} implementation is able to exploit the vector capabilities of a compute core for all of its fundamental linear algebra building blocks.  
\end{definition}

\noindent
For this, a realisation should coalesce memory access subject to sufficiently wide vector instructions.

\paragraph*{Concurrency}

The update in \eqref{eqn:block-Jacobi-IP} can be carried out in parallel for all cells $K$. We therefore postulate that a good implementation should preserve this concurrency in the following sense:

\begin{definition}\label{def:concurrency}
 An \emph{efficient parallel} implementation of the block-Jacobi iteration has a concurrency level that equals (or exceeds) the number of geometric cells.  
\end{definition}

\paragraph*{Assembly overhead}

To ensure efficiency, we also need to take into account storage requirements and the cost of transferring data between memory and the compute unit.  This motivates the following:

\begin{definition}\label{def:matrix_free}
  A \emph{matrix-free} approach does not assemble the global operator $A$ at any point.
\end{definition}

\noindent
Implicitly, working matrix-free eliminates any expensive ``warm-up'' assembly phase predating the actual solve.

\paragraph*{Memory overhead}

While the implementation of the iteration over the mesh in \cref{alg:vanilla-interior-penalty-single-level} is natural, its memory footprint is potentially twice as large as the data carrying the information that we are interested in, namely the solution:
we store both the current solution $\vec{u}^{(c)}$ and the previous iterate $\vec{u}_{\text{old}}^{(c)}$. For high polynomial degrees the number of unknowns associated with a single facet is a factor $p+1$ smaller than the number of degrees-of-freedom associated with a cell. Hence, only the latter is relevant for the memory footprint, and we postulate:

\begin{definition}\label{def:memory_overhead}
  Any \emph{memory overhead} that increases in direct proportion to the memory footprint of the solution $\vec{u}^{(c)}$ as the mesh is refined is significant.
\end{definition}

\paragraph*{Data transfer volume}

In the vanilla implementation given by \cref{alg:vanilla-interior-penalty-single-level}, the numerical flux (encoded in the matrix $A_{K\gets K'}$ for $K'\neq K$) is evaluated twice per facet, as each facet contributes to both of its neighbours. 
While we assume that computations on modern hardware are cheap, the duplicated evaluation also implies that we have to access the data of the two adjacent cells of every interior face twice.
In a parallel implementation, the same argument implies that volumetric halo data has to be exchanged in every smoothing step. This leads to the following:

\begin{definition}\label{def:single_touch}
  A \emph{single-touch} implementation is an implementation where we read and write each piece of data at most once per mesh traversal.
\end{definition}

\noindent
In practice, any stencil-like update of cells cannot avoid repeated data access, as the code has to couple neighbouring cells.
However, as we will see below, temporary data can help to decouple memory accesses and, at the same time, result in a significantly smaller memory footprint.
This motivates a slightly relaxed version of \cref{def:single_touch}:

\begin{definition}\label{def:weak_single_touch}
  In a \emph{weak single-touch} implementation the single-touch policy is only enforced strictly for all volumetric data, i.e. unknowns associated with cells. Facet data in contrast can be read or written more than once.
\end{definition}

\subsection{Assessment of vanilla implementation}

While the implementation of the iteration over the mesh in \cref{alg:vanilla-interior-penalty-single-level} is straightforward, it does not meet all the criteria listed in \cref{section:efficiency}.
The computationally most expensive components of \cref{alg:vanilla-interior-penalty-single-level} are multiplications of small dense local matrices such as $A_{K\gets K}$, $A_{K\gets K'}$ and $A_{K\gets K}^{-1}$ with local vectors such $(\vec{u}^{(c)})^{k}|_K$ and $\vec{r}^{(c)}|_K$. These can be implemented as BLAS routines which implicitly meet the requirements of \cref{def:vectorisation} as long as the degrees of freedom per cell or facet, i.e.~the vectors in \eqref{eqn:block-Jacobi-IP} and subsequent formulae, are stored continuously in memory.
However, for the Gauss-Lobatto basis only a small subset of basis functions are non-zero on the boundary. As a consequence, it is only necessary to access a small number of unknowns from neighbouring cells when evaluating the facet integrals in \eqref{eqn:weak_form_IP}. This introduces scattered, sparse memory access.

While the evaluation of the cell contributions can be parallelised over the cells in line with \cref{def:concurrency}, we note that data is backed up into a helper variable $\vec{u}^{(c)}_{\text{old}}$ (\cref{alg:vanilla-interior-penalty-single-level}, line \ref{alg:vanilla-interior-penalty-single-level:back-up-solution}). The storage of $\vec{u}^{(c)}_{\text{old}}$ increases the memory overhead and violates the requirements in \cref{def:memory_overhead}. Whenever we update the solution in a cell $K$, we use the unknowns associated with all adjacent cells from the previous iteration. Overwriting the solution in the current cell would allow the use of a single vector $\vec{u}^{(c)}$, but this will modify the algorithm since the neighbouring cells $K'$ would use different data when it is their turn to execute inter-cell terms in the loop over $\mathcal{N}(K')$. This introduces data-dependencies and the resulting Gauss-Seidel iteration is no longer inherently parallel. More importantly, a generalisation to non-linear scenarios (which we might want to consider in the future) is problematic: in this case the numerical flux is not guaranteed to be consistent when it is computed across the same facet from the two neighbouring cells in their respective updates.
Unfortunately, the backup process requires some synchronisation if the underlying loop over the cells is parallelised and hence does not strictly exhibit the full concurrency level any more unless we separate the backup process in a completely different mesh sweep of its own:
  
The block-Jacobi iteration can be carried out without assembling the \textit{global} matrix $A$ in \eqref{eqn:matrix_form_IP} (cmp.~\cref{def:matrix_free}): 
we only need the \text{local} matrix blocks $A_{K\gets K}$ and $A_{K\gets K'}$ for $K'\in\mathcal{N}(K')$, as well the inverse of $A_{K\gets K}$ to compute the matrix-vector product $A_{K\gets K}^{-1} \vec{r}^{(c)}|_K$. Note that -- since both $A_{K\gets K}$ and its LU-factors are dense -- there is no advantage in storing these factors instead of $A_{K\gets K}^{-1}$, and multiplication with the latter can be implemented efficiently as a BLAS \texttt{dgemv} operation.  For homogeneous and isotropic problems such as the Poisson equation in \eqref{eqn:introduction:poisson}, none of the local block-matrices will vary across the domain and hence they can be precomputed once for a representative cell at the beginning of the run.

It is worth quantifying the reduction in storage requirements that a matrix-free implementation achieves: bearing in mind that $(p+1)^d$ is the number of DG unknowns per cell and $N_c$ denotes the number of grid cells, then instead of storing and reading $\mathcal{O}(N_c (p+1)^{2d})$ matrix values, only a small number of representative matrices with $\mathcal{O}((p+1)^{2d})$ entries have to be stored. Hence, for sufficiently large grids, the cost of storing and reading the matrix is neglegible compared to the $\mathcal{O}(N_c(p+1)^d)$ storage cost of the unknowns themselves.
On modern hardware architectures, re-computing the small block-matrices on-the-fly every time
might further reduce the runtime, in particular for higher discretisation order \cite{muthing2017high,kronbichler2018performance,bastian2019matrix}. 

Finally, the loop over neighbours (line~\ref{alg:vanilla-interior-penalty-single-level:loop-over-neighbours}) induces repeated reads of volumetric data:
The data in each cell  $\vec{u}^{(c)}_{\text{old}}$ is required up to $2d+1$ times per iteration as input to the residual and update calculations. This violates the weak single-touch criterion in \cref{def:weak_single_touch}.

\section{An efficient implementation of the high-order block-Jacobi smoother}
Having identified the weaknesses of the naive implementation in \cref{alg:vanilla-interior-penalty-single-level}, we now discuss alternative approaches which avoid these issues.
\label{section:single-level}

\subsection{DG with the Interior Penalty method}\label{sec:implementation_with_facet_variables}

Facet integrals such as

\begin{eqnarray}
\langle \jump{v}\avg{n_F\cdot \nabla u}\rangle_F & = &
\frac{1}{2}\int_F \left( v^- (n_F\cdot\nabla u^-)  +  v^- (n_F\cdot\nabla u^+)
- v^+ (n_F\cdot \nabla u^-)  -  v^+ (n_F\cdot\nabla u^+)\right)\;ds
  \label{eqn:single-level:facet-integrals}
\end{eqnarray}

\noindent
in the weak formulation \eqref{eqn:weak_form_IP} are the reason why naive implementations of the block-Jacobi smoother for the interior penalty DG method struggle to exhibit the full concurrency level (\cref{def:concurrency}), impose significant memory overhead (\cref{def:memory_overhead}), and are not (weakly) single touch (\cref{def:single_touch}, \cref{def:weak_single_touch}):
The integral over facet $F$ in \eqref{eqn:single-level:facet-integrals} accepts inputs $u^+$, $u^-$ from two adjacent cells $K^+$, $K^-$ that touch the facet (repeated read of volumetric data), and it writes back into both cells since $v^+$, $v^-$ are the test functions on both sides of $F$ (concurrent write access requiring the back-up of volumetric data and synchronisation). 
To overcome these issues and to design an efficient, single-touch implementation without overhead, we rely on a combination of several techniques. The first two techniques avoid the direct coupling of adjacent cells by introducing temporary fields on the facets:

\begin{technique}
  \label{technique:left-right-face-projections}
  We introduce \emph{projection variables} on the facets to explicitly store the extrapolation $u^+$, $u^-$ of the solution in $K^+$, $K^-$.
  In addition, we might also store the projection of other quantities such as the normal derivatives $n_F\cdot \nabla u^+$ and $n_F\cdot \nabla u^-$. All these additional variables will be stored in the dof-vectors $\vec{u}^{(+)}$, $\vec{u}^{(-)}$, each of which might represent more than one function.
\end{technique}

\begin{technique}
  \label{technique:numerical-flux-outcome}
  We introduce \emph{flux variables} on the facets to store combinations of the projection variables, these might for example represent numerical fluxes such as $\avg{n_F\cdot \nabla u}=\frac{1}{2}((n_F\cdot \nabla u^+)+(n_F\cdot\nabla u^-))$ and jumps in the solution $\jump{u}=u^--u^+$. All flux variables are stored in the single dof-vector $\vec{w}^{(f)}=\vec{w}^{(f)}(\vec{u}^{(+)},\vec{u}^{(-)})$.
\end{technique}

\noindent
\cref{technique:left-right-face-projections} and \cref{technique:numerical-flux-outcome} allow us to replace the direct coupling between neighbouring cells $K^+$ and $K^-$ by an indirect coupling which is realised in three steps: first, data is written to the \emph{projection variables} $\vec{u}^{(+)}$, $\vec{u}^{(-)}$ on facet $F=K^+\cap K^-$. This can be done independently by both cells $K^+$, $K^-$ since the \emph{projection variables} are independent for each pair $(K,F)$. Next, the data is combined  into the \emph{flux variables} $\vec{w}^{(f)}$ which can be done independently on each facet. In the final step data stored in $\vec{w}^{(f)}$ is used to update the fields in $K^+$ and $K^-$.

\subsubsection{Formalisation through spurious facet unknowns}

The additional variables in \cref{technique:left-right-face-projections} and \cref{technique:numerical-flux-outcome} can be constructed for the interior penalty formulation in \eqref{eqn:weak_form_IP}. 
For this we introduce the following fields on each facet in addition to the DG field $u\in \VDG{h}{p}$ in each cell:

\begin{itemize}
  \item Two scalar-valued fields $\widetilde{u}^-,\widetilde{u}^+\in\FDG{h}{p}$ which represent the projection of $u$ from the two neighbouring cells onto the facet.
  \item Two scalar-valued fields $\widetilde{u}'^-,\widetilde{u}'^+\in\FDG{h}{p}$ which represent the projection of the normal derivative $\nabla u\cdot n_F$ from the two neighbouring cells onto the facet.
  \item Two scalar-valued fields $w^f,w'^f\in\FDG{h}{p}$ which are linear combinations of $\widetilde{u}^-$, $\widetilde{u}^+\in\FDG{h}{p}$ and $\widetilde{u}'^-$, $\widetilde{u}'^+\in\FDG{h}{p}$, respectively and which represent the numerical flux $\avg{\nabla u\cdot n_F}$ and jumps $\jump{u}$ that appear in stabilisation terms of the interior penalty DG formulation in \eqref{eqn:weak_form_IP}.
\end{itemize}

\noindent
More specifically, the projected fields $\widetilde{u}^\pm:=\mathcal{P}(u)\in\FDG{h}{p}$ and $\widetilde{u}'^\pm:=\mathcal{P'}(u)\in\FDG{h}{p}$ are defined as
\begin{subequations}
  \begin{align}
    \widetilde{u}^\pm(x) &= \mp \lim_{\varepsilon\rightarrow 0+}u(x\pm \varepsilon n_F),\label{eqn:solution_projection}\\
    \widetilde{u}'^\pm(x)  &= \lim_{\varepsilon\rightarrow 0^+} \nabla u(x\pm\varepsilon n_F) \cdot n_F.
      \label{eqn:gradient_projection}
  \end{align}
\end{subequations}

\noindent
On the boundary, only $\widetilde{u}^-$, $\widetilde{u}'^-$ will be non-zero.
We assume that $w^f$ and $w'^f$ can be expressed as  linear combinations of the corresponding projections:
\begin{equation}
  w^f = \begin{cases}
    \frac{1}{2}\left(\widetilde{u}^- + \widetilde{u}^+\right) \\[1ex]
    \widetilde{u}^-
  \end{cases},\qquad
  w'^f = \begin{cases}
    \frac{1}{2}\left(\widetilde{u}'^-+\widetilde{u}'^+\right) & \text{on interior facets}      \\[1ex]
    \widetilde{u}'^-                               & \text{on the boundary facets}.
  \end{cases}\label{eqn:numerical_fluxes}
\end{equation}
It is easy to see that the weak formulation in \eqref{eqn:weak_form_naive} with the bilinear form $a(u,v)$ in \eqref{eqn:weak_form_IP} is equivalent to
\begin{equation}
  \begin{aligned}
    b(u,v) = a(u,w^f,w'^f,v) = \sum_{K\in\Omega_h} \left(\nabla u\cdot\nabla v\right)_{K}
     & + \sum_{F\in\mathcal{E}_h^\text{i}} \left(
    -\langle\jump{v}w'^f\rangle_F
    +2\theta \langle w^f \avg{\nabla v\cdot n_F}\rangle_F
    + 2\gamma_F \langle w^f\jump{v} \rangle_F
    \right)                                       \\
     & + \sum_{F\in\mathcal{E}_h^\partial} \left(
    -\langle v w'^f\rangle_F
    +\theta \langle w^f (\nabla v\cdot n_F)\rangle_F
    + \gamma_F \langle w^f v \rangle_F
    \right),
  \end{aligned}\label{eqn:weak_form_IP_facets}
\end{equation}
provided the system is closed with the definition of the projections \eqref{eqn:solution_projection}, \eqref{eqn:gradient_projection} and of the numerical flux in \eqref{eqn:numerical_fluxes} which express $w^f$, $w'^f$ in terms of the DG variables $u$ (the factor two front of the penalty terms on interior facets arises since $w^f=\frac{1}{2}\jump{u}$). The bilinear form in \eqref{eqn:weak_form_IP_facets} can be written as a sum over cells
\begin{equation}
  a(u,w^f,w'^f,v) = \sum_{K\in\Omega_h} \left\{\left(\nabla u\cdot\nabla v\right)_{K}
  +\sum_{F\in\mathcal{F}(K)}\sigma \left(
  \langle v^{\sigma} \, w'^f \rangle_{\partial K}
  -\theta \langle w^f \, (\nabla v^{\sigma}\cdot n^{\sigma}) \rangle_{\partial K}
  - \overline{\gamma}_F \langle w^f \, v^{\sigma} \rangle_{\partial K}\right)
  \right\}.
  \label{eqn:weak_form_IP_facet_spaces}
\end{equation}
where $\mathcal{F}(K)\subset \mathcal{E}_h$ is the set of all facets of a given cell $K$. The sign $\sigma=\sigma(K,F)=-n\cdot n_F\in\{+,-\}$ (which implicitly depends on the cell $K$ and facet $F$) is negative if the outward normal $n$ of cell $K$ on facet $F$ is identical to $n_F$ and positive otherwise. The penalty parameter $\overline{\gamma}_F$ is identical to $\gamma_F$ on boundary facets, we have that $\overline{\gamma}_F=2\gamma_F$ on interior facets. In our implementation we set $\theta$ and $\overline{\gamma}_F$ to constant values, which implies that the penalty parameter $\gamma_F$ in \eqref{eqn:weak_form_IP} will differ between interior and boundary facets.

\paragraph{The linear equation system}

Let the \emph{projection variables} from \cref{technique:left-right-face-projections} be represented by the global dof-vector $\vec{u}^{(\pm)}$ of the projections $\widetilde{u}^{+},\widetilde{u}^{-},\widetilde{u}'^{+}$, $\widetilde{u}'^{-}$ on each facet. Let further the \emph{flux variables} from \cref{technique:numerical-flux-outcome} be represented by the global dof-vector $\vec{w}^{(f)}$ of the fluxes $w^f,w'^f$. We can then write \eqref{eqn:weak_form_IP_facet_spaces} in matrix form as

\begin{equation}
  A_{c\gets c} \vec{u}^{(c)} + A_{c\gets f} \vec{w}^{(f)} = \vec{b}^{(c)}\label{eqn:cell_equation}
\end{equation}

\noindent
where $A_{c\gets c}$ is the discretisation of the volume term $(\nabla u\cdot \nabla v)_\Omega$ and $A_{c\gets f}$ describes the couplings from facet-unknowns to cell-unknowns.
Similarly, \eqref{eqn:solution_projection}, \eqref{eqn:gradient_projection} and \eqref{eqn:numerical_fluxes} can be written as

\begin{xalignat}{2}
    \vec{u}^{(\pm)} & = A_{f\gets c} \vec{u}^{(c)},&
    \vec{w}^{(f)}   & = A_{f\gets f}\vec{u}^{(\pm)},\label{eqn:facet_equation}
\end{xalignat}
\noindent
where in the first equation we have multiplied by the inverse of the mass matrix of the space $\FDG{h}{p}$.
It is convenient to combine \eqref{eqn:cell_equation} and \eqref{eqn:facet_equation} into a system of equations

\begin{equation}
  \begin{pmatrix}
    A_{c\mapsfrom c} & 0                & A_{c\mapsfrom f} \\
    A_{f\mapsfrom c} & -id              & 0                \\
    0                & A_{f\mapsfrom f} & -id
  \end{pmatrix}
  \begin{pmatrix}
    \vec{u}^{(c)}   \\
    \vec{u}^{(\pm)} \\
    \vec{w}^{(f)}
  \end{pmatrix}
  =
  \begin{pmatrix}
    \vec{b}^{(c)} \\
    0             \\ 0
  \end{pmatrix}.
  \label{eqn:global_matrix_representation-IP}
\end{equation}

\noindent
In analogy to \eqref{eqn:matrix_form_IP}, the facet dof-vectors $\vec{w}^{(f)}$, $\vec{u}^{(\pm)}$ can be partitioned as $\vec{w}^{(f)}= (\vec{w}^{(f)}|_{F_1}, \vec{w}^{(f)}|_{F_2}, \ldots)$ and $\vec{u}^{(\pm)}= (\vec{u}^{(\pm)}|_{F_1}, \vec{u}^{(\pm)}|_{F_2}, \ldots)$ where the vectors $\vec{w}^{(f)}|_F$ and $\vec{u}^{(\pm)}|_F$ contain the unknowns on a single facet. Similarly, the matrices $A_{c\gets c}$ (which is \emph{not} identical to matrix $A$ in \eqref{eqn:matrix_form_IP}), $A_{c\gets f}$, $A_{f\gets c}$ and $A_{f\gets f}$ decompose into blocks that encode the coupling between two mesh entities (cells or facets). For example, $A_{c\gets f}|_{K\gets F}$ describes how the unknowns on facet $F$ couple to the unknowns in cell $K$. 


The block-Jacobi iteration in \eqref{eqn:block-Jacobi-IP} can be re-written schematically as in \cref{alg:interior-penalty-single-level-multiple-sweeps-and-auxiliary-variables}, which is exactly equivalent to
\cref{alg:vanilla-interior-penalty-single-level}. However, the residual $\vec{r}^{(c)}=\vec{b}^{(c)}-A\vec{u}^{(c)}$ is computed in three stages by using \eqref{eqn:global_matrix_representation-IP}:
\begin{algorithm}[htb]
  \caption{Matrix-free block-Jacobi iteration using auxiliary facet-variables. Input: initial solution $\vec{u}^{(c)}_0$, right hand side $\vec{b}^{(c)}$, relaxation parameter $\omega$, number of iterations $n_{\text{iter}}$. Output: updated solution $\vec{u}^{(c)}$ after $n_{\text{iter}}$ iterations. 
  }
  \label{alg:interior-penalty-single-level-multiple-sweeps-and-auxiliary-variables}
  \begin{algorithmic}[1]
    \State{Set $\vec{u}^{(c)}\gets\vec{u}^{(c)}_0$}
    \For{$k=1,2,\dots,n_{\text{iter}}$}
      \For{every cell $K \in \Omega_h$}
        \label{alg:interior-penalty-single-level-multiple-sweeps-and-auxiliary-variables:first-loop-block:start}
        \For{every facet $F\in \mathcal{F}(K)$ of cell $K$}
          \State {Assemble $A_{f \gets c}|_{F\gets K}$}
          \State {Set $\vec{u}^{(\pm)}|_F = A_{f \gets c}|_{F\gets K} \vec{u}^{(c)}|_K$}
            \Comment {Project solution onto facets}\label{alg:interior-penalty-single-level-multiple-sweeps-and-auxiliary-variables:projection}
        \EndFor
      \EndFor\label{alg:interior-penalty-single-level-multiple-sweeps-and-auxiliary-variables:first-loop-block:end}
      \State {Exchange $\vec{u}^{(\pm)}$ between non-overlapping subdomains.}
    \For{every facet $F\in \mathcal{E}_h$}
      \label{alg:interior-penalty-single-level-multiple-sweeps-and-auxiliary-variables:second-loop-block:start}
      \State {Assemble $A_{f \gets f}|_{F\gets F}$}
      \State {Set $\vec{w}^{(f)}|_F = A_{f \gets f}|_{F\gets F} \vec{u}^{(\pm)}|_F$}
        \Comment{Compute numerical fluxes}
    \EndFor
      \label{alg:interior-penalty-single-level-multiple-sweeps-and-auxiliary-variables:second-loop-block:end}
    \For{every cell $K \in \Omega_h$}\label{alg:interior-penalty-single-level-multiple-sweeps-and-auxiliary-variables:third-loop-block:begin}
    \State {Assemble $A_{c \gets c}|_{K\gets K}$}
    \Comment {On-the-fly assembly of cell-local matrix}
    \State {$\vec{r}^{(c)}|_K \gets \vec{b}^{(c)}|_K -  A_{c\gets c}|_{K\gets K} \vec{u}^{(c)}|_K$}
    \Comment {cell-local contribution to residual}\label{alg:interior-penalty-single-level-multiple-sweeps-and-auxiliary-variables:cell_residual}
    \For{every facet $F\in\mathcal{F}(K)$ of cell $K$}
    \State {Assemble local matrix $A_{c\gets f}|_{K\leftarrow F}$}
    \State {$\vec{r}^{(c)}|_K \gets \vec{r}^{(c)}|_K - A_{c\gets f}|_{K\gets F} \vec{w}^{(f)}|_{F}$}
    \Comment {contribution from facets to residual}\label{alg:interior-penalty-single-level-multiple-sweeps-and-auxiliary-variables:facet_residual}
    \EndFor
    \State {$\vec{u}^{(c)}|_K \gets \vec{u}^{(c)}|_K + \omega  A_{K\gets K}^{-1} \vec{r}^{(c)}|_K$}
      \Comment{update state}
      \label{alg:interior-penalty-single-level-multiple-sweeps-and-auxiliary-variables:update-state}
    \EndFor\label{alg:interior-penalty-single-level-multiple-sweeps-and-auxiliary-variables:iteration:end}
    \EndFor
    \State{\Return $\vec{u}^{(c)}$}
  \end{algorithmic}
\end{algorithm}

\begin{enumerate}
  \item Project cell data $\vec{u}^{(c)}$ onto the faces, $\vec{u}^{(\pm)} = A_{f\gets c}\vec{u}^{(c)}$. For this, compute $\vec{u}^{(\pm)}|_F$ from $\vec{u}^{(c)}|_{K^\pm}$ on each facet $F$ using \eqref{eqn:facet_equation}. Only the cell values from adjacent cells $K^+$, $K^-$ are required on facet $F$. This is the second block line from \eqref{eqn:global_matrix_representation-IP}.
  \item Evaluate the numerical flux $\vec{w}^{(f)} = A_{f\gets f}\vec{u}^{(\pm)}$ based on the projections. This is achieved by computing $\vec{w}^{(f)}|_F$ from $\vec{u}^{(\pm)}|_F$ on each facet $F$ according to last block line from \eqref{eqn:global_matrix_representation-IP}.
  \item Compute the residual $\vec{r}^{(c)} = \vec{b}^{(c)}-A_{c\gets c}\vec{u}^{(c)}-A_{c\gets f}\vec{w}^{(f)}$ by subtracting the expressions on the left-hand side of the first block row of \eqref{eqn:global_matrix_representation-IP} from the right-hand side.
    To compute the residual $\vec{r}^{(c)}|_K$ in cell $K$, only the local value $\vec{u}^{(c)}|_K$ and the numerical flux $\vec{w}^{(f)}|_{F}$ on all facets $F\in\mathcal{F}(K)$ touching the cell $K$ are required.
\end{enumerate}
\noindent
The update $\vec{u}^{(c)}|_K = \vec{u}^{(c)}|_K +\omega A_{K\gets K}^{-1} \vec{r}^{(c)}|_K$ of the local solution is the same as in \cref{alg:vanilla-interior-penalty-single-level}.

In contrast to \cref{alg:vanilla-interior-penalty-single-level}, a single iteration in \cref{alg:interior-penalty-single-level-multiple-sweeps-and-auxiliary-variables} no longer requires the evaluation of terms that directly couple the solution in neighbouring cells. Instead, the operator application is split into volumetric terms and facet contributions. The latter enter the update scheme through the loop over the adjacent facets of a cell. Looking at \cref{alg:interior-penalty-single-level-multiple-sweeps-and-auxiliary-variables} we also find the following:

\begin{observation}
 Due to the introduction of the helper variables $\vec{u}^{(\pm)}$ and $\vec{w}^{(f)}$ on the facets according to \cref{technique:left-right-face-projections} and \cref{technique:numerical-flux-outcome}, we 
 \begin{itemize}
   \item can store the solution along the faces continuously in memory. Hence, the numerical flux calculations result in coalesced memory accesses and reading data from neighbouring cells do not induce scattered memory accesses for Gauss-Lobatto basis functions (\cref{def:vectorisation});
   \item do not have to maintain the backup of a previous iteration (\cref{def:memory_overhead});
   \item no longer have to access the cell solution from the previous iteration to obtain consistent fluxes. This avoids repeated reads of volumetric data (\cref{def:weak_single_touch}).
 \end{itemize}
%
\end{observation}
\begin{observation}
In a distributed-memory parallel setting based on domain decomposition, another advantage of \cref{alg:interior-penalty-single-level-multiple-sweeps-and-auxiliary-variables} is that only data on facets has to be exchanged between neighbouring processors. At higher discretisation order $p>1$ this reduces the communicated data volume by a factor $p+1$, which can lead to notable improvements in parallel scalability. 
\end{observation}

\paragraph{Schur-complement}\label{sec:schur_complement}

The equivalence of the original weak form of the interior penalty formulation in \eqref{eqn:weak_form_IP} on the one hand and \eqref{eqn:weak_form_IP_facet_spaces}, \eqref{eqn:solution_projection}, \eqref{eqn:gradient_projection}, \eqref{eqn:numerical_fluxes} on the other hand can also be expressed on an algebraic level. To see this, eliminate $\vec{w}^{(f)}$ from \eqref{eqn:cell_equation} with the help of \eqref{eqn:facet_equation} to obtain
\[
  \underbrace{\left(A_{c\gets c} + A_{c\gets f} A_{f\gets f} A_{f\gets c}\right)}_{=:S}\vec{u}^{(c)} = \vec{b}^{(c)}.
\]

\noindent
The matrix $S:=A_{c\gets c} + A_{c\gets f} A_{f\gets f} A_{f\gets c}$, which is identical to $A$ in \eqref{eqn:matrix_form_IP}, is the \emph{Schur-complement} of the matrix in \eqref{eqn:global_matrix_representation-IP}
and its block-diagonal determines the iteration matrix in line~\ref{alg:interior-penalty-single-level-multiple-sweeps-and-auxiliary-variables:update-state} of \cref{alg:interior-penalty-single-level-multiple-sweeps-and-auxiliary-variables}.
%
Yet, there is no need to assemble $S$ globally since only its block-diagonal is required. In each cell $K$ the block-diagonal of $A$ is readily constructed as

\[
  A_{K\gets K} = S_{K\gets K} = A_{c\gets c}|_{K\gets K} + \sum_{F\in\mathcal{F}(K)}
  A_{c\gets f}|_{K\gets F} A_{f\gets f}|_{F\gets F} A_{f\gets c}|_{F\gets K}.
\]

\noindent
For homogeneous, isotropic problems on a uniform mesh the matrices $A_{K\gets K}$ will be identical in each cell and $A_{K\gets K}^{-1}$ can be precomputed and stored once at the beginning of the linear solve. However, if the mesh has been refined adaptively to obtain cells $K$ of varying size $h_K$, the matrix $A_{K\leftarrow K} = \sum_{\alpha} h_K^{q_\alpha} B_\alpha$ is linear combination of $h_K$-independent building blocks $B_\alpha$, which can be precomputed once on a reference element. As the facet- terms and cell- terms in \eqref{eqn:block-Jacobi-IP} scale with different powers $q_\alpha$ of the mesh size, the inverse of $A_{K\gets K}$ will differ from cell to cell and needs to be computed on-the-fly as follows:
\begin{technique}
  \label{observation:single-level:precompute-block-inverse}
  For homogeneous, isotropic problems on an adaptively refined mesh, it is only necessary to precompute and store two small reference matrices --- one representing volumetric terms, one hosting face terms --- to construct the cell-local matrix $A_{K \gets K}$ in each cell $K$ and to invert it on-the-fly.  
\end{technique}
\noindent
For problems with inhomogeneous, non-isotropic coefficients, the construction of $A_{K\gets K}$ from pre-computed building blocks according to \cref{observation:single-level:precompute-block-inverse} is not possible. However, it is still possible to maintain a matrix-free implementation by re-assembling $A |_{K\leftarrow K}$ in each cell. We therefore conclude:
\noindent
\begin{observation}
 The implementation in \cref{alg:interior-penalty-single-level-multiple-sweeps-and-auxiliary-variables} is strictly matrix-free.
\end{observation}
\noindent
It should also be pointed out that:
\begin{observation}
  \label{obs:cost_of_local_matrix_inversion}
For high polynomial degrees $p$ the cost of inverting the cell-local matrices $A_{K\gets K}$ if usually a factor $p$ times more expensive than the assembly of $A_{K\gets K}$ itself. Hence, the cost of this inversion is a good measure for the computational overhead that arises when going from a homogeneous, isotropic problem on a uniformly refined mesh to the more complicated setup of an adaptively refined mesh and/or inhomogeneous, non-isotropic problems.
\end{observation}

A closer inspection of $A_{K\gets K}$ reveals that even for homogeneous, isotropic problems on uniformly refined meshes the matrix $A_{K\gets K}$ differs for cells in the interior and adjacent to the boundary of the domain. The same applies for the small reference matrices that are used in \cref{observation:single-level:precompute-block-inverse}. However, ignoring this fact and using a single $A_{K\gets K}^{-1}$ (for homogeneous, isotropic problems) or one set of building blocks $B_\alpha$ (for inhomogeneous, non-isotropic problems) still results in an iterative scheme which converges to the correct solution. This is because the block-Jacobi method in \eqref{eqn:block-Jacobi-IP} can be interpreted as a Richardson iteration preconditioned with the inverse of the block-diagonal of $A$ and modifying the preconditioner will not change the fixed point of the iteration.
\begin{technique}
  \label{technique:approximate-boundary-operators}
  We neglect the fact that the block-inverse $A_{K\gets K}^{-1}$ differs between cells in the interior and adjacent to the boundary of the domain. In our implementation, we only use the block-inverse derived for cells in the interior of the domain.
\end{technique}
\noindent
As our numerical experiments in \cref{section:results} demonstrate, the simplification in \cref{technique:approximate-boundary-operators} still results in rapidly converging multigrid method.
\paragraph{Structure of block-matrices}
The entries of the small block-matrices $A_{c\gets c}|_{K\gets K}$, $A_{f\gets c}|_{F\gets K}$, $A_{c\gets f}|_{K\gets F}$ and $A_{f\gets f}|_{F\gets F}$ in \eqref{eqn:global_matrix_representation-IP} depend on the choice of basis functions for the spaces $\VDG{h}{p}$ and $\FDG{h}{p}$. In particular, if a Gauss-Lobatto basis is chosen for $\VDG{h}{p}$, the $(p+1)^{d-1}\times (p+1)^{d}$ matrices $A_{f\gets c}|_{F\gets K}$ contain only $(p+1)^{2(d-1)}$ non-zero entries since only $(p+1)^{d-1}$ of the cell-wise basis functions are non-zero on the cell boundaries. However, mass matrices are not diagonal in this basis, since $p+2$ Gauss-Lobatto quadrature points are required to integrate functions of degree $2p$ exactly; this is the reason why the number of non-zero entries in each $A_{f\gets c}|_{F\gets K}$ is not $(p+1)^{d-1}$, even if a Gauss-Lobatto basis is also used for $\FDG{h}{p}$. In contrast, Gauss-Legendre basis functions lead to diagonal mass matrices but result in dense $A_{f\gets c}|_{F\gets K}$. These observations have the following implications for the implementation:
\begin{observation}
  \label{observation:GaussLobatto-sparsity}
  For Gauss-Lobatto basis functions, the operation $\vec{v}^{(f)} = A_{f \gets c}\vec{w}^{(c)}$ simply extracts unknowns from the underlying cell data representation of $\vec{w}^{(c)}$ and copies them to $\vec{v}^{(f)}$. In terms of data structures, this corresponds to a strided access to a subarray of an array. In other words, when iterating over the cells of the mesh, the choice of Gauss-Lobatto basis results in a scatter of the data in $\vec{w}^{(c)}$ that is associated with the surface of each cell.
  In contrast, for Gauss-Legendre basis functions multiplication with $A_{f \gets c}$  corresponds to dense matrix-vector products $\vec{v}^{(f)}|_F = A_{f \gets c}|_{F\gets K}\vec{w}^{(c)}|_K \vec{w}^{(f)}|_K$ in each cell-facet pair $(K,F)$.
\end{observation}
\noindent
Both choices of basis functions are popular in different application areas, and we explore them numerically in \cref{sec:numerical_GL_vs_GLL}.

\subsubsection{DG as single level solver}

\begin{algorithm}[htb]
  \caption{
  Wrapper for Matrix-free block-Jacobi iteration with dynamic exit criterion. Input: initial solution $\vec{u}^{(c)}_0$, right hand side $\vec{b}^{(c)}$, relaxation parameter $\omega$, maximal number of iterations $n_{\text{iter}}$ and tolerance $\epsilon$. Output: updated solution $\vec{u}^{(c)}$ after $n_{\text{iter}}$ iterations or convergence to tolerance $\epsilon$.}
  \label{alg:single-level-termination-criterion}
  \begin{algorithmic}[1]
    \State{Set $\vec{u}^{(c)}\gets\vec{u}^{(c)}_0$}\label{alg:single-level-termination-criterion:backup}
    \For{$k=1,2,\dots,n_{\text{iter}}$}
      \For{every cell $K \in \Omega _h$}
        \State{Set $\vec{u}^{(c)}_{\text{old}}\gets\vec{u}^{(c)}$}
          \Comment{Store previous iterate to check exit condition}
        \State [\ldots]
        \State{Compute $\rho_k = \vert|\vec{u}^{(c)}_{\text{old}}-\vec{u}^{(c)}\vert|$}\Comment{preconditioned residual norm}
      \EndFor
        \If{$\rho_k/\rho_1<\epsilon$}
          \State{\textbf{exit loop}}\Comment{check convergence}
        \EndIf      
    \EndFor
  \end{algorithmic}
\end{algorithm}

If the DG scheme is used as a standalone solver rather than a smoother, we might want to supplement the implementation with an early termination criterion to stop the iteration once $\vec{u}^{(c)}$ is sufficiently close to the true solution. 
This is shown in \cref{alg:single-level-termination-criterion}, where ``[\dots]'' stands for all operations in lines \ref{alg:vanilla-interior-penalty-single-level:iteration:start} to \ref{alg:vanilla-interior-penalty-single-level:iteration:end} in \cref{alg:vanilla-interior-penalty-single-level} and lines \ref{alg:interior-penalty-single-level-multiple-sweeps-and-auxiliary-variables:first-loop-block:start} to \ref{alg:interior-penalty-single-level-multiple-sweeps-and-auxiliary-variables:iteration:end} in \cref{alg:interior-penalty-single-level-multiple-sweeps-and-auxiliary-variables} respectively.
Observe that we backup the solution at the previous iteration and store it in $\vec{u}^{(c)}_{\text{old}}$, since the change in solution from one iteration to the next is used to assess convergence (see also the discussion of the preconditioned residual in \cref{section:multigrid}).
Consequently, such a dynamic termination criterion introduces additional volumetric reads and writes.

\subsubsection{Reduction in memory traffic}

Let us assume that  
all local operator matrices are cached. 
Furthermore, the residual data $\vec{r}^{(c)}|_K$ is a local variable that does not have to be held persistently and can therefore also be stored in cache.
In this case only data for the solution and right-hand side vectors needs to be read from and written to memory. A simple performance model reveals that the introduction of additional auxiliary variables with \cref{technique:left-right-face-projections} and \cref{technique:numerical-flux-outcome} reduces the memory traffic for higher polynomial degrees $p$, while the memory traffic is increased for smaller $p$. Since volumetric fields require the storage of $(p+1)^d$ unknowns per cell, \cref{alg:vanilla-interior-penalty-single-level} requires $(2d+5)(p+1)^d$ memory accesses per cell:
\begin{itemize}
  \item $\vec{u}^{(c)}$ is read and $\vec{u}^{(c)}_{\text{old}}$ is written in line \ref{alg:vanilla-interior-penalty-single-level:back-up-solution} $\Rightarrow$ $2(p+1)^d$ memory accesses per cell 
  \item $\vec{b}^{(c)}|_K$ and $\vec{u}^{(c)}_{\text{old}}|_K$ are read in line \ref{alg:vanilla-interior-penalty-single-level:residual} $\Rightarrow$ $2(p+1)^d$ memory accesses per cell
  \item $\vec{u}^{(c)}_{\text{old}}|_{K'}$ is read for all $2d$ face-connected neighbours $K'$ in line \ref{alg:vanilla-interior-penalty-single-level:loop-over-neighbours} $\Rightarrow 2d(p+1)^d$ memory accesses per cell
  \item $\vec{u}^{(c)}|_K$ is written in line \ref{alg:vanilla-interior-penalty-single-level:solution_update} $\Rightarrow$ $(p+1)^d$ memory accesses per cell
\end{itemize}

\noindent
In contrast, \cref{alg:interior-penalty-single-level-multiple-sweeps-and-auxiliary-variables} requires $3(p+1)^d + 7d(p+1)^{d-1}$ memory accesses per cell since fields stored on the facets require only $(p+1)^{d-1}$ unknowns per facet. Bearing in mind that the number of facets is $d$ times larger than the number of cells, this result is obtained with the following counting:
\begin{itemize}
  \item $\vec{u}^{(+)}|_F$ and $\vec{u}^{(-)}|_F$ are read and $\vec{w}^{(f)}|_F$ is written for each facet $F$ in line \ref{alg:interior-penalty-single-level-multiple-sweeps-and-auxiliary-variables:projection} $\Rightarrow$ $3(p+1)^{d-1}$ memory accesses per facet
  \item $\vec{u}^{(c)}|_K$ and $\vec{b}^{(c)}|_K$ are read in each cell $K$ in line \ref{alg:interior-penalty-single-level-multiple-sweeps-and-auxiliary-variables:cell_residual} $\Rightarrow$ $2(p+1)^d$ memory accesses per cell
  \item $\vec{w}^{(f)}|_F$ is read for all $2d$ facets of each cell in line \ref{alg:interior-penalty-single-level-multiple-sweeps-and-auxiliary-variables:facet_residual} $\Rightarrow$ $2d(p+1)^{d-1}$ memory accesses per cell
  \item $\vec{u}^{(c)}|_K$ is written back in line \ref{alg:interior-penalty-single-level-multiple-sweeps-and-auxiliary-variables:update-state} $\Rightarrow$ $(p+1)^d$ memory accesses per cell
\end{itemize}
As a consequence, \cref{alg:interior-penalty-single-level-multiple-sweeps-and-auxiliary-variables} requires fewer memory accesses than \cref{alg:vanilla-interior-penalty-single-level} if $p\ge 2$. The number of memory accesses for different polynomial degrees $p$ and the reduction that results from using \cref{alg:interior-penalty-single-level-multiple-sweeps-and-auxiliary-variables} instead of \cref{alg:vanilla-interior-penalty-single-level} is shown in \cref{tab:memory_references}. 

If \cref{alg:interior-penalty-single-level-multiple-sweeps-and-auxiliary-variables} is used as a standalone solver with a dynamic termination criterion (\cref{alg:single-level-termination-criterion}), the number of memory references increases to $5(p+1)^d + 7d(p+1)^{d-1}$ since $\vec{u}^{(c)}$ is read and $\vec{u}^{(c)}_{\text{old}}$ is written in line \ref{alg:single-level-termination-criterion:backup} of \cref{alg:single-level-termination-criterion}. 

In the limit $p\rightarrow \infty$ the reduction factor is $3.0\times$ ($1.8\times$) in $d=2$ dimensions and $3.7\times$ ($2.2\times$) in $d=3$ dimensions, where the number in brackets are obtained if the \cref{alg:interior-penalty-single-level-multiple-sweeps-and-auxiliary-variables} is used as a standaline solver. As can be seen in \cref{tab:memory_references}, for some lower polynomial degrees \cref{alg:interior-penalty-single-level-multiple-sweeps-and-auxiliary-variables} will require more memory accesses than \cref{alg:vanilla-interior-penalty-single-level}.

\begin{table}[htb]
\caption{
  Number of memory accesses per cell for \cref{alg:vanilla-interior-penalty-single-level} and \cref{alg:interior-penalty-single-level-multiple-sweeps-and-auxiliary-variables} in $d=2$ (top) and $d=3$ (bottom) dimensions.
  \cref{alg:interior-penalty-single-level-multiple-sweeps-and-auxiliary-variables} marked with $\dagger$ presents data if it is used as a standalone solver, i.e.~subject to a dynamic termination criterion.
  The line below the access count shows the arising reduction in memory accesses, being 1 for the baseline.
  \label{tab:memory_references}
}
\footnotesize
\begin{tabular}{cccccccccccc}
\hline
 & degree $p$ & 1 & 2 & 3 & 4 & 5 & 6 & 7 & 8 & 9 & 10\\\hline
\multirow{5}{6ex}{$d = 2$} & \cref{alg:vanilla-interior-penalty-single-level} & $36$ & $81$ & $144$ & $225$ & $324$ & $441$ & $576$ & $729$ & $900$ & $1089$\\
& \cref{alg:interior-penalty-single-level-multiple-sweeps-and-auxiliary-variables} & $40$ & $69$ & $104$ & $145$ & $192$ & $245$ & $304$ & $369$ & $440$ & $517$\\
& reduction  & $  0.90\times$ & $  1.17\times$ & $  1.38\times$ & $  1.55\times$ & $  1.69\times$ & $  1.80\times$ & $  1.89\times$ & $  1.98\times$ & $  2.05\times$ & $  2.11\times$\\
& \cref{alg:interior-penalty-single-level-multiple-sweeps-and-auxiliary-variables}${}^{\dagger}$ & $48$ & $87$ & $136$ & $195$ & $264$ & $343$ & $432$ & $531$ & $640$ & $759$\\
& reduction${}^{\dagger}$  & $  0.75\times$ & $  0.93\times$ & $  1.06\times$ & $  1.15\times$ & $  1.23\times$ & $  1.29\times$ & $  1.33\times$ & $  1.37\times$ & $  1.41\times$ & $  1.43\times$\\
\hline
\multirow{5}{6ex}{$d = 3$} & \cref{alg:vanilla-interior-penalty-single-level} & $88$ & $297$ & $704$ & $1375$ & $2376$ & $3773$ & $5632$ & $8019$ & $11000$ & $14641$\\
& \cref{alg:interior-penalty-single-level-multiple-sweeps-and-auxiliary-variables} & $108$ & $270$ & $528$ & $900$ & $1404$ & $2058$ & $2880$ & $3888$ & $5100$ & $6534$\\
& reduction  & $  0.81\times$ & $  1.10\times$ & $  1.33\times$ & $  1.53\times$ & $  1.69\times$ & $  1.83\times$ & $  1.96\times$ & $  2.06\times$ & $  2.16\times$ & $  2.24\times$\\
& \cref{alg:interior-penalty-single-level-multiple-sweeps-and-auxiliary-variables}${}^{\dagger}$ & $124$ & $324$ & $656$ & $1150$ & $1836$ & $2744$ & $3904$ & $5346$ & $7100$ & $9196$\\
& reduction${}^{\dagger}$  & $  0.71\times$ & $  0.92\times$ & $  1.07\times$ & $  1.20\times$ & $  1.29\times$ & $  1.38\times$ & $  1.44\times$ & $  1.50\times$ & $  1.55\times$ & $  1.59\times$\\
\hline
\end{tabular}
\end{table}

\subsubsection{Simplifaction and extension of numerical fluxes}
For the interior penalty discretisation in \eqref{eqn:weak_form_IP} the flux $\vec{w}^{(f)} = B^+ \vec{u}^{(+)}+B^-\vec{u}^{(-)}$ is a linear combination of $\vec{u}^{(+)}$, $\vec{u}^{(-)}$ for some matrices $B^+$, $B^-$. 
In principle, it would therefore be possible to reduce the storage requirements further by not storing the projections $\vec{u}^{(\pm)}$ at all and directly accumulating into $\vec{w}^{(f)}$. 
We do not exploit this insight here and store the two projections $\vec{u}^{(\pm)}$ in \cref{alg:interior-penalty-single-level-multiple-sweeps-and-auxiliary-variables} explicitly.

For different choices of the numerical flux, other variables than the solution and its normal derivative might have to be stored in $\vec{u}^{(\pm)}$, for example one might also want to project derivatives of higher order. Along the same lines, PDEs that contain non-conservative terms might result in double-valued fluxes: in this case the flux depends not only on the facet $F$ but also on the cell $K\in \{K^+,K^-\}$ from which it is accessed. It will then be necessary to store two vectors $\vec{w}^{(f,+)}$ and $\vec{w}^{(f,-)}$, which doubles the memory footprint relative to the single-valued $\vec{w}^{(f)}$.

\subsubsection{Domain decomposition}
\label{subsubsection:single-level:domain-decomposition}

The algorithmic footprint accommodates a non-overlapping domain decomposition that can be mapped onto a shared or distributed memory systen.
Let the computational domain be subdivided into non-overlapping subdomains where each cell is assigned to a unique processor.
Our code employs the Peano space-filling curve (SFC) to determine this assignment, i.e.~the cells are enumerated along the SFC and this sequence of cells is then subdivided into contiguous subsequences, each of which is assigned to one processor.
This results in connected subpartitions with excellent surface-to-volume ratios \cite{weinzierl2019peano}.
Each individual processor projects the solution $\vec{u}^{(c)}$ onto the facets of all cells that it owns. 
Facets are hence held redundantly for facets along subdomain boundaries.
After the projection phase, the individual $\vec{u}^{(\pm)}|_F$ need to be exchanged between the processors for the subdomain interface facets, to ensure that both have the projections $\vec{u}^{(+)}|_F$ and $\vec{u}^{(-)}|_F$ readily available,
i.e.~we compute $\vec{w}^{(f)}|_F$ redundantly.

\begin{observation}
Compared to cell data, variables stored on the facets require $p+1$ times less storage than adjacent cells.
As we exchange facet data instaed of cell data, the exchanged data volume is reduced by the same factor.
\cref{technique:left-right-face-projections} reduces the communication overhead.
\end{observation}

\subsection{A single-touch grid traversal implementation}
\label{sec:single-touch-grid-traversal} 
In contrast to \cref{alg:vanilla-interior-penalty-single-level}, the improved implementation in \cref{alg:interior-penalty-single-level-multiple-sweeps-and-auxiliary-variables} computes the residual $\vec{r}^{(c)}$ in three stages: it projects the solution $\vec{u}^{(c)}$ onto the faces to obtain $\vec{u}^{(\pm)}$, computes the numerical flux $\vec{w}^{(f)}$, and eventually constructs the residual as sum of facet- and cell contributions.
Schematically, the outer block-Jacobi iteration with loop index $k$ in \cref{alg:interior-penalty-single-level-multiple-sweeps-and-auxiliary-variables} can be written as a repeated application of the three stages:

\[
  \UOLoverbrace{1\rightarrow 2\rightarrow 3}^{k=1}\rightarrow\UOLoverbrace{1\rightarrow 2\rightarrow 3}^{k=2}\rightarrow\UOLoverbrace{1\rightarrow 2\rightarrow 3}^{k=3}\rightarrow 1\rightarrow 2\rightarrow\dots
\]

\noindent
This sequence of operations is not (weakly) single touch, as the cell data is read at least twice. If realised through parallel loops, this also induces three synchronisation points which can limit parallel scalability.
In the following we outline a strategy for overcoming these issues.

\paragraph{Loop fusion into a cell-wise realisation}
\label{sec:single-level:multiple-sweep:fusion-of-face-and-cell-loop}

Our realisation with auxiliary variables in \cref{alg:interior-penalty-single-level-multiple-sweeps-and-auxiliary-variables} separates the smoothing step into a sequence of three mesh traversals.
This sequence can be collapsed into two consecutive loops each of which exclusively runs over the mesh cells: for this we combine the two loops in lines \ref{alg:interior-penalty-single-level-multiple-sweeps-and-auxiliary-variables:second-loop-block:start}-\ref{alg:interior-penalty-single-level-multiple-sweeps-and-auxiliary-variables:second-loop-block:end} and \ref{alg:interior-penalty-single-level-multiple-sweeps-and-auxiliary-variables:third-loop-block:begin}-\ref{alg:interior-penalty-single-level-multiple-sweeps-and-auxiliary-variables:iteration:end} of \cref{alg:interior-penalty-single-level-multiple-sweeps-and-auxiliary-variables} into a single loop over cells, which in turn for each cell $K$ contains an inner loop over the facets $F\in\mathcal{F}(K)$. To achieve this, we introduce a flag $\textsf{touched}(F)$ on each facet which indicates whether $\vec{w}^{(f)}|_F$ has already been computed and can therefore be used for updating the residual $\vec{r}^{(c)}|_K \gets \vec{r}^{(c)}|_K - A_{c\gets f}|_{K\gets F} \vec{w}^{(f)}|_{F}$. At the beginning of the block-Jacobi iteration, this flag is set to $\textsf{false}$ for all facets of the mesh. As we loop over the cells, each cell checks for each adjacent facet $F\in \mathcal{F}(K)$ whether the facet has been touched yet. 
If this is not the case, i.e. if $\textsf{touched}(F)= \textsf{false}$, we compute $\vec{w}^{(f)}|_F = A_{f\gets f}|_{F\gets F}\vec{u}^{(\pm)}|_F$ on this facet and afterwards set $\textsf{touched}(F)\gets\textsf{true}$. This results in an implemenation of \cref{alg:interior-penalty-single-level-multiple-sweeps-and-auxiliary-variables} with two cell loops only, each of which corresponds to a strict cell-wise mesh traversal \cite{weinzierl2019peano}. It should be stressed that we do not implement this variation of \cref{alg:interior-penalty-single-level-multiple-sweeps-and-auxiliary-variables} in our code, but discuss it here to motivate the additional code transformations described in the next paragraph. Observe also that many mesh traversal codes do not require explicit storage of the flag $\textsf{touched}(F)$ since this information is stored implicitly in the ordering of the mesh entities: when processing a given facet, it is possible to infer from the index of the facet whether it has been visited previouly. In a parallel implementation some synchronisation is required to consistently update $\vec{w}^{(f)}|_F$. This can be avoided by computing $\vec{w}^{(f)}|_F$ redundantly on each processor for facets $F$ on subdomain boundaries.

\paragraph{Loop fusion and shifting}

Having combied two of the three mesh-traversals in \cref{alg:interior-penalty-single-level-multiple-sweeps-and-auxiliary-variables} as described in \cref{sec:single-level:multiple-sweep:fusion-of-face-and-cell-loop}, we finally fuse all operations into a single loop by combining 
the two loops in lines \ref{alg:interior-penalty-single-level-multiple-sweeps-and-auxiliary-variables:second-loop-block:start}-\ref{alg:interior-penalty-single-level-multiple-sweeps-and-auxiliary-variables:second-loop-block:end} and \ref{alg:interior-penalty-single-level-multiple-sweeps-and-auxiliary-variables:third-loop-block:begin}-\ref{alg:interior-penalty-single-level-multiple-sweeps-and-auxiliary-variables:iteration:end} for the \textit{current} block-Jacobi iteration $k$ with the loop over cells in lines \ref{alg:interior-penalty-single-level-multiple-sweeps-and-auxiliary-variables:first-loop-block:start}-\ref{alg:interior-penalty-single-level-multiple-sweeps-and-auxiliary-variables:first-loop-block:end} in the \textit{next} iteration $k+1$:
\begin{technique}
  \label{technique:shift-and-fuse}
  Each block-Jacobi iteration can be written as a single mesh traversal provided we \emph{shift} the operator evaluation: the fields are updated in the order $\vec{u}^{(\pm)}\xrightarrow{2}\vec{w}^{(f)} \xrightarrow{3} \vec{u}^{(c)} \xrightarrow{1} \vec{u}^{(\pm)}$ instead of $\vec{u}^{(c)}\xrightarrow{1} \vec{u}^{(\pm)}\xrightarrow{2} \vec{w}^{(f)}\xrightarrow{3} \vec{u}^{(c)} $.
\end{technique}
This is illustrated in the following diagram:
\[
  \UOLoverbrace{1\rightarrow}[2\rightarrow 3]^{k=1} 
  \UOLunderbrace{\rightarrow}[1]_{\text{fused}} 
  \UOLoverbrace{\rightarrow}[2\rightarrow 3]^{k=2} 
  \UOLunderbrace{\rightarrow}[1]_{\text{fused}} 
  \UOLoverbrace{\rightarrow}[2\rightarrow 3]^{k=3}
  \UOLunderbrace{\rightarrow 1}_{\text{fused}} 
  \rightarrow2\rightarrow\dots
\]
\noindent

\begin{algorithm}
  \caption{Matrix-free block-Jacobi iteration using auxilliary facet-variables and loop fusion. $\text{BlockJacobi}(\vec{u}^{(c)}_0, \vec{b}^{(c)};\omega,n_{\text{iter}},\epsilon)$ Input: initial solution $\vec{u}^{(c)}_0$, right hand side $\vec{b}^{(c)}$, relaxation parameter $\omega$, number of iterations $n_{\text{iter}}$, tolerance $\epsilon$. Output: solution $\vec{u}^{(c)}$ and its facet-projection $\vec{u}^{(\pm)}=A_{f\gets c}\vec{u}^{(c)}$ after $n_{\text{iter}}$ iterations.}
  \label{alg:interior-penalty-single-level-one-mesh-sweep}
  \begin{algorithmic}[1]
    \State{Set $\vec{u}^{(c)}\gets \vec{u}^{(c)}_0$}
    \For{every cell $K \in \Omega _h$}
      \label{alg:interior-penalty-single-level-one-mesh-sweep:projection-preamble-start}
      \For{every facet $F\in \mathcal{F}(K)$ of cell $K$}
        \State {Assemble $A_{f \gets c}|_{K\gets F}$}
        \State {Set $\vec{u}^{(\pm)}|_F = A_{f \gets c}|_{F\gets K} \vec{u}^{(c)}|_K$}
          \Comment {Project solution onto facets}
      \EndFor
    \EndFor
      \label{alg:interior-penalty-single-level-one-mesh-sweep:projection-preamble-end}
    \For{$k=1,2,\dots,n_{\text{iter}}$}
    \State {Exchange $\vec{u}^{(\pm)}$ between non-overlapping subdomains.}
    \State {Set $\textsf{touched}(F) = \textsf{false}$ for all $F\in\mathcal{E}_h$}
    \Comment{Mark all facets as untouched}
    \For{every cell $K \in \Omega _h$}\label{alg:interior-penalty-single-level-one-mesh-sweep:fused_loop:start}
    \State {Assemble $A_{c \gets c}|_{K\gets K}$}
    \Comment {On-the-fly assembly of cell-local matrix}
    \State {$\vec{r}^{(c)}|_K \gets \vec{b}^{(c)}|_K -  A_{c\gets c}|_{K\gets K} \vec{u}^{(c)}|_K$}
    \Comment {cell-local contribution to residual}
    \For{every facet $F\in\mathcal{F}(K)$ of cell $K$}
    \If{$\textsf{touched}(F) = \textsf{false}$}
    \State {Assemble $A_{f \gets f}|_{F\gets F}$}
    \State {Set $\vec{w}^{(f)}|_F = A_{f \gets f}|_{F\gets F} \vec{u}^{(\pm)}|_F$}
    \Comment{Compute numerical fluxes}
    \State{Set $\textsf{touched}(F) = \textsf{true}$}
    \Comment{mark facet $F$ as touched}
    \EndIf
    \State {Assemble local matrix $A_{c\gets f}|_{K\leftarrow F}$}
    \State {$\vec{r}^{(c)}|_K \gets \vec{r}^{(c)}|_K - A_{c\gets f}|_{K\gets F} \vec{w}^{(f)}|_{F}$}
    \Comment {contribution from facets to residual}
    \EndFor
    \State {$\vec{u}^{(c)}|_K \gets \vec{u}^{(c)}|_K + \omega  A_{K\gets K}^{-1} \vec{r}^{(c)}|_K$}
    \Comment{update state}\label{alg:interior-penalty-single-level-one-mesh-sweep:update}
    \For{every facet $F\in \mathcal{F}(K)$ of cell $K$}\label{alg:interior-penalty-single-level-one-mesh-sweep:projection:start}
    \State {Assemble $A_{f \gets c}|_{K\gets F}$}
    \State {Set $\vec{u}^{(\pm)}|_F = A_{f \gets c}|_{F\gets K} \vec{u}^{(c)}|_K$}
    \Comment {Project solution onto facets}
    \EndFor\label{alg:interior-penalty-single-level-one-mesh-sweep:projection:finish}
    \EndFor \label{alg:interior-penalty-single-level-one-mesh-sweep:fused_loop:finish}
    \EndFor
    \State{\Return $\vec{u}^{(c)}, \vec{u}^{(\pm)}=A_{f\gets c}\vec{u}^{(c)}$}
  \end{algorithmic}
\end{algorithm}

\noindent 
In the code, the fusion of all three loops can be achieved by projecting the current solution $\vec{u}^{(c)}|_K$ onto all facets $F\in\mathcal{F}(K)$ of the cell $K$ as soon as it becomes available. This results in \cref{alg:interior-penalty-single-level-one-mesh-sweep}, which is mathematically equivalent to \cref{alg:vanilla-interior-penalty-single-level} and \cref{alg:interior-penalty-single-level-multiple-sweeps-and-auxiliary-variables}.
The projections $\vec{u}^{(\pm)}$, which are required to start the shifted sequence of operator evaluations according to \cref{technique:shift-and-fuse}, are computed in a warm-up mesh traversal (lines 
\ref{alg:interior-penalty-single-level-one-mesh-sweep:projection-preamble-start}--\ref{alg:interior-penalty-single-level-one-mesh-sweep:projection-preamble-end}) prior to the block-Jacobi iteration.
After that, each mesh traversal implements a smoothing step, i.e.~$n_{\text{iter}}$ smoothing steps can be realised with a total of $n_{\text{iter}}+1$ mesh traversals.

\begin{observation}
  Following the introduction of helper variables (\cref{technique:left-right-face-projections} and \cref{technique:numerical-flux-outcome}) and shift of operator evaluations (\cref{technique:shift-and-fuse}) the three mesh iterations in \cref{alg:interior-penalty-single-level-multiple-sweeps-and-auxiliary-variables} can be fused into a single mesh traversal in \cref{alg:interior-penalty-single-level-one-mesh-sweep}; the algorithm is weakly single touch in the sense of \cref{def:weak_single_touch}. 
\end{observation}

\noindent
The same domain decomposition as before can be used in \cref{alg:interior-penalty-single-level-one-mesh-sweep} and the projections $\vec{u}^{(\pm)}$ need to be exchanged between neighbouring processors before the start of the fused loop in lines \ref{alg:interior-penalty-single-level-one-mesh-sweep:fused_loop:start}-\ref{alg:interior-penalty-single-level-one-mesh-sweep:fused_loop:finish}.

\subsection{Task graphs}

\begin{table}
\caption{
 Types of tasks used in the block-Jacobi iteration in \cref{alg:interior-penalty-single-level-multiple-sweeps-and-auxiliary-variables} and \cref{alg:interior-penalty-single-level-one-mesh-sweep}.
 \label{tab:tasks}
}
  \begin{center}
\begin{tabular}{lcl}
  \hline
task type & used mesh entities & operation\\
\hline
Projection & $K \rightarrow F$ & $\vec{u}^{(\pm)}|_F = A_{f\gets c}|_{F\gets K}\vec{u}^{(c)}|_K$ \\
Numerical flux & $F \rightarrow F$ & $\vec{w}^{(f)}|_F = A_{f\gets f}|_{F\gets F}\vec{u}^{(\pm)}|_F$\\
Cell-residual & $K \rightarrow K$ & $\vec{r}^{(c)}|_K \gets \vec{b}^{(c)}|_K - A_{c\gets c}|_{K\gets K}\vec{u}^{(c)}|_K$\\
Facet-residual & $F \rightarrow K$ & $\vec{r}^{(c)}|_K \gets \vec{r}^{(c)}|_K - A_{c\gets f}|_{K\gets F}\vec{w}^{(f)}|_F$\\
Solution update & $K \rightarrow K$ & $\vec{u}^{(c)}|_K\mapsto \vec{u}^{(c)}|_K + \omega K_{K\gets K} \vec{r}^{(c)}|_K$\\
Matrix assembly & $K\;\text{or}\;F$ & Assemble $A_{c\gets c}|_{K\gets K}$, $A_{c\gets f}|_{K\gets F}$ or $A_{f\gets f}|_{F\gets F}$\\
Matrix inversion & $K$ & Invert $A_{K\gets K}$\\\hline
\end{tabular}
\end{center}
\end{table}

The two implementations of the block-Jacobi iteration in \cref{alg:interior-penalty-single-level-multiple-sweeps-and-auxiliary-variables} and \cref{alg:interior-penalty-single-level-one-mesh-sweep} can be written down as a task graph, where each task corresponds to an operation on a single mesh entity, i.e.~a cell or a facet (\cref{tab:tasks}). 
The resulting directed graph expresses the dependencies between different tasks.
Each task can (but does not have to) be executed once it is ``ready'', i.e. once all the other tasks that it depends on have completed. While the task graph formalism can be applied at a finer granularity, e.g.~by breaking down the matrix-vector products into further smaller tasks \cite{badia2010scheduling} or by applying task paradigms to assembly steps as well~\cite{murray2021assembly}, we refrain from such a fine-granular decomposition here and only consider tasks that operate on data associated with entire cells or facets shown in \cref{tab:tasks}: these tasks naturally map onto BLAS routines, which form building blocks of appropriate granularity for our application (cmp.~\cref{def:vectorisation}).

\begin{figure}[htb]
  \begin{center}
    \includegraphics[width=0.9\textwidth]{\figdir/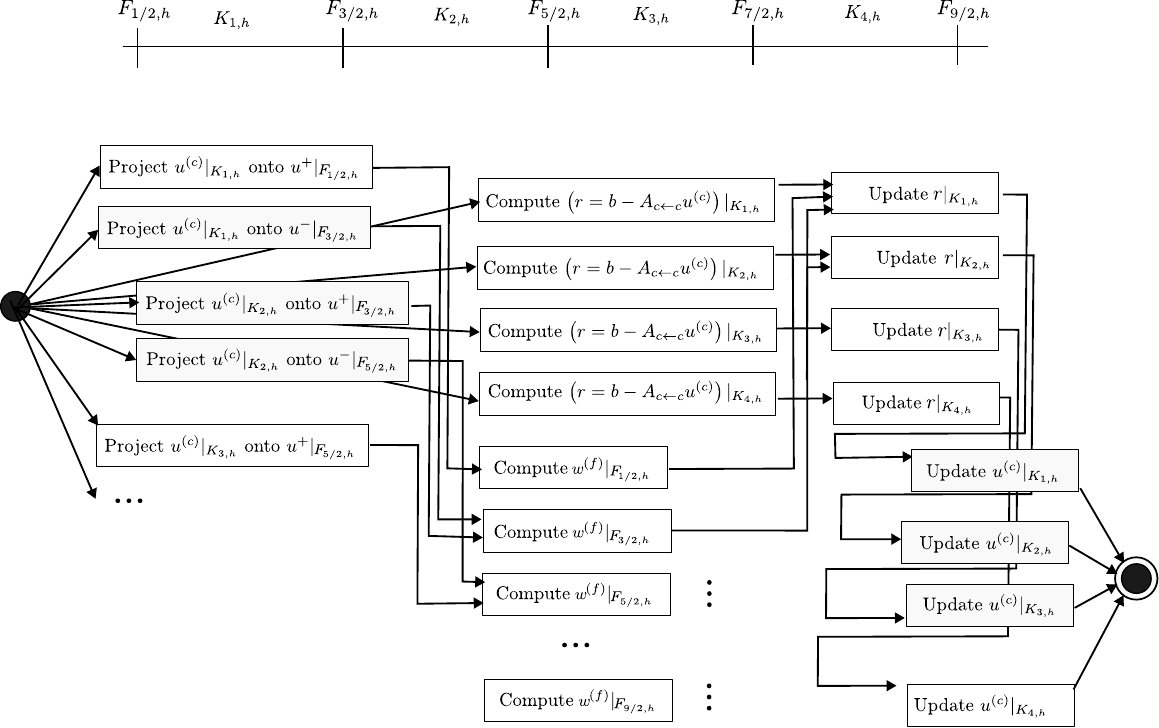}
  \end{center}
  \caption{
    Schematic task graph (bottom) for a single block-Jacobi iteration applied to a 1d problem with four mesh cells (top). 
    Fractional indices are used to number the facets sitting in-between two cells. 
    \cref{tab:tasks} breaks down the algorithm into more detailed tasks, which are omitted here for brevity. 
    \label{figure:task-formalism:single-time-step-vanilla-penality-method}
  }  
\end{figure}

\subsubsection{Direct mapping of algorithms onto a task language}

\paragraph{Structure and character of the task graph}
It is instructive to visualise the task graph by considering the order in which the tasks are spawned by \cref{alg:interior-penalty-single-level-multiple-sweeps-and-auxiliary-variables}. 
If we ignore the assembly of local matrices for a moment, each of the three mesh traversals spawns bursts of tasks of a particular type. 
The first mesh sweep over cells in lines \ref{alg:interior-penalty-single-level-multiple-sweeps-and-auxiliary-variables:first-loop-block:start}--\ref{alg:interior-penalty-single-level-multiple-sweeps-and-auxiliary-variables:first-loop-block:end}
generates exclusively ``Projection''-type tasks that map a solution $\vec{u}^{(c)}$ onto the facets to obtain $\vec{u}^{(\pm)}$. 
The second traversal over facets in lines \ref{alg:interior-penalty-single-level-multiple-sweeps-and-auxiliary-variables:second-loop-block:start}--\ref{alg:interior-penalty-single-level-multiple-sweeps-and-auxiliary-variables:second-loop-block:end} spawns ``Numerical flux''-type tasks to compute $\vec{w}^{(f)}$ from $\vec{u}^{(\pm)}$. 
The third and final iteration over mesh cells in lines \ref{alg:interior-penalty-single-level-multiple-sweeps-and-auxiliary-variables:third-loop-block:begin} - \ref{alg:interior-penalty-single-level-multiple-sweeps-and-auxiliary-variables:iteration:end} spawns ``Cell-residual''- and ``Facet-residual''-type tasks to update $\vec{r}^{(c)}$ and ``Solution update''-type tasks to compute the new state $\vec{u}^{(c)}$. The task graph that is naturally drawn based on this order of spawning tasks is sketched for a particular one-dimensional setup in \cref{figure:task-formalism:single-time-step-vanilla-penality-method}. Although as topological objects the task-graphs of \cref{alg:interior-penalty-single-level-multiple-sweeps-and-auxiliary-variables} and \cref{alg:interior-penalty-single-level-one-mesh-sweep} are identical, the order in which they are constructed by the two algorithms differs. The task creation pattern in \cref{alg:interior-penalty-single-level-one-mesh-sweep} is more complex in the  sense that it results in the spawning of more tasks of the same type in a single loop: while the first mesh traversal in lines \ref{alg:interior-penalty-single-level-one-mesh-sweep:projection-preamble-start} - \ref{alg:interior-penalty-single-level-one-mesh-sweep:projection-preamble-end} (which can be considered as a ``warm-up'' phase) exclusively spawns task of the ``Projection''-type, the subsequent mesh sweeps in lines \ref{alg:interior-penalty-single-level-one-mesh-sweep:fused_loop:start} - \ref{alg:interior-penalty-single-level-one-mesh-sweep:fused_loop:finish} spawn a mixture of all task types (as given in \cref{tab:tasks}) due to the application of \cref{technique:shift-and-fuse}.

\paragraph{Concurrency analysis}

In contrast to \cref{alg:vanilla-interior-penalty-single-level}, the volumetric cell operations and therefore the associated tasks are decoupled from each other due to the projection of the solution onto facets (\cref{technique:left-right-face-projections} and \cref{technique:numerical-flux-outcome}): the operations in one particular cell do not depend directly on the solutions in adjacent cells. 
More generally:

\begin{observation}
 \label{observation:tasking:independent-tasks-of-same-type}
 The introduction of auxiliary variables makes all tasks of one particular type given in \cref{tab:tasks} independent of each other.
\end{observation}

\paragraph{Two extreme strategies for task execution} Having constructed the task graph, we need to decide how and in which order to schedule the tasks for execution, bearing in mind their mutual dependencies. In a task language, executing \cref{alg:interior-penalty-single-level-multiple-sweeps-and-auxiliary-variables} and \cref{alg:interior-penalty-single-level-one-mesh-sweep} line by line is equivalent to executing the tasks in the order in which they are encountered while looping over the mesh entities. In this case, tasking can be considered as a logical abstraction: there is no need to construct the task graph explicitly since all dependencies are implicitly satisfied. We can interpret the mesh traversals in \cref{alg:interior-penalty-single-level-multiple-sweeps-and-auxiliary-variables} and \cref{alg:interior-penalty-single-level-one-mesh-sweep} as task producers that loop over mesh entities and spawn ready tasks, i.e.~tasks with no pending dependencies as all dependencies are implicitly fulfilled, and then execute them immediately. In this sense, the mesh traversals fulfil the dual roles of task-producer and task-scheduler. This can be seen as a particular way of scheduling: the tasks in the graph are executed in a fixed deterministic order. The other extreme would be to assemble the task graph and then leave the execution of the tasks in the graph to a completely separate runtime system such as OpenMP or TBB. Instead of executing the tasks immediately, the mesh traversals in \cref{alg:interior-penalty-single-level-multiple-sweeps-and-auxiliary-variables} and \cref{alg:interior-penalty-single-level-one-mesh-sweep} will not process the instructions they encounter, but instead map each operation to a  ``physical'' task, i.e. a set of instructions together with the dependencies of the input variables on other variables. The separate runtime system then executes the physical tasks, possibly in a non-deterministic order, bearing in mind the dependencies encoded in the task graph.

The latter execution strategy, which involves the explicit construction of the task graph, exposes the maximum concurrency level (cmp.~\cref{def:concurrency}) of the numerical scheme. The downside of this approach is that we (i) break the single touch semantics (\cref{def:weak_single_touch}), (ii) introduce some memory overhead (\cref{def:memory_overhead}) and (iii) introduce additional costs for task dependency management, which is particularly disadvantageous for computationally inexpensive tasks such as the numerical flux calculations.

Notably, since we have no control of the execution order of ready tasks and their assignment to cores, the volumetric data required for the cell residual calculations and the assembly might be moved accross the bus multiple times due to capacity, coherence and conflict cache misses.

\begin{observation}
  \label{obs:task_management_overhead}
We find empirically that a complete separation between task-graph generation and task-execution with a runtime system results in non-competitive performance. The additional cost from task-management is too high and there are too many computationally inexpensive tasks with a disproportionally large overhead.
\end{observation}
\noindent
This observation is in line with other studies.
Notably large bursts of tiny, interdependent tasks such as those of ``Projection''-, ``Numerical flux''- and ``Facet-residual''-type challenge modern runtimes \cite{Tuft:2024:DetrimentalTaskPatterns}.

To avoid the issues raised in \cref{obs:task_management_overhead} while still exploiting the advantages of the task-based approach we develop a hybrid execution model which is an intermediate between the two extreme scheduling strategies described above. To motivate this, observe that data parallelism can be also be interpreted in the context of a task language.
To see this, consider \cref{alg:interior-penalty-single-level-multiple-sweeps-and-auxiliary-variables} and assume that the ``Solution-update''-, ``Cell-residual''- and ``Facet-residual''-type tasks, i.e.~all tasks arising in the loop in lines \ref{alg:interior-penalty-single-level-multiple-sweeps-and-auxiliary-variables:first-loop-block:start}-\ref{alg:interior-penalty-single-level-multiple-sweeps-and-auxiliary-variables:first-loop-block:end} are fused into a single task operating on a cell.
In this case, each of the three loops in lines \ref{alg:interior-penalty-single-level-multiple-sweeps-and-auxiliary-variables:second-loop-block:start}-\ref{alg:interior-penalty-single-level-multiple-sweeps-and-auxiliary-variables:second-loop-block:end}, \ref{alg:interior-penalty-single-level-multiple-sweeps-and-auxiliary-variables:third-loop-block:begin}-\ref{alg:interior-penalty-single-level-multiple-sweeps-and-auxiliary-variables:iteration:end} and \ref{alg:interior-penalty-single-level-multiple-sweeps-and-auxiliary-variables:first-loop-block:start}-\ref{alg:interior-penalty-single-level-multiple-sweeps-and-auxiliary-variables:first-loop-block:end} can be executed in parallel. This is possible since in each loop operations on different mesh entities can be executed independently without any write conflicts according to~\cref{observation:tasking:independent-tasks-of-same-type}. The organisation into three separate loops implicitly imposes global synchronisation points at the end of each mesh traversal. In the context of a task language, this corresponds to spawning all tasks during the mesh traversal, but waiting for the runtime system to execute all tasks before proceeding to the next mesh traversal. This now requires only minimal overhead from task management since all tasks are ready and can be executed independently: in practice, it is not necessary to explicitly construct the task graph.

\begin{definition}
 \label{definition:task-graph-assembly-free}
 A \emph{task-graph-construction-free} execution model logically employs a task graph, but introduces global synchronisation points and spawns the tasks in an order which guarantees that all dependencies are implicitly fulfilled: all spawned tasks are ``ready'' by definition and can be executed independently. 
\end{definition}

\noindent
While this approach eliminates overheads from task management and results in very high concurrency per mesh traversal (\cref{def:concurrency}), it still has a serious drawback:

\begin{observation}\label{obs:task_synchronisation}
Imposing global synchronisation points at the end of each mesh traversal limits parallel scalability.
\end{observation}

\noindent
Because of this drawback we do not pursue the data-parallel execution model based on \cref{alg:interior-penalty-single-level-multiple-sweeps-and-auxiliary-variables} any further here. Instead, we use the derived insights and propose a task-based execution strategy that is almost task-graph-construction-free. The resulting approach balances between modelling calculations as tasks with dependencies and the direct execution of tasking during the mesh traversals. 

\subsubsection{Almost task-graph construction-free execution model}

The block-Jacobi iteration written down in \cref{alg:interior-penalty-single-level-one-mesh-sweep} collapses all three mesh traversals into a single loop. This improves data locality and reduces the number of synchronisation points by a factor three, thereby addressing the issues in \cref{obs:task_synchronisation}.

We could again avoid task management overheads by employing a task-graph construction-free execution model in the sense of \cref{definition:task-graph-assembly-free}. For this, some of the individual task in lines \ref{alg:interior-penalty-single-level-one-mesh-sweep:fused_loop:start}-\ref{alg:interior-penalty-single-level-one-mesh-sweep:fused_loop:finish} of \cref{alg:interior-penalty-single-level-one-mesh-sweep} would need to be combined into larger ``meta-tasks'' in such a way that each of the resulting tasks can be executed independently of all other tasks. Since then by construction all tasks are ready and can be executed in parallel, it is again not necessary to construct the task graph.

\paragraph{Hybrid execution model}
The final execution model we consider is based on \cref{alg:interior-penalty-single-level-one-mesh-sweep}, which reduces the number of synchronisation points and improves data locality compared to \cref{alg:interior-penalty-single-level-multiple-sweeps-and-auxiliary-variables}. However, instead of applying a pure task-graph construction-free approach (\cref{definition:task-graph-assembly-free}), we organise the tasks into two categories: during the mesh traversal in lines \ref{alg:interior-penalty-single-level-one-mesh-sweep:fused_loop:start}-\ref{alg:interior-penalty-single-level-one-mesh-sweep:fused_loop:finish} of \cref{alg:interior-penalty-single-level-one-mesh-sweep}, most tasks are executed immediately when they are encountered and a domain decomposition strategy is used for parallelisation. In contrast, tasks that fall into the second category are spawned and passed on to the runtime system, which is responsible for their execution according to the dependency graph. However, the task graph that needs to be managed by the runtime system is relatively (and possibly trivial). As a consequence, the issues in \cref{obs:task_management_overhead} are avoided, in particular if the tasks in the second category are chosen such that they are computationaly expensive and have a disproportionally small management overhead.

More specifically, we employ the following:
\begin{technique}
 \label{technique:residual-and-inversion-tasks}
 Only two types of computationally expensive tasks, the computation of $A_{K\gets K}^{-1}$ (``Matrix inversion'') and the ``Cell-residual'' calculation, are spawned as ``physical'' tasks $T_K^{\text{(inv)}}$ and $T_K^{\text{(residual)}}$ respectively. All other tasks are executed immediately during the mesh traversal.
\end{technique}

\noindent
Since they are independent, the tasks $T_K^{\text{(inv)}}$ and $T_K^{\text{(residual)}}$ can be executed in parallel for different cells $K$ due to \cref{observation:tasking:independent-tasks-of-same-type}. They are by definition ready, i.e.~the approach remains task-graph construction-free (\cref{definition:task-graph-assembly-free}). For homogeneous systems $T_K^{\text{inv}}$ only needs to be executed once at the beginning of the simulation since $A_{K\gets K}$ does not vary across the domain and can be precomputed.

\paragraph{Elimination of global synchronisation points}
A straightforward implementation with \cref{technique:residual-and-inversion-tasks} results in one global synchronisation point at the end of the loop in lines \ref{alg:interior-penalty-single-level-one-mesh-sweep:fused_loop:start}-\ref{alg:interior-penalty-single-level-one-mesh-sweep:fused_loop:finish} of \cref{alg:interior-penalty-single-level-one-mesh-sweep}. As the temporal shifts in \cref{technique:shift-and-fuse} imply that the outcomes of $T_K^{\text{(inv)}}$ and $T_K^{\text{(residual)}}$ are not required prior to the next mesh sweep, this explicit global synchronisation point is not necessary. Instead, we can wait for $T_K^{\text{(inv)}}$ and $T_K^{\text{(residual)}}$ to complete when processing cell $K$ in the next mesh traversal:
\begin{algorithm}
  \caption{Matrix-free block-Jacobi iteration using loop fusion plus tasking for the cell-residual calculation and local matrix inversion. $\text{BlockJacobi}(\vec{u}^{(c)}_0, \vec{b}^{(c)};\omega,n_{\text{iter}},\epsilon)$ Input: initial solution $\vec{u}^{(c)}_0$, right-hand side $\vec{b}^{(c)}$, relaxation parameter $\omega$, number of iterations $n_{\text{iter}}$. Output: solution $\vec{u}^{(c)}$ and its facet-projection $\vec{u}^{(\pm)}=A_{f\gets c}\vec{u}^{(c)}$ after $n_{\text{iter}}$ iterations.}
  \label{alg:interior-penalty-single-level-tasks}
  \begin{algorithmic}[1]
    \For{every cell $K \in \Omega _h$}
      \If {$A_{K\gets K}$ is not constant in time and space}
        \State {spawn $task(\text{compute}\; A_{K\gets K}^{-1})=:T_{K}^{(\text{inv})}$}\Comment{Local matrix inversion}\label{line_alg:tasking_spawn_inversion_preamble}
      \EndIf
      \State {spawn $task \left( \vec{r}^{(c)}|_K \gets \vec{b}^{(c)}|_K -  A_{c\gets c}|_{K\gets K} \vec{u}^{(c)}|_K \right)=:T_K^{(\text{residual})}$}\label{line_alg:tasking_spawn_cell_residual_preamble}
        \Comment {Cell-residual}
      \For{every facet $F\in \mathcal{F}(K)$ of cell $K$}
        \State {Assemble $A_{f \gets c}|_{K\gets F}$}
        \State {Set $\vec{u}^{(\pm)}|_F = A_{f \gets c}|_{F\gets K} \vec{u}^{(c)}|_K$}
          \Comment {Project solution onto facets}
      \EndFor      
    \EndFor
    \For{$k=1,2,\dots,n_{\text{iter}}$}
      \State {Exchange $\vec{u}^{(\pm)}$ between non-overlapping subdomains.}
      \State {Set $\textsf{touched}(F) = \textsf{false}$ for all $F\in\mathcal{E}_h$}
        \Comment{Mark all facets as untouched}
      \For{every cell $K \in \Omega _h$}
        \State {Wait for tasks $T_K^{(\text{inv})}$ and $T_K^{(\text{residual})}$ to compute $A_{K\gets K}^{-1}$ and $\vec{r}^{(c)}|_K$ in cell $K$}\label{line_alg:tasking_wait}
    \For{every facet $F\in\mathcal{F}(K)$ of cell $K$}
    \If{$\textsf{touched}(F) = \textsf{false}$}
    \Comment{Only compute numerical flux once on each facet}
    \State {Assemble $A_{f \gets f}|_{F\gets F}$}
    \State {Set $\vec{w}^{(f)}|_F = A_{f \gets f}|_{F\gets F} \vec{u}^{(\pm)}|_F$}
    \Comment{Compute numerical fluxes}
    \State{Set $\textsf{touched}(F) = \textsf{true}$}
    \Comment{mark facet $F$ as touched}
    \EndIf
    \State {Assemble local matrix $A_{c\gets f}|_{K\leftarrow F}$}
    \State {$\vec{r}^{(c)}|_K \gets \vec{r}^{(c)}|_K - A_{c\gets f}|_{K\gets F} \vec{w}^{(f)}|_{F}$}
    \Comment {contribution from facets to residual}
    \EndFor
    \State {$\vec{u}^{(c)}|_K \gets \vec{u}^{(c)}|_K + \omega  A_{K\gets K}^{-1} \vec{r}^{(c)}|_K$}
    \Comment{update state}
    \For{every facet $F\in \mathcal{F}(K)$ of cell $K$}
      \State {Assemble $A_{f \gets c}|_{K\gets F}$}
      \State {Set $\vec{u}^{(\pm)}|_F = A_{f \gets c}|_{F\gets K} \vec{u}^{(c)}|_K$}
        \Comment {Project solution onto facets}
    \EndFor
      \If {$A_{K\gets K}$ is not constant in time and space}
      \State {spawn $task(\text{compute}\;A_{K\gets K}^{-1})=:T_K^{(\text{inv})}$}
        \Comment {Local matrix inversion}\label{line_alg:tasking_spawn_inversion}
      \EndIf
      \State {spawn $task \left( \vec{r}^{(c)}|_K \gets \vec{b}^{(c)}|_K -  A_{c\gets c}|_{K\gets K} \vec{u}^{(c)}|_K \right)=:T_K^{(\text{residual})}$}
        \Comment {Cell-residual}\label{line_alg:tasking_spawn_cell_residual}
    \EndFor
    \EndFor
    \State{\Return $\vec{u}^{(c)}, \vec{u}^{(\pm)}=A_{f\gets c}\vec{u}^{(c)}$}
  \end{algorithmic}
\end{algorithm}

\noindent
\begin{technique}
 \label{technique:split-and-postpone}
 We split the computations per cell into computations whose output is required at the end of the present mesh sweep and execute these immediately. The remaining calculations are spawned as separate, physical tasks and handed to the runtime system. In the next traversal, we wait for the completion of these tasks before executing the cell calculations that depend on them. This way, we \emph{split and postpone} some calculations through a task formalism. 
\end{technique}
This results in the implementation shown as \cref{alg:interior-penalty-single-level-tasks}, which is identical to \cref{alg:interior-penalty-single-level-one-mesh-sweep} except for lines \ref{line_alg:tasking_spawn_inversion_preamble}, \ref{line_alg:tasking_spawn_inversion} and \ref{line_alg:tasking_spawn_cell_residual_preamble}, \ref{line_alg:tasking_spawn_cell_residual} which spawn $T_K^{\text{(inv)}}$ and $T_K^{\text{(residual)}}$ respectively and line \ref{line_alg:tasking_wait} which waits for the tasks to complete.
\noindent
\cref{technique:split-and-postpone} is an antagonist to the fusion in \cref{technique:shift-and-fuse}:
We break up big volumetric, cell-wise tasks resulting from the loop fusion in \cref{alg:interior-penalty-single-level-one-mesh-sweep} and execute the
``Facet residual'', ``solution update'' and ``Facet projection'' tasks without a task framework. While many individual tasks are now again executed immediately during a mesh traversal, it is left to the runtime system to decide when to complete the remaining tasks $T_K^{\text{(inv)}}$ and $T_K^{\text{(residual)}}$, as long as the outcome is available before we carry out further calculations in cell $K$ in the next mesh traversal after line \ref{line_alg:tasking_wait}.

\paragraph{Properties}

Although the realisation of \cref{technique:split-and-postpone} in \cref{alg:interior-penalty-single-level-tasks} breaks with some of the paradigms described in \cref{sec:implementation_with_facet_variables} and \cref{sec:single-touch-grid-traversal}, it improves important aspects of the implementation:

\begin{enumerate}
  \item The implementation is no longer based on a pure task-based approach, as some ``urgent'' calculations are immediately executed during the mesh traversal. In this respect it is similar to an approach which does not employ any task-based modelling at all. The hybrid approach in \cref{alg:interior-penalty-single-level-tasks} avoids the management of computationally inexpensive tasks with a disproportionally large overhead, yet it results in an increased level of concurrency as advocated by \cref{def:concurrency}.
  \item In contrast, the approach is not entirely task-free since some parts of a task graph are assembled. All spawned tasks are ready by construction, but we have to take care when executing other operations that depend on these tasks. To achieve this, the task runtime is queried in a subsequent mesh traversal to check whether the ``Matrix inversion'' and ``Cell-residual'' tasks have completed. This results in some overhead due to task management, but this is not as sizeable as in a purely task-based execution model.
  \item As we compute the inverse of $A_{K\gets K}$ in parallel to all other operations in a mesh traversal and need to hold these values until they are needed to update the solution $\vec{u}^{(c)}$, additional temporary variables are introduced. This increases the memory footprint non-deterministically (cmp.~\cref{def:memory_overhead}). An analoguous argument holds for the computation of the cell-residual.
    \item All explicit, global task synchronisation points are eliminated. At the end of a mesh traversal, all computations required for this traversal have implicitly completed, while some of the fused calculations (which logically belong into the subsequent traversal) are mapped onto tasks, whose completion is not required to initiate the next mesh traversal.
\end{enumerate}

\noindent
The last property implies that the implementation has weak synchronisation points: when a mesh traversal completes there might still be pending tasks within the system that have not completed yet. This is not an issue since the output will only be required at some later point in the subsequent mesh traversal. As a consequence, we never reduce the concurrency level to one.

\paragraph{Domain decomposition and parallelisation}

To use a traditional domain decomposition approach in a distributed memory setting, the values $\vec{u}^{(\pm)}$ on the facets between adjacent subdomains have to remain consistent. In a purely task based approach where the execution of all tasks is handled by the runtime system, this would require careful attention and a dedicated synchronisation mechanisms which involves parallel communication: the data exchange on subdomain boundaries cannot be triggered before all the projections have finished. Even if the domain decomposition is implemented on a shared memory system a purely task-based approach will introduce complicated synchronisation issues.

\begin{observation}
 The hybrid execution model employed here only leaves the execution of purely volumetric tasks to the runtime system. Since these tasks are independent of $\vec{u}^{(\pm)}$, the approach does not interfer with domain decomposition.
\end{observation}
\noindent
In a distributed memory setting, the runtime system that handles the tasks $T_K^{\text{(inv)}}$ and $T_K^{\text{(residual)}}$ can run independently on each processor and does not require parallel communication.

\section{Multigrid}
\label{section:multigrid}
Although the stationary block-Jacobi iteration \eqref{eqn:block-Jacobi-IP} implemented in \cref{alg:interior-penalty-single-level-one-mesh-sweep} converges for suitable values of $\omega$, the convergence rate is mesh-dependent and deteriorates as the resolution increases: the finer the mesh the slower the convergence. The reason for this is that components of the error which vary slowly over the grid will not be reduced efficiently by the block Jacobi-iteration which is inherently local.
To address this issue, we use the $hp$-multigrid solver described in \cite{bastian2019matrix,bastian2012algebraic}: the block-Jacobi iterations on the finest level reduce high-frequency components of the error, while the slowly varying error components are eliminated by solving the residual equation on a hierarchy of lower dimensional subspaces. 
Since the block-Jacobi smoother is efficient at eliminating components of the error that fluctuate within each cell, the first coarsening step reduces the polynomial degree into the lowest order continuous space over the same mesh. After this initial $p$-coarsening step, we follow a traditional $h$-coarsening strategy which increases the grid spacing by a constant factor of three, i.e.~we exploit the space-tree structure of the mesh.

\begin{technique}
 \label{technique:weaker-block-smoother}
 For the present $hp$-multigrid scheme, an \emph{overlapping} block-Jacobi smoother is required to guarantee $p$-robustness \cite{bastian2012algebraic}.
 Here we substitute it with the weaker non-overlapping block-Jacobi iteration from \eqref{eqn:block-Jacobi-IP} which has a smaller memory movement imprint and which has also been used in \cite{bastian2019matrix}.
\end{technique}

\subsection{Two-level method}

A two-level method uses the fine level $\VDG{h}{p}$ plus the coarse level $\VCG{h}{1}\subset \VDG{h}{p}$, i.e. the subspace of continuous piecewise linear functions on the same mesh.
Since the function spaces are nested, the prolongation 

\[
  P:\vec{u}^{(\text{coarse})} \mapsto \vec{u}^{(c)}=P\vec{u}^{(\text{coarse})}
\]

\noindent
of a dof-vector $\vec{u}^{(\text{coarse})}$ on the coarse level onto a dof-vector $\vec{u}^{(c)}$ on the fine level is naturally defined by requiring that $\vec{u}^{(c)}$ and $\vec{u}^{(coarse)}$ represent the same function.

The corresponding restriction $R = P^\top$ for dual vectors is given by the transpose of $P$. Similar to the stiffness matrix, we can partition the prolongation matrix into local blocks that couple the unknowns associated with a cell and its vertices. With this we can write the prolongation in each cell $K$ as follows
\begin{equation}
  \vec{u}^{(c)}|_K = \sum_{V\in\mathcal{V}(K)} P|_{K\gets V} \vec{u}^{(\text{coarse})}|_V,
  \label{equation:multigrid:prolongation}
\end{equation}
where the small matrix $P|_{K\gets V}$ maps the unknowns associated with the vertex $V$ to the unknowns in cell $K$. 

For a given $\vec{u}^{(c)}$ the error $A^{-1}\vec{r}^{(c)}$ with $\vec{b}^{(c)}-A\vec{u}^{(c)}$ can be approximated by the coarse level correction $\delta \vec{u}^{(c)}:=P(A^{(\text{coarse})})^{-1}P^\top \vec{r}^{(c)} = P(P^\top A P)^{-1}P^\top \vec{r}^{(c)}$, which can be used to improve the current solution $\vec{u}^{(c)}$. The coarse level correction $\delta \vec{u}^{(c)}$ can be computed in four phases as follows:
\begin{enumerate}
  \item Compute the DG residual $\vec{r}^{(c)} = \vec{b}^{(c)}-A\vec{u}^{(c)}$;
  \item restrict the residual to the coarse level to obtain $\vec{b}^{(\text{coarse})}=P^\top\vec{r}^{(c)}$;
  \item solve the coarse level equation $A^{(\text{coarse})}\vec{e}^{(\text{coarse})}=\vec{b}^{(\text{coarse})}$ for $\vec{e}^{(\text{coarse})}$;
  \item prolongate the coarse level solution back to the fine level to obtain $\delta \vec{u}^{(c)}=P \vec{e}^{(\text{coarse})}$.
\end{enumerate}

\begin{algorithm}[htb]
  \caption{$\text{CoarseGridCorrection}(\vec{u}^{(c)},\vec{u}^{(\pm)}, \vec{b}^{(c)})$. Compute coarse grid correction $\delta\vec{u}^{(c)}=P\vec{e}^{(\text{coarse})}$ with $A^{(\text{coarse})}\vec{e}^{(\text{coarse})}=P^\top \vec{r}^{(c)}$.  Input: solution $\vec{u}^{(c)}$ and its projection to facets $\vec{u}^{(\pm)}=A_{f\gets c}\left(\vec{u}^{(c)}\right)$, right hand side $\vec{b}^{(c)}$. Output: coarse grid correction $\delta\vec{u}^{(c)}$.}
  \label{alg:coarse-grid correction}
  \begin{algorithmic}[1]        
    \State {Set $\vec{b}^{(\text{coarse})}=0$}    
    \State {Set $\textsf{touched}(F) = \textsf{false}$ for all $F\in\mathcal{E}_h$}
    \Comment{Mark all facets as untouched}    
    \For{every cell $K \in \Omega _h$}
    \State {Assemble $A_{c \gets c}|_{K\gets K}$}
    \Comment {On-the-fly assembly of cell-local matrix}
    \State {$\vec{r}^{(c)}|_K \gets \vec{b}^{(c)}|_K -  A_{c\gets c}|_{K\gets K} \vec{u}^{(c)}|_K$}
    \Comment {cell-local contribution to residual}
    \For{every facet $F\in\mathcal{F}(K)$ of cell $K$}
    \If{$\textsf{touched(F)} = \textsf{false}$}
    \Comment{Only compute numerical flux once on each facet}
    \State {Assemble $A_{f \gets f}|_{F\gets F}$}
    \State {Set $\vec{w}^{(f)}|_F = A_{f \gets f}|_{F\gets F} \vec{u}^{(\pm)}|_F$}
    \Comment{Compute numerical fluxes}
    \State{Set $\textsf{touched(F)} = \textsf{true}$}
    \Comment{mark facet $F$ as touched}
    \EndIf
    \State {Assemble local matrix $A_{c\gets f}|_{K\leftarrow F}$}
    \State {$\vec{r}^{(c)}|_K \gets \vec{r}^{(c)}|_K - A_{c\gets f}|_{K\gets F} \vec{w}^{(f)}|_{F}$}
    \Comment {contribution from facets to residual}
    \EndFor
    \For{every vertex $V\in\mathcal{V}(K)$ of cell $K$}
    \State{Update $\vec{b}^{(\text{coarse})}|_V \gets \vec{b}^{(\text{coarse})}|_V + P|_{K\gets V} \vec{r}^{(c)}|_K$}
    \Comment{Restrict residual}
    \EndFor
    \EndFor
    \State{(Approximately) solve $A^{\text{(coarse)}} \vec{e}^{(\text{coarse})} = \vec{b}^{(\text{coarse})}$ for $\vec{e}^{(\text{coarse})}$}
      \Comment{Coarse grid solve}
      \label{alg:coarse-grid correction:coarse-grid-solve}
    \State{Set $\delta \vec{u}^{(c)}\gets 0$}
    \For{every cell $K \in \Omega _h$}
    \For{every vertex $V\in\mathcal{V}(K)$ of cell $K$}
    \State{Update $\delta \vec{u}^{(c)}|_K \gets \delta\vec{u}^{(c)}|_K + P|_{K\gets V} \vec{e}^{(\text{coarse})}|_V$}
    \Comment{Prolongate}
    \EndFor    
    \EndFor
    \State{\Return $\delta\vec{u}^{(c)}$}
  \end{algorithmic}
\end{algorithm}

\noindent
To obtain an iterative two-level solver, the block-Jacobi smoother and the coarse grid correction in \cref{alg:coarse-grid correction} are interleaved: After $\nu$ block-Jacobi iterations the coarse grid solution computed from the residual $\vec{r}^{(c)}=\vec{b}^{(c)}-A\vec{u}^{(c)}$ is used to construct an improved solution $\vec{u}^{(c)}+\delta \vec{u}^{(c)}$. This process is repeated iteratively as shown in \cref{alg:multiplicative_multigrid}, which is the well-known (multiplicative) multigrid algorithm with $\nu$ pre- and zero post-smoothing steps.
\begin{algorithm}
  \caption{Multiplicative multigrid. $\text{MGMult}(\vec{u}^{(c)}_0, \vec{b}^{(c)};\omega,\nu,n_{\text{iter}})$ Input: initial solution $\vec{u}^{(c)}_0$, right hand side $\vec{b}^{(c)}$, relaxation parameter $\omega$, number of smoothing steps $\nu$, number of iterations $n_{\text{iter}}$, tolerance $\epsilon$. Output: solution $\vec{u}^{(c)}$ after $n_{\text{iter}}$ iterations or convergence to tolerance $\epsilon$.}
  \label{alg:multiplicative_multigrid}
  \begin{algorithmic}[1]
    \State{Set $\vec{u}^{(c)}\gets \vec{u}^{(c)}_0$}
    \For{$k=1,2,\dots,n_{\text{iter}}$}
    \State{Set $\vec{u}^{(c)}_{\text{old}}\gets\vec{u}^{(c)}$}
    \State{$\vec{u}^{(c)},\vec{u}^{(\pm)}\gets\text{BlockJacobi}(\vec{u}^{(c)}, \vec{b}^{(c)};\omega,\nu,0)$}
    \Comment{Smoothing}
    \State{$\vec{u}^{(c)}\gets\vec{u}^{(c)}+\text{CoarseGridCorrection}(\vec{u}^{(c)},\vec{u}^{(\pm)},\vec{b}^{(c)})$}    
    \Comment{Add coarse grid correction}
    \State{Compute $\rho_k = \vert|\vec{u}^{(c)}_{\text{old}}-\vec{u}^{(c)}\vert|$}\Comment{preconditioned residual norm}
    \If{$\rho_k/\rho_1<\epsilon$}
    \State{\textbf{exit loop}}\Comment{check convergence}
    \EndIf
    \EndFor
    \State{\Return $\vec{u}^{(c)}$}
  \end{algorithmic}
\end{algorithm}

\subsubsection{Efficient implementation: auxiliary facet variables and loop fusion}

By introducing additional variables on the facets of the mesh (\cref{technique:left-right-face-projections} and \cref{technique:numerical-flux-outcome}), the four phases required to compute the coarse level correction $\delta\vec{u}^{(c)}$ can be realised as traversals over the cells of the mesh (\cref{alg:coarse-grid correction}).
A matrix-free, parallel, efficient realisation of the coarse grid solve follows (degenerated) DG techniques \cite{weinzierl2011peano}. As the following discussion shows, the number of mesh traversals can be minimised by employing loop fusion.

\paragraph{Restriction}

The block-Jacobi scheme in \cref{alg:interior-penalty-single-level-one-mesh-sweep} and \cref{alg:interior-penalty-single-level-multiple-sweeps-and-auxiliary-variables} includes the construction of $\vec{r}^{(c)}=\vec{b}^{(c)}-A\vec{u}^{(c)}$. We can therefore directly re-use one of these algorithms to compute the residual that is to be restricted to the coarse level if we omit the step that updates the solution $\vec{u}^{(c)}\gets \vec{u}^{(c)}+\omega A_{K\gets K}^{-1}\vec{r}^{(c)}$. More specifically, we can use \cref{alg:interior-penalty-single-level-multiple-sweeps-and-auxiliary-variables} without the update of the solution in line \ref{alg:interior-penalty-single-level-multiple-sweeps-and-auxiliary-variables:update-state}, or remove the solution update in line \ref{alg:interior-penalty-single-level-one-mesh-sweep:update} and the projection in lines \ref{alg:interior-penalty-single-level-one-mesh-sweep:projection:start} - \ref{alg:interior-penalty-single-level-one-mesh-sweep:projection:finish} from \cref{alg:interior-penalty-single-level-one-mesh-sweep}; analogous modifications can be made to \cref{alg:interior-penalty-single-level-tasks}.

Using \cref{alg:interior-penalty-single-level-multiple-sweeps-and-auxiliary-variables} (which requires three mesh traversals per iteration) to perform $\nu$ block-Jacobi iterations followed by one residual calculation requires $3(\nu+1)$ mesh traversals overall ($3\nu$ mesh traversals in the block-Jacobi smoother and three mesh traversals to compute the residual). To perform the same sequence of operations with \cref{alg:interior-penalty-single-level-one-mesh-sweep} requires $\nu+2$ mesh traversals: lines \ref{alg:interior-penalty-single-level-one-mesh-sweep:fused_loop:start}--\ref{alg:interior-penalty-single-level-one-mesh-sweep:fused_loop:finish} are executed $\nu+1$ times ($\nu$ times for the smoother and once to compute the residual) and lines \ref{alg:interior-penalty-single-level-one-mesh-sweep:projection-preamble-start}--\ref{alg:interior-penalty-single-level-one-mesh-sweep:projection-preamble-end} are executed once at the beginning to compute the projections $\vec{u}^{(c)}$ from $\vec{u}^{(c)}$. In the final mesh traversal it is not necessary to execute lines \ref{alg:interior-penalty-single-level-one-mesh-sweep:projection:start}--\ref{alg:interior-penalty-single-level-one-mesh-sweep:projection:finish} since the projections $\vec{u}^{(\pm)}$ are not required to restrict the residual. Analogous arguments apply for \cref{alg:interior-penalty-single-level-tasks}.

\paragraph{Prolongation}

When using \cref{alg:interior-penalty-single-level-multiple-sweeps-and-auxiliary-variables} for the fine-level block-Jacobi smoother, the prolongation of the coarse grid correction does not require an additional mesh traversal. Instead, this can be integrated into the loop in lines \ref{alg:interior-penalty-single-level-multiple-sweeps-and-auxiliary-variables:first-loop-block:start}--\ref{alg:interior-penalty-single-level-multiple-sweeps-and-auxiliary-variables:first-loop-block:end} for the first subsequent smoother application in the \textit{next} multigrid iteration: before projecting the solution to the facets in line \ref{alg:interior-penalty-single-level-multiple-sweeps-and-auxiliary-variables:projection} of \cref{alg:interior-penalty-single-level-multiple-sweeps-and-auxiliary-variables}, we prolongate $\vec{e}^{(\text{coarse})}$ according to \eqref{equation:multigrid:prolongation} in each cell $K$ to obtain the coarse grid correction $\delta \vec{u}^{(c)}|_K$, which is added to the current solution $\vec{u}^{(c)}$. In the final multigrid iteration, in which the prolongation is not followed by another block-Jacobi step, an additional mesh traversal is required to prolongate the coarse grid solution. We conclude that a total of $3n_{\text{iter}}(\nu+1)+1$ mesh traversals is required to perform $n_{\text{iter}}$ multigrid iterations with \cref{alg:multiplicative_multigrid} if the fine level smoother implementation is based on \cref{alg:interior-penalty-single-level-multiple-sweeps-and-auxiliary-variables}.

Loop fusion of the prolongation step can be applied in a very similar way when \cref{alg:interior-penalty-single-level-one-mesh-sweep} is used instead of \cref{alg:interior-penalty-single-level-multiple-sweeps-and-auxiliary-variables}: the prolongation of the coarse grid solution can be combined with the projection in lines \ref{alg:interior-penalty-single-level-one-mesh-sweep:projection-preamble-start}--\ref{alg:interior-penalty-single-level-one-mesh-sweep:projection-preamble-end} of \cref{alg:interior-penalty-single-level-one-mesh-sweep} in the first block-smoother application of the \textit{next} multigrid iteration: in each cell $K$ the correction $\delta \vec{u}^{(c)}|_K$ is computed from $\vec{e}^{(\text{coarse})}$ according to \eqref{equation:multigrid:prolongation} and added to the current solution $\vec{u}^{(c)}$ before computing $\vec{u}^{(\pm)}$. Again, the final multigrid iteration, which is not followed by another smoothing step, needs to be treated differently: here \eqref{equation:multigrid:prolongation} has to be executed in every cell $K$ in a separate mesh traversal. Altogether this results in $n_{\text{iter}}(\nu+2)+1$ mesh traversals if \cref{alg:multiplicative_multigrid} is implemented with \cref{alg:interior-penalty-single-level-one-mesh-sweep} (or \cref{alg:interior-penalty-single-level-tasks}).


\paragraph{Exit criterion}

An iterative solver such as the multigrid iteration, which computes the new iterate $\vec{u}^{(c)} := (\vec{u}^{(c)})^{k+1}$ from the previous iterate $\vec{u}^{(c)}_{\text{old}} := (\vec{u}^{(c)})^{k}$, is usually subject to a dynamic exit criterion: instead of performing a fixed number of steps, the iteration is terminated once the approximate solution $(\vec{u}^{(c)})^k$ is sufficiently close to the true solution $\vec{u}^{(c)}_{\text{true}}$. Unfortunately, since we do not know $\vec{u}^{(c)}_{\text{true}}$, it is not possible to compute the norm of the error $\vec{e}^{(c)} = \vec{u}^{(c)}-\vec{u}^{(c)}_{\text{true}}$ directly. Although often used in practice, the norm of the residual $\vec{r}^{(c)} = \vec{b}^{(c)} - A\vec{u}^{(c)} = A\vec{e}^{(c)}$ leads to a poor termination criterion if the matrix $A$ is ill-conditioned: a small values of the residual norm $\|\vec{r}^{(c)}\|$ does not necessarily imply the smallness of the error itself. This is the case for the interior penalty discretisation of the Poisson equation that used here and the same applies for many other problems of practical interest. 

A better exit criterion is the \emph{preconditioned} residual $\vec{r}^{(c)}_{\text{prec}}$, which for \cref{alg:multiplicative_multigrid} is given by the difference between the solutions at two subsequent iterations :
\begin{equation}
  \vec{r}^{(c)}_{\text{prec}} = \vec{u}^{(c)}-\vec{u}^{(c)}_{\text{old}}
   =  \mathcal{P}^{-1} \vec{r}^{(c)}
   = \mathcal{P}^{-1} A (\vec{u}^{(c)}_{\text{true}}-\vec{u}^{(c)}_{\text{old}}).
  \label{eqn:preconditioned_residual}  
\end{equation}
Here $\mathcal{P}^{-1}$ stands for one application of the preconditioner that corresponds to a single multigrid cycle. Observe in particular that $\vec{r}^{(c)}_{\text{prec}}$ is $\mathcal{P}^{-1} A$ times the error $\vec{u}^{(c)}_{\text{true}}-\vec{u}^{(c)}_{\text{old}}$. For a good preconditioner such as multigrid the matrix $\mathcal{P}^{-1} A$ is well-conditioned and hence the preconditioned residual is a good proxy for the error itself.
Maintaining it requires us to introduce an additional volumetric field $\vec{u}_{\text{old}}^{(c)}$. This introduces a memory overhead (\cref{def:memory_overhead}). Fortunately, the data is written and read only once per multigrid cycle, i.e.~not per mesh traversal or smoothing step (\cref{def:weak_single_touch}).

\subsection{Extension to $\boldsymbol{hp}$- multigrid}\label{sec:hp_multigrid}

It remains to find an efficient solver for the coarse grid equation $A^{(\text{coarse})}\vec{e}^{(\text{coarse})}=\vec{b}^{(\text{coarse})}$ in \cref{alg:coarse-grid correction} line \ref{alg:coarse-grid correction:coarse-grid-solve}.
In our purely geometric multigrid approach with rediscretisation, the matrix $A^{(\text{coarse})}$ arises from a piecewise linear finite element discretisation of the Poisson equation \eqref{eqn:introduction:poisson}.
As a direct solve of this problem might still be excessively expensive, some papers use an algebraic multigrid (e.g.~\cite{bastian2019matrix,bastian2012algebraic}).
For the problem at hand, geometric multigrid methods \cite{reusken2008introduction,hackbusch2013multi} are significantly simpler and show comparable performance. 
Extensions to geometric-algebraic approaches which preserve the geometric nesting of the function spaces while employing more complex algebraic operators are known \cite{weinzierl2018quasi} but have not been used for the present work.

Since the mesh is constructed through a spacetree based upon three-partitioning, it is coarsened recursively by combining blocks of $3^d$ grid cells (\cref{figure:discretisation:spacetree}). 
This induces a hierarchy of nested continuous Galerkin (CG) function spaces $\VCG{h}{1}\supset\VCG{3h}{1}\supset\dots$. 
On each level, the solution is smoothed with a simple point-Jacobi method before restricting the residual to the next coarser grid. There the algorithm is applied recursively to solve the coarse grid equation, before prolongating the solution back to the next-finer level and applying a small number of post-smoothing steps. 

The $h$-multigrid algorithm for the correction problem is highly efficient since it reduces the error on all length scales. Hence, usually only a single V-cycle is applied to obtain an approximate solution of the coarse grid equation in \cref{alg:coarse-grid correction}. 
The overall algorithm is classic $hp$-multigrid since it combines $p$-coarsening in the polynomial degree ($\VDG{h}{p}\rightarrow \VCG{h}{1}$) with $h$-coarsening of the grid and associated function spaces ($\VCG{h}{1}\rightarrow\VCG{3h}{1}$) on the coarser levels.

\section{Numerical results}
\label{section:results}

\begin{figure}
 \begin{center}
  \includegraphics[width=0.9\textwidth]{\figdir/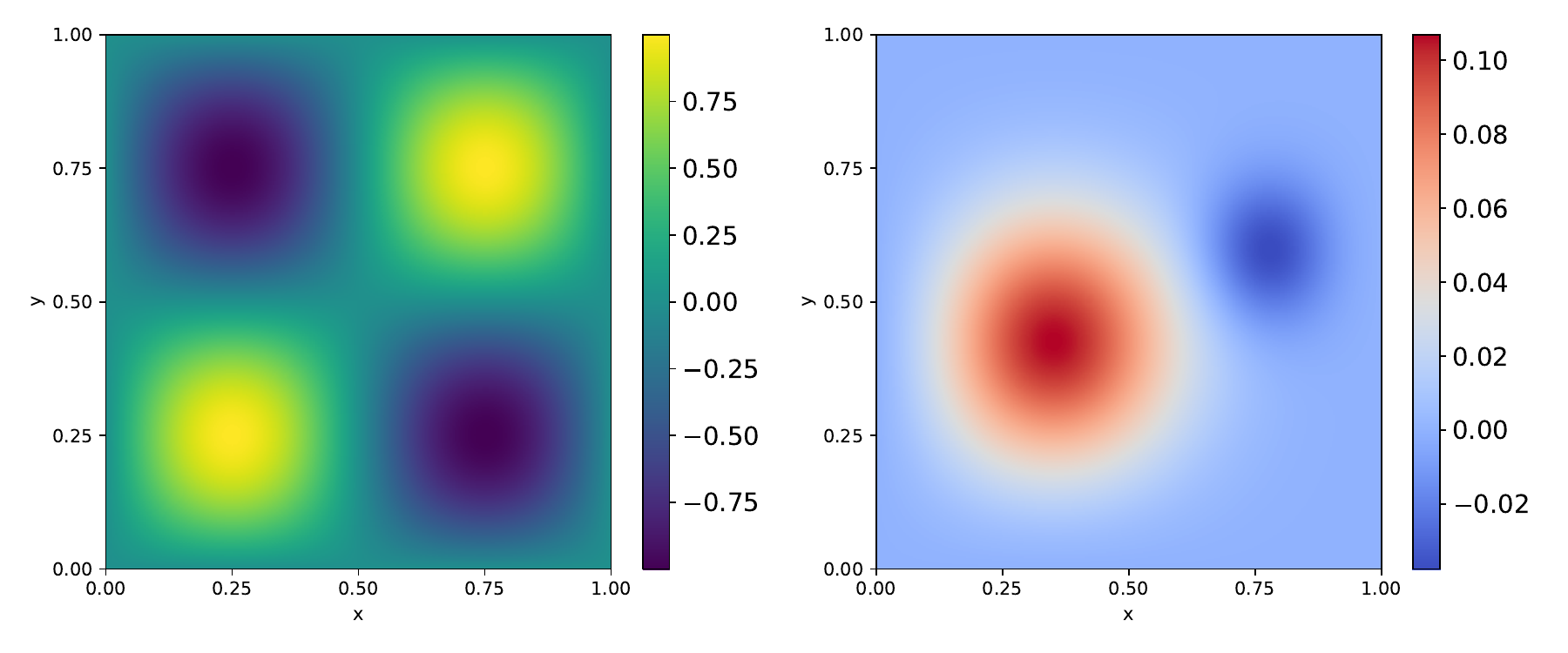}
  \caption{
    Visualisation of ``sin-product'' reference solution $u_1^{\textrm{ref.}}$ (left) as defined in \eqref{eq:results:sin_product_exact} and ``two-peak'' reference solution $u_2^{\textrm{ref.}}$ (right) as defined in \eqref{eq:results:two_peak_exact}.
    \label{fig:exact_solutions}
   }
 \end{center}      
\end{figure}

We consider the Poisson equation \eqref{eqn:introduction:poisson} with homogeneous Dirichlet boundary conditions on the unit square $\Omega = [0,1]\times[0,1]$ for two setups with manufactured analytical solutions:
\begin{subequations}
	\begin{align}
  u_1^{\textrm{ref.}} (x,y) & =  \sin(2 \pi x) \sin(2 \pi y)
    \qquad \text{or}
	\label{eq:results:sin_product_exact}
	\\
  u_2^{\textrm{ref.}} (x,y) & =  x(1-x)y(1-y) 
	\left( 
	2 \exp\left(-\frac{(x - x_{01})^2 + (y - y_{01})^2}{2\sigma_1^2}\right) -
	\exp\left(-\frac{(x - x_{02})^2 + (y - y_{02})^2}{2\sigma_2^2}\right) 
	\right)
	\label{eq:results:two_peak_exact}.
\end{align}
\end{subequations} 
The values of the parameters are set to $x_{01} = 0.3$, $y_{01} = 0.4$, $\sigma_1 = 0.2$, $x_{02} = 0.8$, $y_{02} = 0.6$, $\sigma_2 = 0.1$. Analytical expressions for the corresponding right-hand sides $f_i = - \Delta u_i^{\textrm{ref.}}$ are obtained by applying the Laplacian to the expressions in \eqref{eq:results:sin_product_exact} and \eqref{eq:results:two_peak_exact}. In what follows, we will also refer to $u_1^{\textrm{ref.}}$ as the ``sin-product'' and $u_2^{\textrm{ref.}}$ as the ``two-peak'' reference solution (\cref{fig:exact_solutions}).

\subsection{Discretisation error and mesh convergence}
Let $\vec{u}^{\textrm{ref.}}_{h,p}$ be the vector of unknowns that is obtained by interpolating the exact solution in \eqref{eq:results:sin_product_exact} or \eqref{eq:results:two_peak_exact} onto $\VDG{h}{p}$. The dof-vector of the corresponding numerical solution of the discretised equation \eqref{eqn:matrix_form_IP} for a given grid spacing $h$ and polynomial degree $p$ in the Gauss-Lobatto basis is denoted by $\vec{u}_{h,p}$. We compute $\vec{u}_{h,p}$ with the two-grid DG solver in \cref{alg:multiplicative_multigrid} and use \eqref{eqn:preconditioned_residual} to converge to a tolerance $\vec{r}^{(c)}_{\text{prec}} \leq \epsilon = 10^{-10}$ on the relative (preconditioned) residual to ensure that the error induced by the iterative solver is negligible. The relative discretisation error norm can be defined as

\[
	E_{h,p} = \| \vec{u}_{h,p} - \vec{u}^{\textrm{ref.}}_{h,p} \| / \| \vec{u}^{\textrm{ref.}}_{h,p} \|,
\]

\noindent
where $\|\cdot\|$ denotes either the $\ell_2$ or the $\ell_\infty$ norm defined by $\|\vec{x}\|_2 = (\sum_{j} x_j)^{1/2}$ or $\|\vec{x}\|_\infty = \max_j |x_j|$, respectively.

\begin{figure}
	\begin{center}
		\includegraphics[width=0.48\textwidth]{\figdir/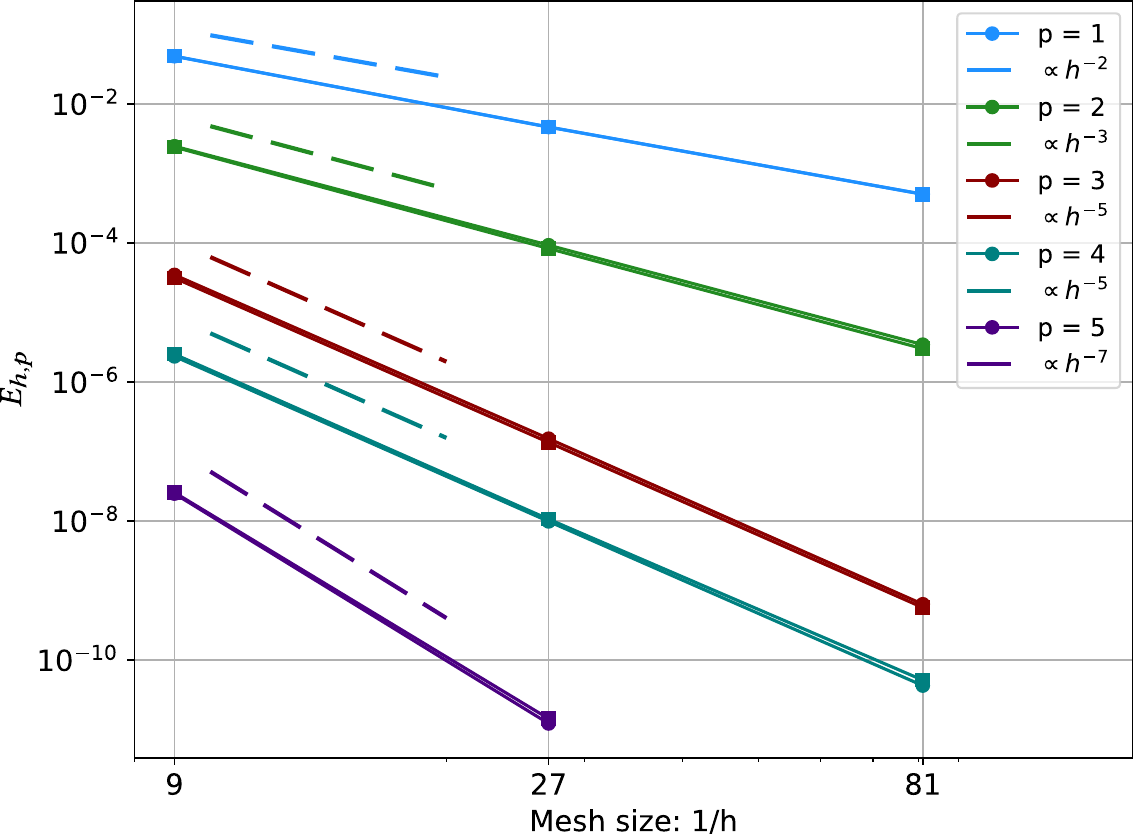}
		\hspace{1em}
		\includegraphics[width=0.48\textwidth]{\figdir/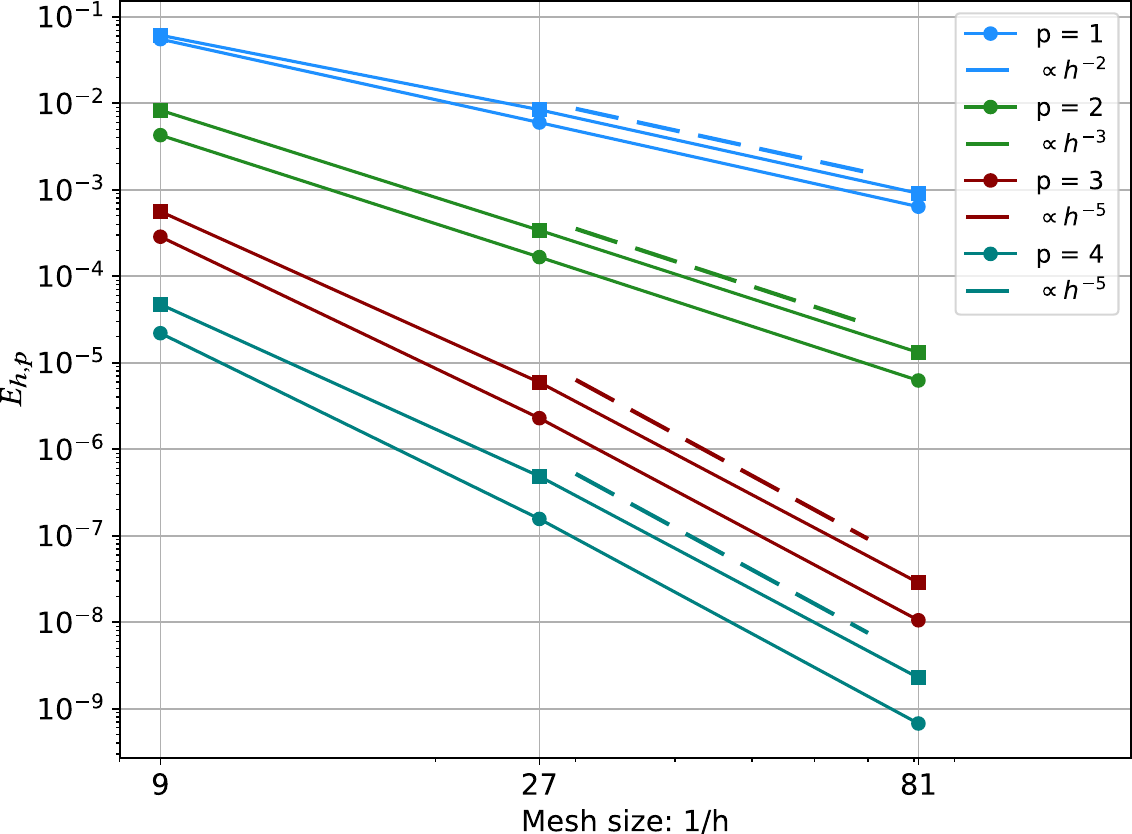}
	\end{center}
	\caption{
		Error $E_{h,p}$ in $\ell_2$ (circles $\bullet$) and $\ell_\infty$ (squares \scalebox{0.6}{$\blacksquare$}) norms for various choices of $p$ and $h$. Results are shown for the ``sin-product'' reference solution $u_1^{\textrm{ref.}}$ in \eqref{eq:results:sin_product_exact} (left) and the ``two-peak'' reference solution $u_2^{\textrm{ref.}}$ in \eqref{eq:results:two_peak_exact} (right).
		\label{fig:results:discretisation_error}
	}
\end{figure}

Empirically, the error decreases with $E_{h,p}\propto h^{p+1}$ (\cref{fig:results:discretisation_error}) for both manufactured analytical solutions in \eqref{eq:results:sin_product_exact} and \eqref{eq:results:two_peak_exact}. 
This exponential dependence of the error on the polynomial degree $p$ makes the interior penalty discretisation computationally efficient.
Compared to low order methods, significantly fewer unknowns are required to reduce the error below a given threshold. The results also confirm that with the given tolerance on the preconditioned residual, the error introduced by the iterative solver is indeed negligible compared to the discretisation error.

\subsection{Comparison of different solver variants}

Next, we explore the numerical efficiency of the different solver algorithms from \cref{section:efficiency} and \cref{section:multigrid}. For this, we consider the following configurations:
\begin{enumerate} 
	\item The standalone single-level DG solver (\cref{alg:interior-penalty-single-level-one-mesh-sweep}), and
	\item the two-grid algorithm (\cref{alg:multiplicative_multigrid}), where two strategies are used to compute the correction in the CG subspace:
	\begin{enumerate}
		\item Solve the CG equation up to a relative tolerance of $10^{-14}$ in the $\ell_\infty$ norm with repeated applications of $h$-multigrid V-cycles (since this tolerance is comparable to the truncation error in double precision floating point arithmetics the coarse level solver can be considered to be exact), or
		\item apply a single $h$-multigrid V-cycle to approximate the solution to the coarse level correction.
	\end{enumerate}
\end{enumerate}

\noindent
To compare the efficiency of the different iterative solvers, we track the evolution of the normalised residual norm. 

\begin{figure}
	\begin{center}
		\includegraphics[width=0.49\textwidth]{\figdir/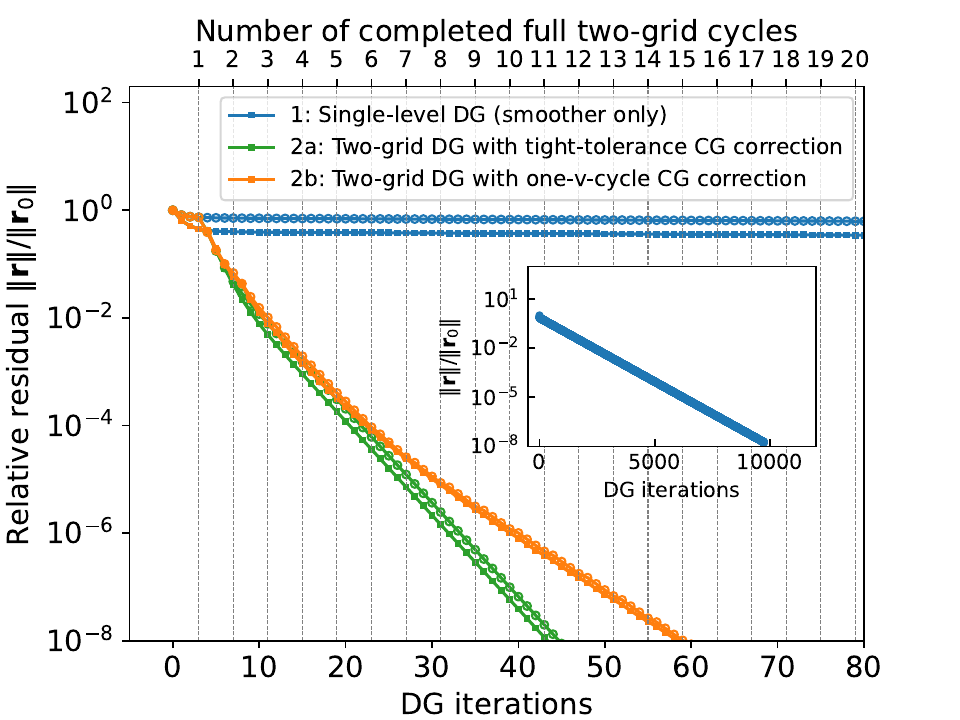}
		\includegraphics[width=0.49\textwidth]{\figdir/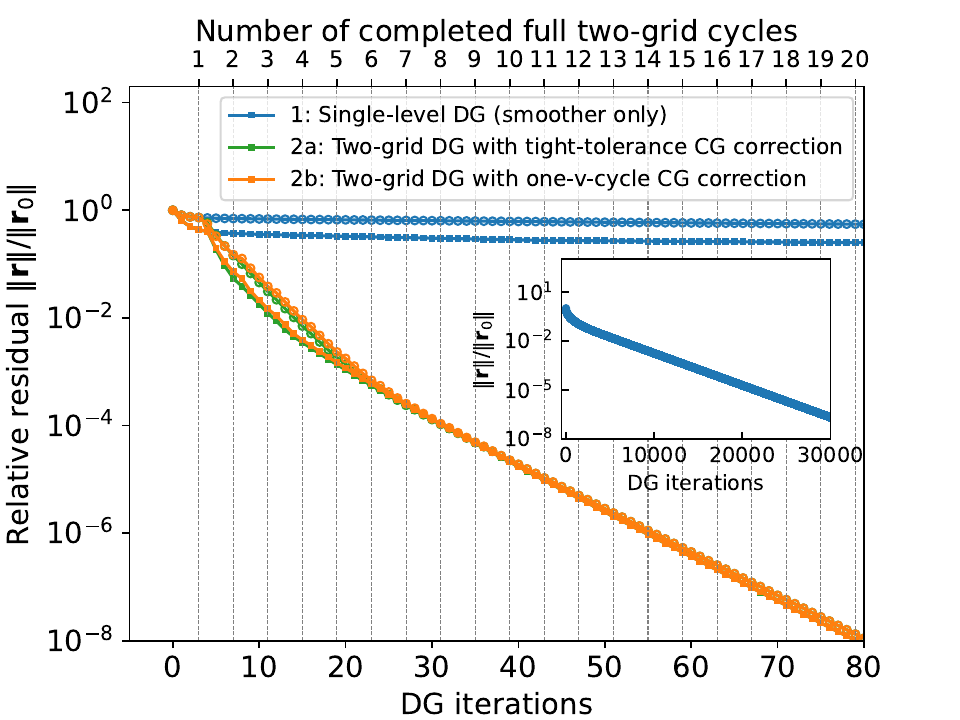}
	\end{center}
	\caption{
		 Evolution of the relative residual in the $\ell_{\infty}$ (squares) and $\ell_2$- (circles) norms computed for the ``sin-product'' reference solution $u_1^{\textrm{ref.}}$ in \eqref{eq:results:sin_product_exact} (left) and the ``two-peak'' reference solution $u_2^{\textrm{ref.}}$ in \eqref{eq:results:two_peak_exact} (right). In both cases we choose $p = 2$ on a mesh with $27\times27$ cells. Vertical dashed lines separate full two-grid cycles.
		 The inset plots show the zoomed-out curve for the single-level case (solver ``1'') in $\ell_2$-norm.
		\label{fig:results:residual_evolution}
	}
\end{figure}

For the single level DG smoother convergence is extremely slow (\cref{fig:results:residual_evolution}): thousands of iterations are required to reduce the residual norm below a reasonably small tolerance. 
We conclude that using the single level block-Jacobi iteration in \cref{alg:interior-penalty-single-level-one-mesh-sweep} as a standalone solver is inadequate for any practical applications. In contrast, the multigrid algorithm converges rapidly, reducing the relative residual by one order of magnitude for every 2-3 two-grid cycles. As expected, convergence is fastest if the coarse level equation $A^{(\text{coarse})}\vec{e}^{(\text{coarse})}=\vec{r}^{(\text{coarse})}$ is solved exactly. 
However, this advantage, compared to the approximate solve of the coarse level equation with a single geometric multigrid V-cycle, is negligible for the more realistic ``two-peak'' setup in \eqref{eq:results:two_peak_exact}. Even for the ``sin-product'' setup in \eqref{eq:results:sin_product_exact}, the reduction in the number of iterations does not justify the significantly higher cost of the ``exact'' CG solve.  In the following, we therefore always use a single $h$-multigrid V-cycle to approximately solve the error correction equation in the coarse CG space.

\subsection{Robustness of the multigrid solver}\label{sec:multigrid_robustness}
To assess the robustness of the $hp$-multigrid with respect to changes in the grid-resolution $h$ and polynomial degree $p$, we consider the evolution of both the unpreconditioned residual  $\vec{r}^{(c)} = \vec{u}^{(c)} - A\vec{u}^{(c)}$ and the preconditioned residual $\vec{r}^{(c)}_{\text{prec}}$ defined in \eqref{eqn:preconditioned_residual}. The initial residual values, which are used to normalise the relative residuals, are computed at different stages of the algorithm for the two cases: for the unpreconditioned residual, we use the value $\vec{r}_0^{(c)} = \vec{b}^{(c)} - A \vec{u}^{(c)}_0$, where $\vec{u}_0^{(c)}$ is the initial guess of the solution, computed before the first cycle starts. 
For the preconditioned case, the first residual vector is computed at the end of the first two-grid cycle as $\vec{r}^{(c)}_{\mathrm{prec,1}} = \vec{u}^{(c)}_1 - \vec{u}^{(c)}_0$. 

While initially the preconditioned residual norm decreases faster than the unpreconditioned residual norm, the asymptotic convergence rates are comparable.

\begin{figure}
	\begin{center}
		\includegraphics[width=0.49\textwidth]{\figdir/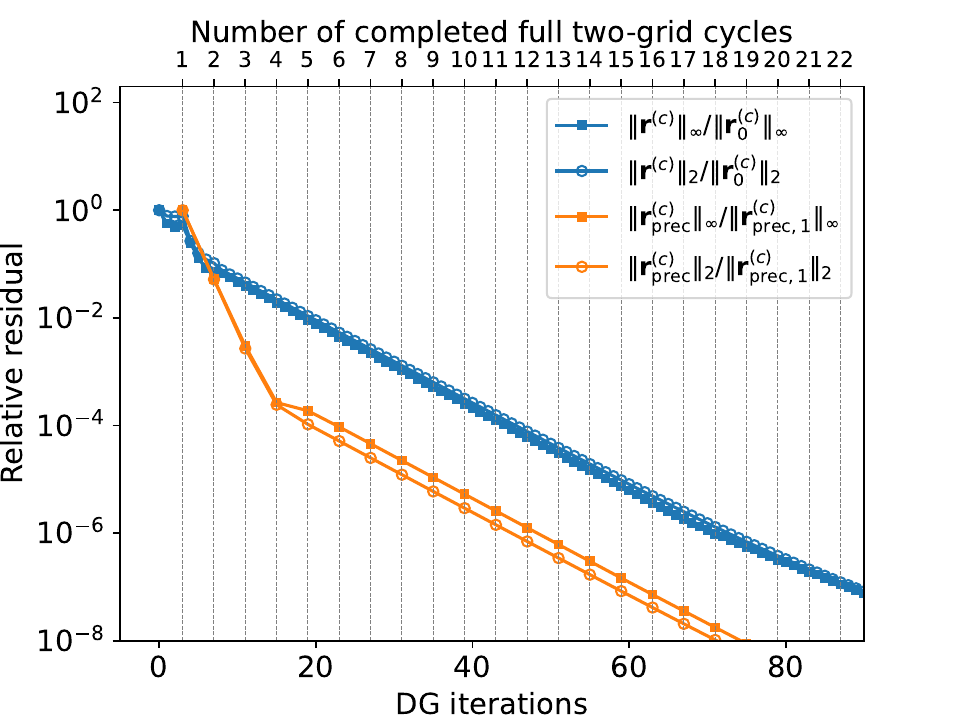}
		\includegraphics[width=0.49\textwidth]{\figdir/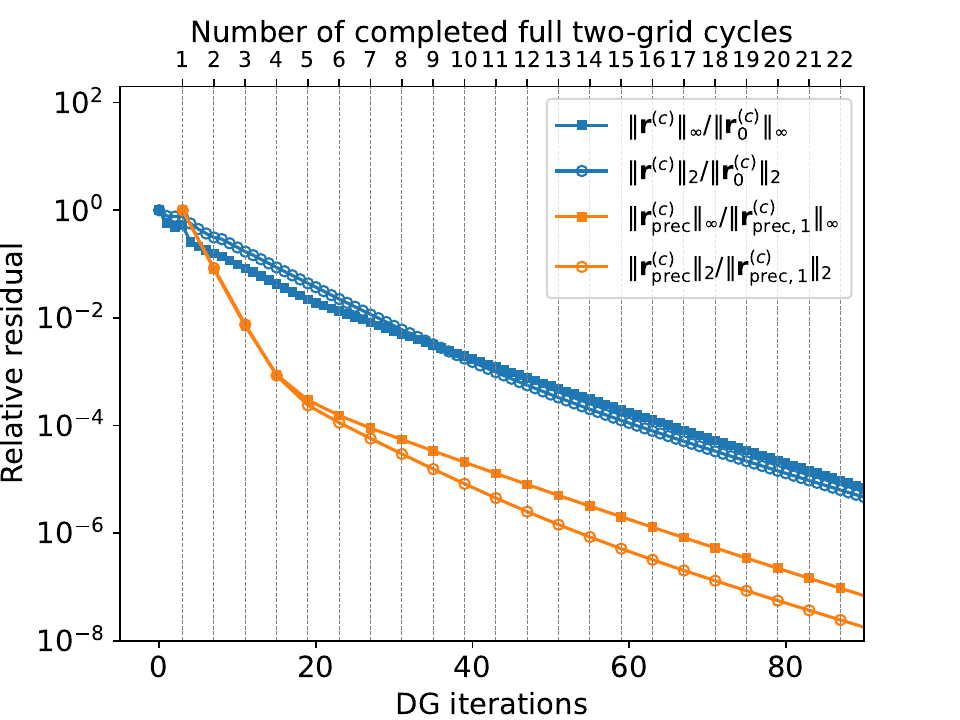}
	\end{center}
	\caption{
		Evolution of the preconditioned and unpreconditioned residuals. The relative residual is shown in the $\ell_\infty$- (squares) and $\ell_2$- (circles) norm computed for the ``sin-product'' reference solution $u_1^{\textrm{ref.}}$ in \eqref{eq:results:sin_product_exact} (left) and the ``two-peak'' reference solution $u_2^{\textrm{ref.}}$ in \eqref{eq:results:two_peak_exact} (right). In both cases we choose $p = 3$ on a mesh with $27\times27$ cells and only consider solver configuration 2(b), i.e. the $hp$-multigrid algorithm. Vertical dashed lines separate subsequent two-grid cycles.
		\label{fig:results:prec_residual_evolution}
	}
\end{figure}

The choice of exit criterion has an impact on the error, i.e. the norm of the difference between the approximate solution $\widehat{\vec{u}}^{(c)}$ obtained with the solver and the exact solution $\vec{u}^{(c)}_{\text{exact}}$ of $A\vec{u}^{(c)}=\vec{b}^{(c)}$. One would expect that for a given $\epsilon$ the norm of this error can be bounded by a constant times the tolerance: $\|\widehat{\vec{u}}^{(c)}-\vec{u}^{(c)}_{\text{exact}}\|/\|\vec{u}^{(c)}_{\text{exact}}\| < C\cdot \epsilon$. As shown in \cref{sec:precond_vs_unprecond}, if the exit criterion is based on the preconditioned residual, this estimate is robust in the sense that empirically $C$ does not depend on the resolution or polynomial degree. Therefore $\epsilon$ can be considered as a robust proxy for the size of the error. In contrast, for the unpreconditioned residual $C=C(h,p)$ depends on both the grid spacing $h$ and the degree $p$. Nevertheless, in the literature it is common to use the unpreconditioned residual. This has the advantage that the solver does not require additional storage and memory accesses (\cref{def:memory_overhead}) to evaluate the expression in \eqref{eqn:preconditioned_residual}. However, in this case care has to be taken in assessing the quality of the numerical solution for varying $h$ and $p$ if $\epsilon$ is kept fixed. In fact, \cref{fig:results:prec_residual_evolution} and \cref{fig:er_vs_ez} both suggest that the use of an unpreconditioned residual will result in an unnecessary increase in the number of iterations, which will offset any performance gain achieved by not storing the previous iterate required in \eqref{eqn:preconditioned_residual}.
%


\begin{table}[h!]
	\caption{Number of $hp$-multigrid cycles required to reduce the relative $\ell_2$-norm of the unpreconditioned residual $\vec{r}^{(c)}=\vec{b}^{(c)}-A\vec{u}^{(c)}$ by a factor of $10^{-7}$. Results are shown both for the ``sin-product'' reference solution $u_1^{\text{ref.}}$ in \eqref{eq:results:sin_product_exact} and the ``two-peak'' reference solution in \eqref{eq:results:two_peak_exact}.}\label{tab:results:residual_unprec}				
	\begin{center}		
		\begin{tabular}{ccccccccccccc}\hline
			& & \multicolumn{5}{c}{``sin-product''} &
			& \multicolumn{5}{c}{``two-peak''}\\
			degree & & $2$ & $3$ & $4$ & $5$ & $6$ &
				   & $2$ & $3$ & $4$ & $5$ & $6$
			\\ \hline
			$9\times9$ & & 12 & 24 & 43 & 63 & 89 &
			& 19 & 36 & 59 & 89 & 125
			\\ 
			$27\times27$ & & 13 & 23 & 41 & 61 & 86 & 
			& 18 & 33 & 55 & 82 & 116
			\\ 
			$81\times81$ & & 13 & 22 & 41 & 61 & 85 &
			& 17 & 31 & 52 & 78 & 111
			\\ 
			$243\times243$ & & 13 & 22 & 41 & 61 & 85 &
			& 16 & 30 & 50 & 75 & 106 			
			\\ \hline
		\end{tabular}
	\end{center}
\end{table}

\begin{table}[h!]
	\caption{Number of $hp$-multigrid cycles required to reduce the relative $\ell_2$-norm of the preconditioned residual $\vec{r}^{(c)}_{\text{prec}}$ in \eqref{eqn:preconditioned_residual} by a factor of $10^{-7}$. Results are shown both for the ``sin-product'' reference solution $u_1^{\text{ref.}}$ in \eqref{eq:results:sin_product_exact} and the ``two-peak'' reference solution in \eqref{eq:results:two_peak_exact}.}\label{tab:results:residual_prec}
	\begin{center}
		\begin{tabular}{ccccccccccccc}\hline
			& & \multicolumn{5}{c}{``sin-product''} &
			& \multicolumn{5}{c}{``two-peak''}\\
			degree & & $2$ & $3$ & $4$ & $5$ & $6$ & & $2$ & $3$ & $4$ & $5$ & $6$
			\\ \hline
			$9\times9$ & & 11 & 20 & 32 & 46 & 62 &
			& 16 & 27 & 43 & 61 & 82
			\\ 
			$27\times27$ & & 9 & 15 & 25 & 35 & 47 &
			& 12 & 19 & 29 & 41 & 55
			\\ 
			$81\times81$ & & 7 & 12 & 19 & 26 & 35 &
			& 9 & 14 & 21 & 29 & 39
			\\ 
			$243\times243$ & & 7 & 9 & 13 & 18 & 22 &
			& 8 & 10 & 15 & 20 & 26
			\\ \hline
		\end{tabular}
	\end{center}
\end{table}

\cref{tab:results:residual_unprec} shows the number of multigrid cycles required to reduce the relative unpreconditioned residual norm by seven orders of magnitude for the two setups in ``sin-product'' setup in \eqref{eq:results:sin_product_exact} and for the ``two-peak'' setup in \eqref{eq:results:two_peak_exact}; the corresponding results for the preconditioned residual norm are shown in \cref{tab:results:residual_prec}. 

The method is $h$-robust.
The number of iterations even decreases on finer meshes, which might be a result of the boundary operator approximation in \cref{technique:approximate-boundary-operators}.
For coarse resolutions, an inaccurate treatment of the cells at the domain boundary has a relatively greater impact than for fine resolutions. These effects only show up in the preconditioned residual and require further in-depth studies. 

The method is not totally $p$-robust, i.e.~the number of iterations grows for higher polynomial degrees. 
This is to be expected since our non-overlapping block-Jacobi smoother is known to be not $p$-robust (\cref{technique:weaker-block-smoother}) if it is combined with the agressive $p$-coarsening to the lowest order CG function space. 
Instead of implementing a better smoother and no longer enforcing the desirable properties described in \cref{def:memory_overhead} or \cref{def:weak_single_touch}, one could reduce the polynomial degree gradually by constructing a nested sequence of spaces $\VDG{h}{p} \supset \VDG{h}{p-1} \supset \VDG{h}{p-2} \supset\dots\supset \VCG{h}{1}$ before transitioning into h-multigrid \cite{kronbichler2018performance}.

However, the non-overlapping cell-wise smoothers that we use here are beneficial from an HPC point-of-view, as they do not require non-local data accesses and complex synchronisation.
A gradual reduction of polynomial degree preserves this advantageous character of our implementation and should be the subject of future studies.

\subsection{Choice of nodal basis}\label{sec:numerical_GL_vs_GLL}

All previously discussed results were obtained with a nodal Gauss-Lobatto basis for both DG function spaces. In this case, those nodal points of $\VDG{h}{p}$ which are associated with the surface of a cell coincide with nodal points of $\FDG{h}{p}$ on the facets. This simplifies the projection (cmp.~\cref{observation:GaussLobatto-sparsity}). However, the choice of Gauss-Legendre basis functions results in a diagonal mass matrix, which can be advantageous in certain applications. Both choices of basis functions are used in the literature. As can be seen from \cref{tab:results:two_peak_unprec_gauss_legendre} (which should be compared to \cref{tab:results:residual_unprec} and \cref{tab:results:residual_prec}), both choices of basis functions result in very similar convergence behaviour.

\begin{table}[h!]
	\caption{Number of $hp$-multigrid cycles required to reduce the relative $\ell_2$-norm of the unpreconditioned residual $\vec{r}^{(c)}=\vec{b}^{(c)}-A\vec{u}^{(c)}$ (left) and preconditioned residual $\vec{r}^{(c)}_{\text{prec}}$ (right) by a factor of $10^{-7}$. Results are shown for the ``two-peak'' reference solution in \eqref{eq:results:two_peak_exact}.	In contrast to \cref{tab:results:residual_unprec} and \cref{tab:results:residual_prec}, Gauss-Legendre nodes are used to construct the DG basis functions.}
	\label{tab:results:two_peak_unprec_gauss_legendre}
	\begin{center}
			\begin{center}
				\begin{tabular}{cccccccccccc}\hline
					& \multicolumn{5}{c}{unpreconditioned $\vec{r}^{(c)}$} &
			& \multicolumn{5}{c}{preconditioned $\vec{r}^{(c)}_{\text{prec}}$}\\
					degree & $2$ & $3$ & $4$ & $5$ & $6$ & & $2$ & $3$ & $4$ & $5$ & $6$
					\\ \hline
					$9\times9$ & 20 & 37 & 62 & 92 & 131 & & 16 & 27 & 43 & 61 & 82\\ 
					$27\times27$ & 19 & 34 & 57 & 85 & 122 & & 12 & 19 & 29 & 41 & 55 \\ 
					$81\times81$ & 18 & 33 & 55 & 82 & 117 & & 9 & 13 & 21 & 29 & 38 \\ 
					$243\times243$ & 17 & 31 & 53 & 78 & --- & & 8 & 10 & 15 & 20 & ---\\ \hline
				\end{tabular}
			\end{center}
	\end{center}
\end{table}


\section{Performance evaluation}
\label{section:runtime}

The experiments to assess the computational efficiency of our algorithms use the AMD K17 (Zen2) architecture, namely a pair of AMD EPYC 7702 64-Core processors, where the 2$\times $64 cores per node are spread over two sockets. 
Each core has access to 32 kB exclusive L1 cache, and 512 kB L2 cache. 
The shared L3 cache is (physically) split into chunks of 16 MB associated with four cores each, while the internal setup of the chip gives each group of 16 cores access to two memory channels. This results in four NUMA domains per socket or eight per node. 
The code is compiled with Intel's oneAPI C++ Compiler \texttt{icpx} 2025.0.1. 
All  tasking relies exclusively on Threading Building Blocks (TBB) \cite{voss2019proTBB} subject to dynamic task graph extension as submitted to the 
\texttt{uxlfoundation/oneAPI-spec} repository under tag \texttt{ed26d0c}.


For the Poisson equation, which is solved in all benchmarks presented in this paper, the cell-local matrices are constant across the mesh.
While they need to be scaled by appropriate powers of the grid spacing to account for the changing resolution of different meshes in the multigrid hierarchy, the finite element matrices $A_{c\gets c}|_{K\gets K}$, $A_{f\gets c}|_{F\gets K}$ etc. can be assembled once on the unit reference cell and kept in cache for the entire run (\cref{observation:single-level:precompute-block-inverse}); the same applies to the block-diagonal matrix $A_{K\rightarrow K}$ whose inverse is applied to the residual in the block-Jacobi update \eqref{eqn:block-Jacobi-IP} to compute $\vec{u}^{(c)}|_K \gets \vec{u}^{(c)}|_K+\omega A_{K\gets K}^{-1}\vec{r}^{(c)}|_K$.  

For problems with an inhomogeneous, anisotric diffusion coefficient and/or on adaptively refined meshes the cell-local matrices vary across the domain. Pre-assembling and storing these matrices is undesireable since they will need to be loaded from memory and their storage requirements will limit the size of the problems that can be solved (see discussion in \cite[Section 4.4]{bastian2019matrix}); this issue is particularly pronouced for higher polynomial degrees. In a matrix-free implementation, the cell-local matrices hence have to be re-assembled on the fly. 
Unless techniques as in \cite{bastian2019matrix} (where the matrices $A_{K\gets K}$ are inverted iteratively) are used, this implies that the inverse $A_{K\gets K}^{-1}$ needs to be computed in every cell to obtain the increment $\omega A_{K\gets K}^{-1}\vec{r}^{(c)}|_K$ to the current iterate. It is reasonable to assume that the on-the-fly assembly of $A_{K\gets K}$ itself remains cheap compared to the computation of $A_{K\gets K}^{-1}$, especially for higher polynomial degrees $p$. 
If this is not the case, iterative integration can hide the assembly cost behind the solve \cite{murray2021assembly}.
While the benchmark considered here is homogeneous and isotropic, in the following we also extrapolate (``mimic'') the performance to such more challenging scenarios
by using \cref{obs:cost_of_local_matrix_inversion}: for this, we re-assemble the volumetric matrices and re-compute the inverse of the local matrix $A_{K\gets K}$ in each grid cell.

\subsection{Single-core performance of the DG block-Jacobi iteration}

\begin{table}
  \caption{Performance of the single-level DG solver on a single core of an AMD EPYC processor for a range of polynomial degrees $p$. Results are shown for both for the naive algorithm which involves multiple mesh sweeps (\cref{alg:interior-penalty-single-level-multiple-sweeps-and-auxiliary-variables}, top) and for the improved algorithm with loop fusion which only requires a single iteration of the computational grid (\cref{alg:interior-penalty-single-level-one-mesh-sweep}, bottom), as well as for variants where the inverse of the block-diagonal $A_{K\gets K}$ is pre-computed once at the start of the run (which is sufficient for solving the Poisson equation) vs.~variants where $A_{K\gets K}^{-1}$ is recomputed in each cell.
    \label{table:single-level-hw-counters-AMD-sweeps}
  }
  \begin{center}
    {\footnotesize

\begin{tabular}{rrrrrrr}
\hline\multicolumn{7}{c}{\cref{alg:interior-penalty-single-level-multiple-sweeps-and-auxiliary-variables}  (no loop fusion, multiple mesh sweeps), recompute $A_{K\gets K}^{-1}$ in each cell}\\\hline
\multicolumn{1}{c}{degree} & \multirow{2}{*}{t/dof [ns]} & \multirow{2}{*}{MFLOPs/s} & bandwidth & data volume & \multicolumn{2}{c}{L3 cache} \\
\multicolumn{1}{c}{$p$} &  &  & [MBytes/s] & [GBytes] & request rate & miss ratio \\\hline 
1 &  4236.23 &    64.72 &   580.80 &    97.65 &  0.44130 &  0.00280 \\ 
2 &  1984.52 &   191.82 &   542.36 &    95.14 &  0.44440 &  0.00550 \\ 
3 &  1267.10 &   580.85 &   596.40 &   112.48 &  0.45160 &  0.00920 \\ 
4 &   984.85 &  1485.79 &   825.48 &   175.89 &  0.48960 &  0.01520 \\ 
5 &  1761.03 &  1976.99 &   526.49 &   220.61 &  0.41320 &  0.01120 \\ 
6 &  2550.40 &  2532.70 &   176.08 &   130.79 &  0.43030 &  0.01730 \\ 
7 &  3981.48 &  2751.14 &   126.64 &   183.85 &  0.44750 &  0.03820 \\ 
8 &  6184.81 &  2802.77 &    86.46 &   238.80 &  0.45660 &  0.02670 \\ 
9 &  9388.70 &  2774.93 &    65.22 &   330.89 &  0.47530 &  0.02500 \\ 
\\
\hline\multicolumn{7}{c}{\cref{alg:interior-penalty-single-level-multiple-sweeps-and-auxiliary-variables}  (no loop fusion, multiple mesh sweeps), precompute $A_{K\gets K}^{-1}$ once}\\\hline
\multicolumn{1}{c}{degree} & \multirow{2}{*}{t/dof [ns]} & \multirow{2}{*}{MFLOPs/s} & bandwidth & data volume & \multicolumn{2}{c}{L3 cache} \\
\multicolumn{1}{c}{$p$} &  &  & [MBytes/s] & [GBytes] & request rate & miss ratio \\\hline 
1 &  4229.70 &    51.87 &   569.97 &    96.26 &  0.44170 &  0.00280 \\ 
2 &  1957.22 &   134.20 &   619.07 &   107.79 &  0.44330 &  0.00480 \\ 
3 &  1153.74 &   275.80 &   743.98 &   133.55 &  0.44880 &  0.00770 \\ 
4 &   769.74 &   548.90 &   662.20 &   120.73 &  0.45560 &  0.01340 \\ 
5 &   585.94 &   914.55 &   700.21 &   134.02 &  0.46260 &  0.01960 \\ 
6 &   493.60 &  1355.97 &   979.34 &   204.84 &  0.48340 &  0.02630 \\ 
7 &   401.95 &  2002.69 &   761.34 &   163.96 &  0.51450 &  0.03730 \\ 
8 &   385.50 &  2585.65 &  1312.47 &   324.51 &  0.53940 &  0.04240 \\ 
9 &   371.86 &  3254.02 &  2036.14 &   565.32 &  0.60210 &  0.04740 \\ 
\\
\hline\multicolumn{7}{c}{\cref{alg:interior-penalty-single-level-one-mesh-sweep} (loop fusion, single mesh sweep), recompute $A_{K\gets K}^{-1}$ in each cell}\\\hline
\multicolumn{1}{c}{degree} & \multirow{2}{*}{t/dof [ns]} & \multirow{2}{*}{MFLOPs/s} & bandwidth & data volume & \multicolumn{2}{c}{L3 cache} \\
\multicolumn{1}{c}{$p$} &  &  & [MBytes/s] & [GBytes] & request rate & miss ratio \\\hline 
1 &  4100.16 &    37.69 &   510.76 &   172.33 &  0.44130 &  0.00230 \\ 
2 &  1894.18 &   110.16 &   545.36 &   192.58 &  0.44330 &  0.00400 \\ 
3 &  1122.45 &   341.33 &   593.77 &   217.34 &  0.44760 &  0.00700 \\ 
4 &   781.29 &   914.95 &   540.89 &   211.99 &  0.46950 &  0.01170 \\ 
5 &   937.61 &  1515.47 &   334.34 &   206.00 &  0.41920 &  0.01070 \\ 
6 &  1157.29 &  2162.59 &   391.62 &   385.12 &  0.43240 &  0.01660 \\ 
7 &  1635.91 &  2544.20 &   184.32 &   323.15 &  0.44850 &  0.03590 \\ 
8 &  2417.49 &  2701.30 &   211.49 &   685.27 &  0.45540 &  0.02700 \\ 
9 &  3582.94 &  2728.87 &    81.98 &   475.26 &  0.47480 &  0.02520 \\ 
\\
\hline\multicolumn{7}{c}{\cref{alg:interior-penalty-single-level-one-mesh-sweep} (loop fusion, single mesh sweep), precompute $A_{K\gets K}^{-1}$ once}\\\hline
\multicolumn{1}{c}{degree} & \multirow{2}{*}{t/dof [ns]} & \multirow{2}{*}{MFLOPs/s} & bandwidth & data volume & \multicolumn{2}{c}{L3 cache} \\
\multicolumn{1}{c}{$p$} &  &  & [MBytes/s] & [GBytes] & request rate & miss ratio \\\hline 
1 &  4129.87 &    29.60 &   508.80 &   174.22 &  0.44140 &  0.00240 \\ 
2 &  1870.02 &    77.06 &   553.09 &   193.17 &  0.44290 &  0.00370 \\ 
3 &  1080.54 &   162.04 &   559.96 &   199.12 &  0.44590 &  0.00570 \\ 
4 &   696.46 &   332.64 &   536.35 &   192.35 &  0.44860 &  0.00990 \\ 
5 &   501.86 &   577.50 &   505.70 &   188.29 &  0.45510 &  0.01460 \\ 
6 &   394.66 &   898.62 &   794.17 &   308.80 &  0.46640 &  0.02040 \\ 
7 &   311.94 &  1359.85 &   960.25 &   380.96 &  0.48940 &  0.02770 \\ 
8 &   268.74 &  1896.97 &   973.02 &   421.97 &  0.51120 &  0.03390 \\ 
9 &   240.75 &  2512.91 &   463.85 &   283.78 &  0.55330 &  0.03930 \\ 
\\
\end{tabular}
    }
  \end{center}
\end{table}

We start with performance measurements on a regular two-dimensional grid with $729\times 729 = 531,441$ cells and evaluate the hardware performance counters on our system.
We focus on the DG block-Jacobi iteration in \cref{alg:interior-penalty-single-level-multiple-sweeps-and-auxiliary-variables} and \cref{alg:interior-penalty-single-level-one-mesh-sweep} since this is expected to be the bottleneck of the multigrid method in \cref{alg:multiplicative_multigrid}. 
 
In this section we exclusively study the performance on a single compute core. 
Because of this, we do not include the task-based \cref{alg:interior-penalty-single-level-tasks} in the comparison which intrinscially requires a multicore system.
 On our AMD hardware, the Stream TRIAD \cite{McCalpin1995,McCalpin2007} as shipped with \texttt{likwid-bench} \cite{treibig2012likwid} reports a memory bandwidth of 
1,883.56 MB/s per core on a fully populated node (resulting in a total bandwidth of 241,095.62 MB/s per node). 
This is equivalent to a cost of $ 3.28 \cdot 10^{-10}\text{s} \leq t_{\text{mem}} \leq 4.25 \cdot 10^{-9}\text{s}$ per double precision number transferred through the whole memory subsystem. 
The same benchmark reports 8,721.87 MFlops/s on a single core if we exclusively employ scalar operations.
However, once vectorisation with FMA and AVX512 kicks in, this increases to 45,261.50 MFlops/s.
This is equivalent to a cost of $2.21 \cdot 10^{-11}s \leq t_{\text{flop}} \leq 1.15 \cdot 10^{-10}s$ per floating point operation in double precision arithmetic. We assume that these two ``corridors'' that bound the cost $t_{\text{mem}}$ of one floating point operation and the cost $t_{\text{mem}}$ of a memory access ultimately limit the performance of our implementation.

\cref{table:single-level-hw-counters-AMD-sweeps} shows measurements of key performance counters for implementations of \cref{alg:interior-penalty-single-level-multiple-sweeps-and-auxiliary-variables} and \cref{alg:interior-penalty-single-level-one-mesh-sweep}. In each case we report two sets of results: in the first case, the inverse of the matrix $A_{K\gets K}$ which is needed in the update $\vec{u}^{(c)}|_K \gets \vec{u}^{(c)}|_K + \omega A_{K\gets K}^{-1} \vec{r}^{(c)}|_K$, is calculated once at the beginning of the run. Hence, in each iteration we only have to compute the matrix-vector product $A_{K\gets K}^{-1} \vec{r}^{(c)}|_K$ with the BLAS \texttt{dgemv} routine. In the second setup the inversion of $A_{K\gets K}$ and the matrix-vector product $A_{K\gets K}^{-1} \vec{r}^{(c)}|_K$ are both performed in every cell. 
The cost of the LU-factorisation, which is required to invert $A_{K\gets K}$ grows with the third power of the matrix size, i.e. $(p+1)^{3d}$ in our case. The matrix-vector product with \texttt{dgemv} only incurs a cost of $\mathcal{O}((p+1)^{2d})$. We would therefore expect that the homogeneous setup, which does not require a matrix-inversion per grid cell, is significantly more efficient for higher polynomial degrees $p$.

\paragraph{Experimental data for naive implementation with multiple mesh sweeps.}
First, we consider the results for the block-Jacobi iteration from \cref{alg:interior-penalty-single-level-multiple-sweeps-and-auxiliary-variables}. 
The runtime for a single smoothing step grows with increasing polynomial degree $p$ for both setups. However, the relative cost per degree of freedom, i.e.~the time for a single solver iteration divided by the number of grid cells and the number $(p+1)^2$ of DG unknowns per cell, decreases as $p$ increases up to around $p\approx 4$. After that, the relative cost rises again if we perform the inversion of $A_{K\gets K}$ in each cell.
If we rely on the precomputed $A_{K\gets K}^{-1}$, and hence exclusively apply mat-vecs (BLAS \texttt{dgemv}), the relative cost continues to decrease.

The delivered MFLOPs/s grow with the polynomial degree $p$. Not very surprisingly, precomputing the inverse $A_{K\gets K}^{-1}$ once at the beginning of the simulation reduces the runtime significantly, while it slightly increases the total memory moved over the bus.
Almost every second instruction hits the L3 cache. However, only a minority of these hits cannot be served by L3 and hence lead to a data transfer. This ratio only slightly increases with $p$.

\paragraph{Experimental data for improved implementation with single mesh sweep.}

Shifting the operator evaluation and fusing the three mesh traversals into a single one in \cref{alg:interior-penalty-single-level-one-mesh-sweep} reduces the wallclock time but has no significant impact on the L3 access characteristics.
Its impact on the volume of data moved as well as the bandwidth is not immediately clear, but seems to depend on how we deal with the inversion of $A_{K\gets K}$: With an on-the-fly matrix inversion, bandwidth usage and memory transferred decrease, while the use of a precomputed matrix $A_{K\gets K}$ reverses this trend.

The decrease in runtime compared to \cref{alg:interior-penalty-single-level-multiple-sweeps-and-auxiliary-variables} is observable for low polynomial degrees and becomes significant for larger values of $p$.
It is more pronounced if the matrix $A_{K\gets K}^{-1}$ is precomputed, in which case \cref{alg:interior-penalty-single-level-one-mesh-sweep} results in a speedup of more than a factor of two. Interestingly, the MFlop rate does not increase in line with the savings of runtime, which is difficult to explain given that both realisations compute exactly the same operations:
For either variant, with precomputed operators or not, a reduction of runtime due to loop fusion should, as we perform exactly the same arithmetic operations, lead to a higher FLOPs rate.
This effect requires further investigration.

\paragraph{Discussion.}
Since we employ a higher-order DG method, we apply dense (stiffness) matrices per cell and hence get a better ALU usage as $p$ increases (cmp.~\cref{def:vectorisation}).
While the code can make more efficient use of the hardware, the cost per inversion of $A_{K\gets K}$ grows faster than the efficiency gains. We hence find that increasing $p$ beyond a certain value results in an increase of cost per degree of freedom, i.e.~the improved ALU usage (vectorisation) cannot compensate for the growing cost anymore. Obviously, the cheapest, i.e.~best solution is always to avoid any on-the-fly matrix inversion and multiply with a precomputed inverse of $A_{K\gets L}$. In this case, the relative cost per degree of freedom decreases continuously with $p$ as the code can make increasingly better use of the ALU.

However, this usually works only for simple, homogeneous problems such as the Poisson equation in \eqref{eqn:introduction:poisson}. If the problem is locally homogenous,~i.e. the parameters are constant in parts of the domain, it is reasonable to consider multiple precomputed local matrix inversions or approximate inverses. To which degree such an approach remains stable---it is expected that the smoother and multigrid algorithm is less efficient, cmp.~\cite[Example 4.1]{bastian2012algebraic}---has to be subject of further investigations.

Our code base is written in a cache-oblivious way \cite{weinzierl2011peano} and hence manages to keep most face data in the L3 cache.
For \cref{alg:interior-penalty-single-level-one-mesh-sweep}, the additional storage required due to \cref{technique:left-right-face-projections} and \cref{technique:numerical-flux-outcome} hence does not manifest in increased memory stress
and the loop fusion exclusively induces volumetric data transfers, which validates the distinction between \cref{def:single_touch} and \cref{def:weak_single_touch}.
The shift-and-fuse approach of \cref{technique:shift-and-fuse} therefore pays off robustly.
We crossvalidated data from the strongly hierarchical EPYC architecture to an Intel Xeon Platinum 8480 (Sapphire Rapid) CPU that features two sockets with 56 cores per socket.
Its total L3 cache offers 105 MB per core, and the total memory of $2 \cdot 256$GB is split over only two NUMA domains.
No qualitiatively different results are obtained (not shown).

It is not clear why even the per-cell matrix inversion variants run only at around 50\% of the scalar peak of the core or 10\% of the theoretical peak performance.
We assume that this is due to the cell-face and face-face operators that slot into the sequence of expensive calculations.
These are intrinsic to the use of DG and the use of temporary data due to \cref{technique:left-right-face-projections} and \cref{technique:numerical-flux-outcome}.

\subsection{Domain decomposition with perfect balancing}
\label{subsection:results:domain-decomposition-perfect}

%
%
We next study the strong scalability of the different block-Jacobi implementations subject to non-overlapping domain decomposition in a (logically) distributed memory model. For this, the mesh is split along the Peano space-filling curve \cite{weinzierl2011peano} into chunks of approximately equal size per core: the number of cells in the resulting subdomains differs by at most one. In the strong scaling setting the mesh (and hence the problem size) remains fixed, while we increase the number of subdomains until we eventually reach the full core count; this process is repeated for different polynomial degrees $p$. 

\begin{figure}[htb]
  \begin{center}
    \includegraphics[width=0.48\textwidth]{\figdir/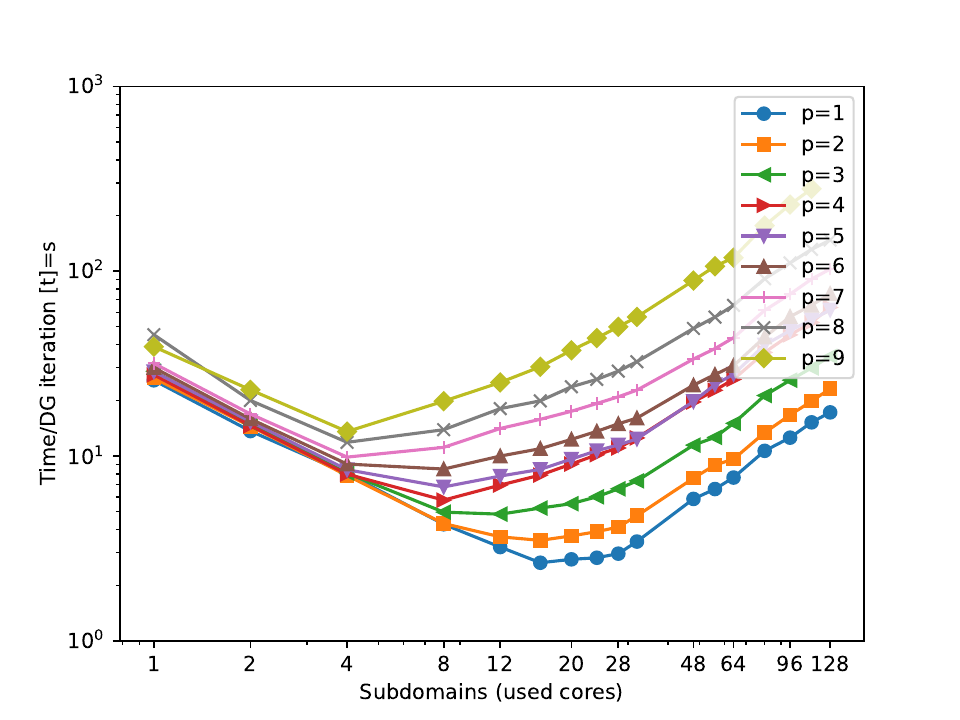}
    \includegraphics[width=0.48\textwidth]{\figdir/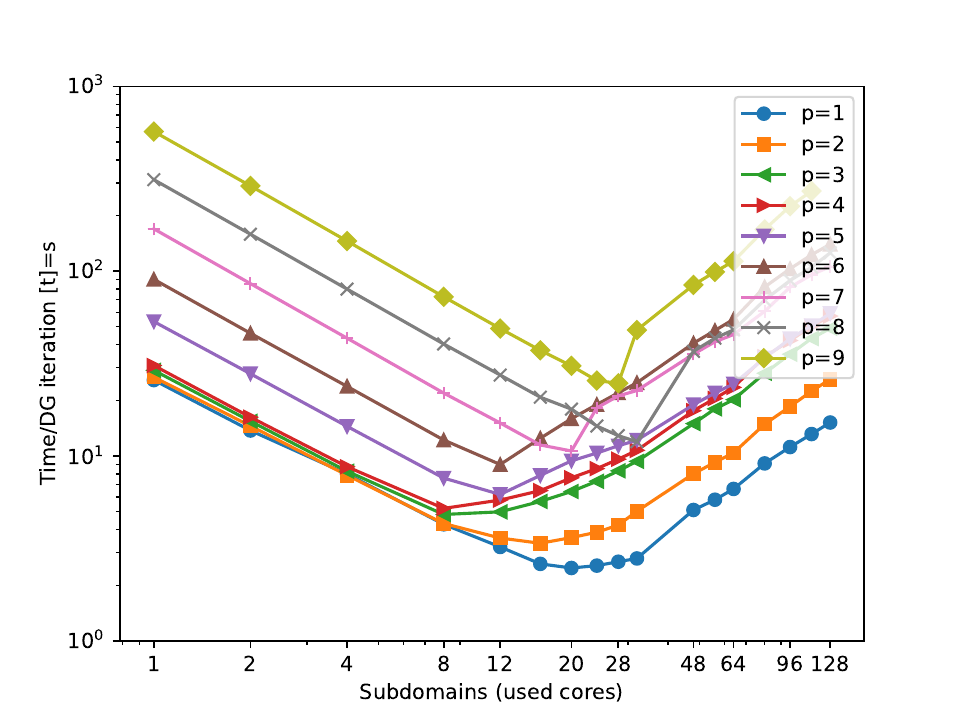}
  \end{center}
  \caption{
    Scalability with precomputed $A_{K\gets K}^{-1}$ (left) and inversion of $A_{K\gets K}$ in each grid cell (right) on an AMD EPYC node. Each smoothing step is mapped onto three mesh sweeps (\cref{alg:interior-penalty-single-level-multiple-sweeps-and-auxiliary-variables}).
    \label{figure:single-level-scalability:sweeps}
  }
\end{figure}

\paragraph{Experimental data for naive implementation with multiple mesh sweeps.}
\cref{figure:single-level-scalability:sweeps} shows the change in runtime as the number of cores increases for fixed problem size. Results are shown for a range of polynomial degrees $p$, and we consider both the case where the matrix $A_{K\gets K}$ is computed and inverted once at the beginning of the run (left) and the setup where $A_{K\gets K}^{-1}$ is re-computed in each cell of the mesh (right).

Our code exhibits ideal scaling, as long as we use a moderate number of subpartitions.
Once we increase the number of subpartitions beyond a certain threshold, the runtime increases again. 
The deterioration of scaling on higher core counts is caused by the fact that communication costs for exchanging the solution between subdomains eventually become dominant.
Due to \cref{technique:left-right-face-projections}, we still have to keep the $\vec{u}^{(\pm)}|_F$ values along the domain boundares consistent (\cref{subsubsection:single-level:domain-decomposition}).

If the matrix $A_{K\gets K}^{-1}$ is precomputed once at the beginning of the run, strong scaling breaks down significantly earlier. However, in the region where the code scales well, the runtime is less dependent on the polynomial degree. Both results are not very surprising since the inversion of $A_{K\gets K}$ is computationally expensive and depends more strongly on the polynomial degree $p$: the cost of a matrix inversion is $\mathcal{O}((p+1)^{3d})$ whereas matrix-vector products scale with $\mathcal{O}((p+1)^{2d})$ in $d$ dimensions. For higher polynomial degrees this will lead to a more favourable computation/communication ratio which improves scalability. Interestingly, if $A_{K\gets K}^{-1}$ is precomputed, scaling breaks down earlier for higher polynomial degrees. If $A_{K\gets K}$ is inverted in each cell the opposite behaviour is observed. As a consequence, the performance of the setup where $A_{K\gets K}^{-1}$ is precomputed once can be (almost) matched by the implementation which inverts $A_{K\gets K}$ in each grid cell, provided the code is run on a larger number of compute cores.

\paragraph{Experimental data for shifted and fused implementation.}

\begin{figure}[htb]
  \begin{center}
    \includegraphics[width=0.48\textwidth]{\figdir/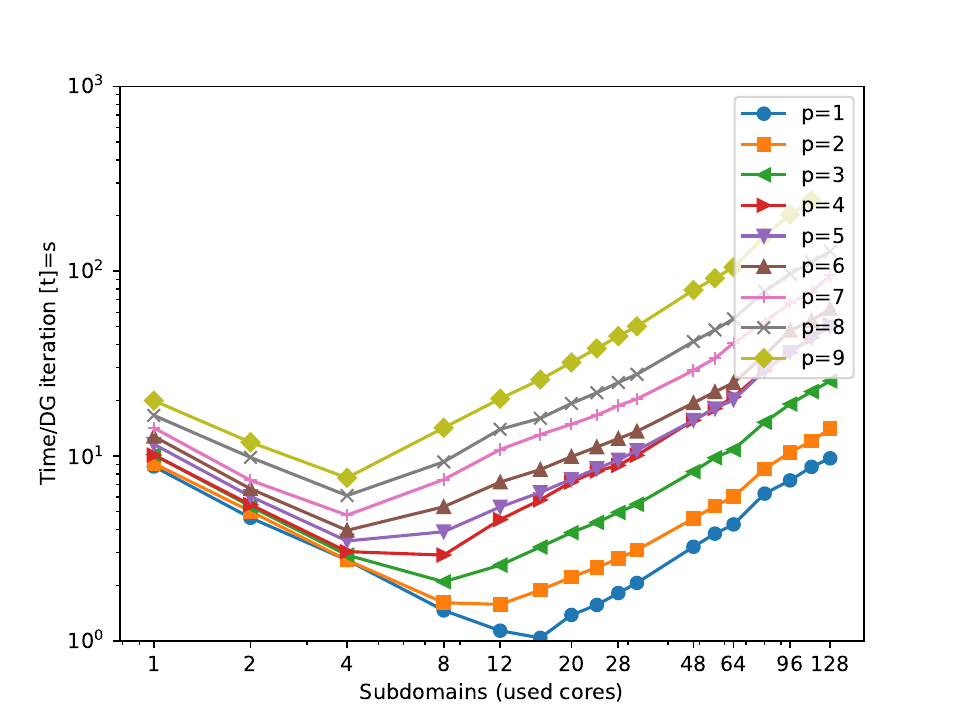}
    \includegraphics[width=0.48\textwidth]{\figdir/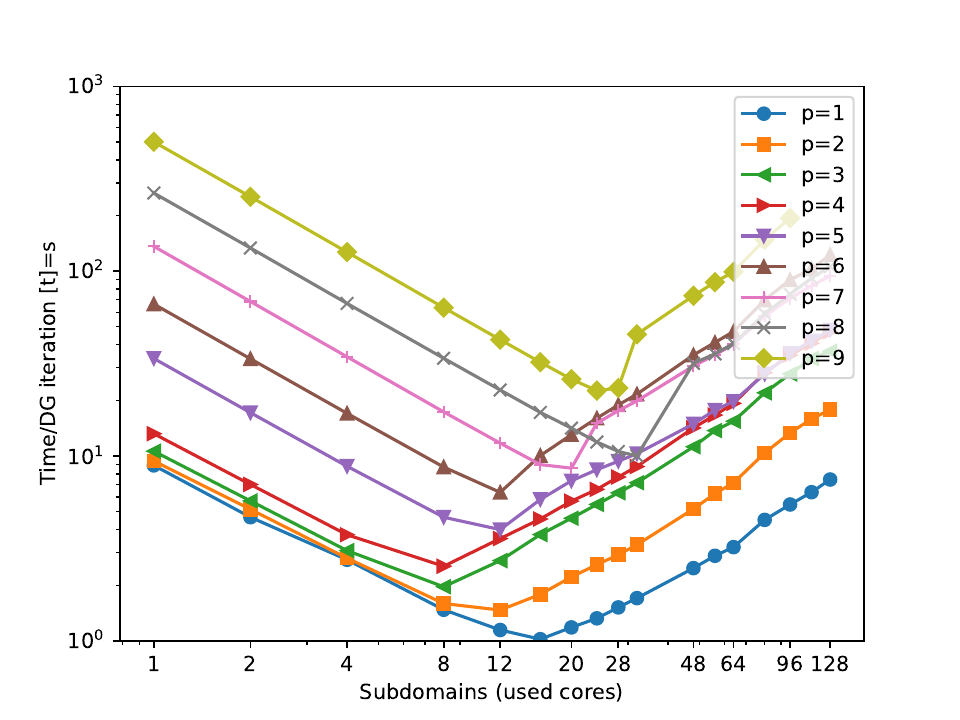}
  \end{center}
  \caption{
    Scalability of \cref{alg:interior-penalty-single-level-one-mesh-sweep} with precomputed $A_{K\gets K}^{-1}$ (left) and inversion of $A_{K\gets K}$ in each mesh cell (right) on an AMD EPYC node. In contrast to \cref{figure:single-level-scalability:sweeps}, only on mesh traversal is required per block-Jacobi iteration.
    \label{figure:single-level-scalability:shifted}
  }
\end{figure}

The techniques used in \cref{alg:interior-penalty-single-level-one-mesh-sweep} reduce the runtime significantly (compare \cref{figure:single-level-scalability:shifted} to the corresponding \cref{figure:single-level-scalability:sweeps}).
This is consistent with the single-core results (\cref{table:single-level-hw-counters-AMD-sweeps}). The benefit is most pronounced for the setup where $A_{K\gets K}^{-1}$ is precomputed once at the beginning of the run.
Here, the runtime is easily reduced by a factor three. If $A_{K\gets K}$ is inverted in every grid cell, the gain is significantly smaller for higher polynomial degrees. This is not surprising since in this case most of the time is spent in the local matrix inversion which will not benefit from loop fusion. Qualitatively, the curves in \cref{figure:single-level-scalability:sweeps} and \cref{figure:single-level-scalability:shifted} are very similar. However, as expected, \cref{alg:interior-penalty-single-level-one-mesh-sweep} shows slightly better scaling for the lowest polynomial degrees where the computation/communication ratio is expected to be particularly poor.

\paragraph{Discussion.}
We empirically confirmed that our DG implementation is well-suited for parallelisation with domain decomposition. The introduction of auxilliary variables requires only the exchange of lower-dimensional data on the facets, which improves scalability. Comparing the results obtained with \cref{alg:interior-penalty-single-level-multiple-sweeps-and-auxiliary-variables} and \cref{alg:interior-penalty-single-level-one-mesh-sweep}, it is enouraging to see that for lower polynomial degrees loop-fusion significantly improves performance of the smoother, which is likely to be the bottleneck in the multigrid algorithm. However, we also observe that for the fixed problem size chosen here strong scaling breaks down once we use at the order of 10 subpartitions per node. 
The scaling improves for larger problem sizes, i.e.~scales weakly in the sense that the turnaround point shifts towards higher core counts (not shown).
We continue to investigate to which degree task-based parallelism can shift this turnaround point, too.


\subsection{Task-based realisation over well-balanced domain decompositions with homogeneous workload}

\cref{alg:interior-penalty-single-level-tasks} can be executed in parallel through a combination of standard non-overlapping domain decomposition and task-based asynchronous processing: the evaluation of the cell-based residual $\vec{r}^{(c)}_K = \vec{b}^{(c)} - A_{c\gets c}|_{K\gets K}\vec{u}^{(c)}$ and the computation of $A_{K\gets K}^{-1}$ (in setups where the inverse of $A_{K\gets K}$ is re-computed in each cell) have been outsourced into tasks that are executed by the runtime system. All other operations, such as the computation of numerical fluxes, facet-based residual updates, solution updates and projections are executed in a standard mesh-traversal which is performed independently in the different subdomains.

It spawned ``Cell-residual''-type and ``Matrix-inversion''-type tasks on top of the domain decomposition introduce two levels of parallelism increasing the total concurrency, yet require careful balancing: we need to decide how many cores to dedicate to the mesh-traversal and how many cores are reserved for the asynchronous execution of the spawned ``Cell-residual''- and ``Matrix-inversion''-type tasks. 
The distribution of cores between the two different modes of execution (grid-traversal and task-processing) can be controlled by varying the number of subdomains, each of which is traversed by a single compute core, while using the remaining cores on a full node to execute the spawned tasks. 
It is natural to ask whether for a particular distribution it is possible to beat the performance of \cref{alg:interior-penalty-single-level-one-mesh-sweep} without tasking.

\begin{figure}[htb]
  \begin{center}
    \includegraphics[width=0.48\textwidth]{\figdir/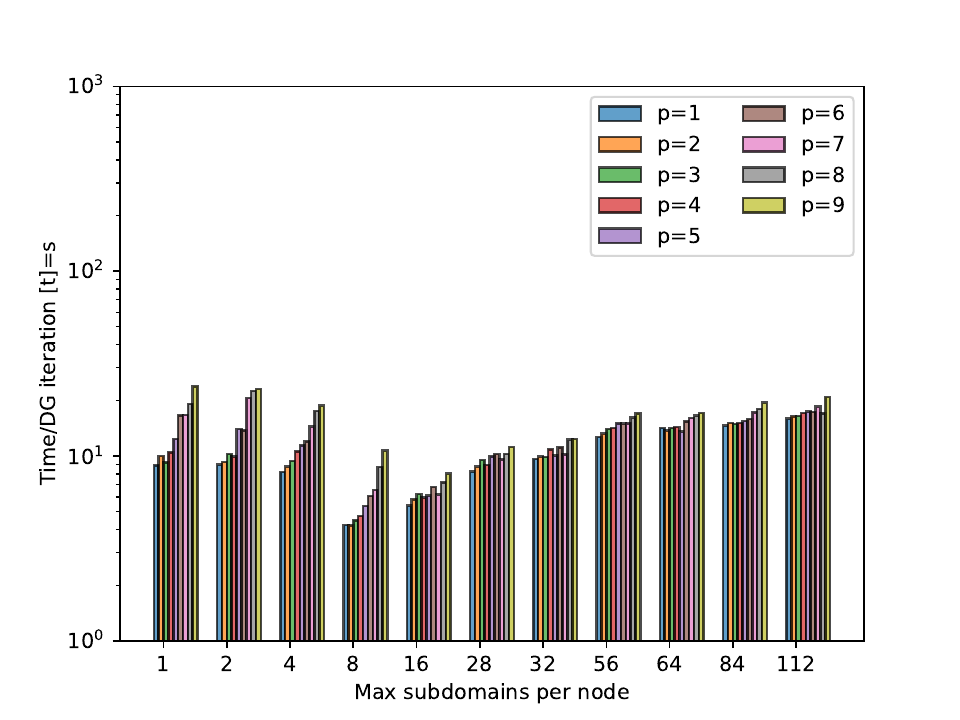}
    \includegraphics[width=0.48\textwidth]{\figdir/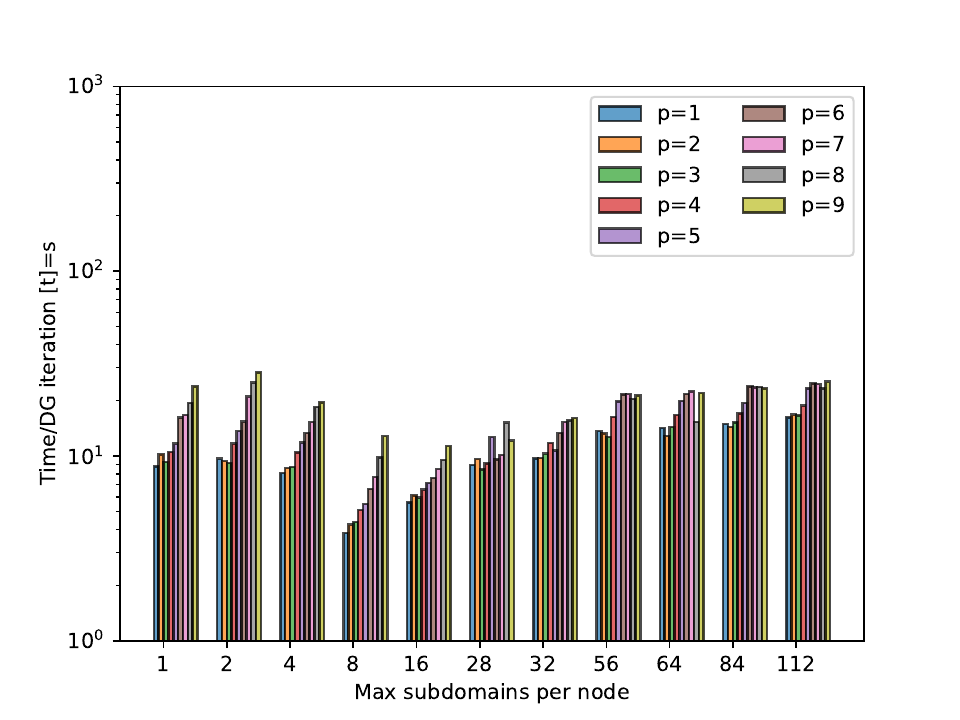}
  \end{center}
  \caption{
    Performance of \cref{alg:interior-penalty-single-level-tasks} on a full 128-core AMD EPYC node with precomputed $A_{K\gets K}^{-1}$ (left) and inversion of $A_{K\gets K}$ in every grid cell (right). The number of subdomains (and therefore the number of tasks responsible for traversing the mesh) increases from 1 to 112.	
    \label{figure:tasks:perfect-lb:tasks}
  }
\end{figure}

\paragraph{Experimental data}
\cref{figure:tasks:perfect-lb:tasks} shows the results of a numerical experiment in which we assigned different numbers of compute cores to the mesh traversal in \cref{alg:interior-penalty-single-level-tasks}, while again balancing the size of the resulting subdomains along the space-filling curve perfectly. Compared to the results in \cref{figure:single-level-scalability:shifted}, which are obtained with the pure domain-decomposition approach to parallelisation in \cref{alg:interior-penalty-single-level-one-mesh-sweep}, the total runtime shown in \cref{figure:tasks:perfect-lb:tasks} depends only weakly on the polynomial degree $p$. The optimal choice of subdomains does not depend strongly on the polynomial degree either, and typically about 8-16 subdomains lead to optimal results. 
For a given polynomial degree, performance is not very sensitive to the exact choice of the number of subdomains. For small $p$, the performance of \cref{alg:interior-penalty-single-level-tasks} is a factor of $3\times$ - $4\times$ worse than the optimal result obtained with \cref{alg:interior-penalty-single-level-one-mesh-sweep}. In contrast, for the largest polynomial degree $p=9$ running the task-based \cref{alg:interior-penalty-single-level-tasks} with eight subdomains improves the performance by a factor of $2\times$ - $3\times$ in the setup where the inverse of $A_{K\gets K}$ is re-computed in each grid cell.

\paragraph{Discussion.}
The hybrid execution model in \cref{alg:interior-penalty-single-level-tasks} avoids overwhelming the runtime system with a large number of small tasks \cite{Tuft:2024:DetrimentalTaskPatterns}. 
Instead, the very small tasks are merged into the mesh traversal and hence executed immediately without spawning overhead. The remaining two types of tasks, which are executed asynchronously, are computationally expensive and will therefore incur a much smaller relative overhead from task-management. Nevertheless, a drawback of any task parallelism is the potential loss of cache coherency and data affinity: the mesh traversals update the solution and spawn ``Cell-residual'' and ``Matrix-inversion''-type tasks, i.e.~have to read data into their caches, but then these tasks might be executed by different cores, which require data transfers. The numerical results in \cref{figure:tasks:perfect-lb:tasks} show that performance can be optimised by balancing the number of processors assigned to different types of tasks. It might be coincidence that the optimal number of subdomains is identical to the  number of NUMA domains, as we do not explicitly pin the producer tasks to the respective domains. 
Overall however, the tasking fails to improve the performance signifcantly over a well-balanced data decomposition. 

The robustness of \cref{alg:interior-penalty-single-level-tasks} with respect to the number of subdomains however suggests that it is well-prepared to compensate for imbalances in the domain decomposition.
The spawned ``Cell-residual''-type tasks all have a very similar computational cost. 
The same applies to the ``Matrix-inversion''-type tasks.
While we continue to focus on geometric imbalances, complex domains requiring adaptive numerical integration can introduce imbalanced, too \cite{murray2021assembly}, and we therefore expect to end up processing a large number of very similar tasks. 

\subsection{Imbalanced domain decomposition}
We also study the performance of \cref{alg:interior-penalty-single-level-tasks} for a scenario in which the subdomains have significantly different sizes.
This can arise for example when the mesh is refined dynamically guided by an error estimator.
To control the amount of spatial imbalance in local domain size, we consider a mesh with $N=531,441$ cells and partition it such that the first subdomain consists of $N/2$ cells, the second subdomain contains $N/4$ cells, the third $N/8$ cells and so forth; the final subdomain consists of all cells that have not been distributed in this way yet.
If the parallelisation strategy is purely based on domain-decomposition, we would expect a speedup of no more than two for this artificially imbalanced setup. 
\begin{figure}[htb]
  \begin{center}
    \includegraphics[width=0.48\textwidth]{\figdir/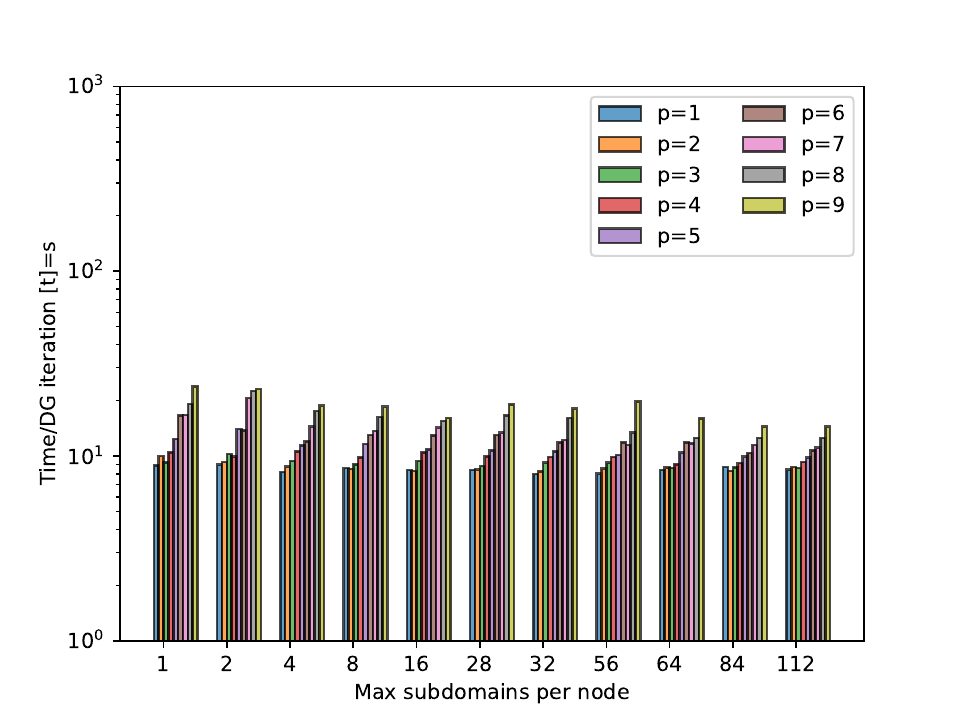}
    \includegraphics[width=0.48\textwidth]{\figdir/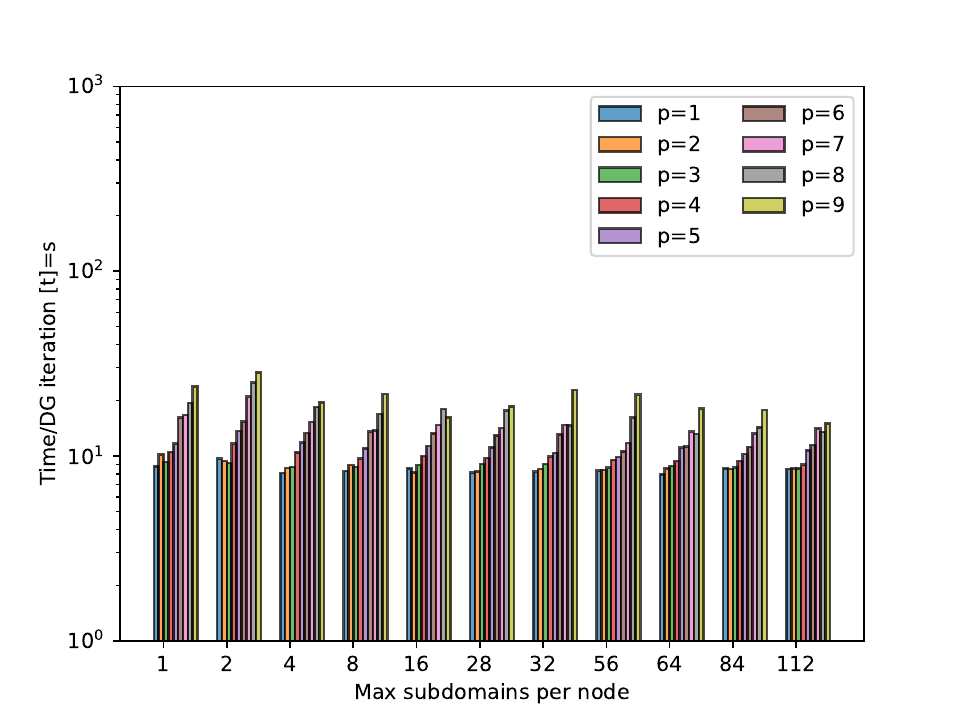}
  \end{center}
  \caption{
    Performance of \cref{alg:interior-penalty-single-level-tasks} when executed on a full AMD EPYC node with precomputed $A_{K\gets K}^{-1}$ (left) and inversion of $A_{K\gets K}$ in each grid cell (right) for an artificially imbalanced domain decomposition. As for \cref{figure:tasks:perfect-lb:tasks}, one core is assigned to each subdomain while the remaining cores simultaneously execute the spawned ``Cell-residual''- and ``Matrix-inversion''-type tasks.
    \label{figure:tasks:imbalanced:tasks}
  }
\end{figure}
\paragraph{Experimental data.}
However, this can be improved by using tasking from \cref{alg:interior-penalty-single-level-tasks} to leverage additional concurrency (\cref{figure:tasks:imbalanced:tasks}).
Even for the very imbalanced domain decomposition, the performance of \cref{alg:interior-penalty-single-level-tasks} is virtually independent of the number of subdomains. The additional concurrency provided by the asynchronous execution of some tasks compensates for the fact that the domain decomposition is (artificially) inefficient. 
Despite the impact of tasking, the imbalance continues to have a negative effect on overall performance, notably for smaller polynomial degrees $p$.
\paragraph{Discussion and implications.}
The overhead from task management is still expensive, and only large per-cell workloads allow us to compensate for this. 
However, for sufficiently large $p$, it is reasonable to employ tasking once the domain decomposition is exhausted.

\section{Conclusion}
\label{section:conclusion}

The scientific computing community suffers from a lack of papers on implementation idioms that help scientists translate original ideas into working code. 
Our work helps bridge this gap by introducing explicit techniques to address particular challenges of DG implementations in the context of multigrid algorithms. 
As is common with such implementation techniques, they do not uniformly pay off.
In some situations, they are of great value; in others, they are detrimental to performance. 
A prime example of such techniques is the use of a task paradigm within the present work: 
while it is beneficial to phrase complex calculations in a task language, the construction of an efficient execution schedule requires careful attention, and the powerful modelling technique does not necessarily always manifest in a task-based implementation that performs well.
Exploring this implementation space is highly context-dependent yet benefits from a formal write-up of implementation techniques and observations, as it allows us to combine individual ingredients more systematically.
We expect hence collections of techniques to unfold their full impact as part of other implementations.

There are two natural follow-up directions for the present research. 
On the one hand, our work does not implement a particularly sophisticated multigrid flavour. 
Modern $hp$-multigrid typically combines various function spaces and refrains from overly aggressive $p$-coarsening,
combines algebraic operator construction with geometric multigrid principles, employs more sophisticated smoothers, and pays particular attention to the design of the actual multilevel cycle. 
While the implementation techniques presented here streamline the development of such more sophisticated multigrid variants, they will in turn give raise to implementation challenges and hence lead to new techniques.
A prime example is the discussion around smoothers.
Indeed, the literature suggests that more complex PDEs require actually more complex, sophisticated block operators rather than simpler, cell-local approximations. Vertex patch smoothers for the ill-conditioned bi-harmonic problem $-\Delta^2 u(x) = b(x)$ are discussed and optimised in \cite{witte2025tensor}. The authors of \cite{farrell2021pcpatch} describe a general framework for patch-based multigrid smoothers for linear elasticity problems and the Stokes- and Semilinear Allen-Cahn equations; similarly, monolithic multigrid smoothers are also applied to the Stokes problem in in \cite{rafiei2025achieving}.
It is not clear how such numerical developments can be translated into techniques facilitating efficient implementations in the tradition of \cref{def:vectorisation}--\cref{def:weak_single_touch}.

On the other hand, we think that our techniques in themselves open the door to totally new numerics and implementation flavours. 
A leading motif behind our techniques is the decoupling and localisation of solution updates. 
While we use this ambition to minimise memory transfers, we note that the decoupling also facilitates totally decoupled, asynchronous solvers \cite{Yamazaki2019Async,Wolfson2025Asynchronous}: 
as we hold data representations such as $\vec{u}^{(+)}$ and $\vec{u}^{(-)}$ redundantly between subdomains and as an auxiliary yet first-class data structure, it is possible to let individual subdomains iterate independently of each other. However it remains unclear how such asynchronicity translates into a multiscale setup. 
The localisation and atomic character of the individual tasks furthermore facilitates flexible load balancing, including the offloading to accelerators, while we assume that the approach is of great value for resilient algorithms.

\begin{acks}
Our work has been supported by the \grantsponsor{id-epsrc}{Engineering and Physical Sciences research Council (EPSRC)}{https://www.ukri.org/councils/epsrc/} through Grant Nos. \grantnum{id-epsrc}{EP/W026775/1} and \grantnum{id-epsrc}{EP/X019497/1}. Software development relied on the DiRAC@Durham facility managed by the Institute for Computational Cosmology on behalf of the STFC DiRAC HPC Facility
(\href{www.dirac.ac.uk}{www.dirac.ac.uk}); DiRAC is part of the National e-Infrastructure. The equipment was funded by BEIS capital funding via the \grantsponsor{id-stfc}{Science and Technology Facilities Council (STFC)}{https://www.ukri.org/councils/stfc/} through grant numbers \grantnum{id-stfc}{ST/K00042X/1}, \grantnum{id-stfc}{ST/P002293/1}, \grantnum{id-stfc}{ST/R002371/1}, \grantnum{id-stfc}{ST/S002502/1} and \grantnum{id-stfc}{ST/R000832/1}. Numerical experiments and performance measurements also made use of the facilities of the Hamilton HPC Service of Durham University. The authors would like to express their particulars thanks to Intel's Academic Centre of Excellence at Durham University.
\end{acks}



\appendix
\section{Preconditioned versus unpreconditioned residual as error estimator}\label{sec:precond_vs_unprecond}
In \cref{section:multigrid} we argued that for ill-conditioned problems the preconditioned residual in \eqref{eqn:preconditioned_residual} is a better proxy for the error than $\vec{r}^{(c)} = \vec{b}^{(c)}-A\vec{u}^{(c)}$. In the following we provide numerical evidence for this claim.

\begin{figure}[h]
\centering
 \includegraphics[width=0.6\textwidth]{\figdir/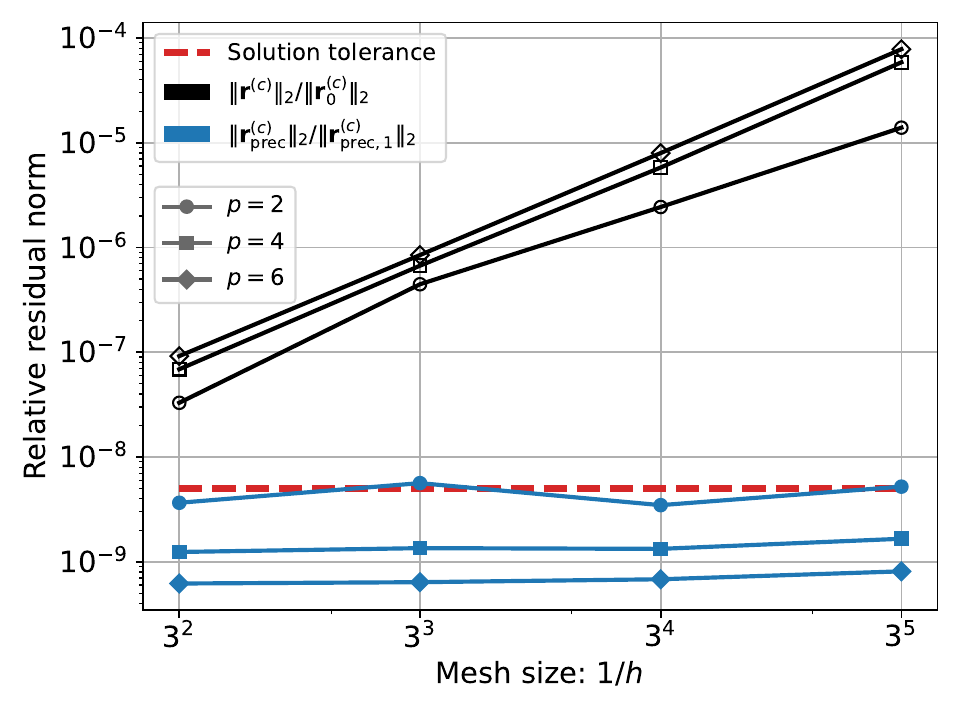}
 \caption{
   Preconditioned residual norm \mbox{$\| \vec{r}^{(c)}_{\text{prec}} \|_{2} / \| \vec{r}^{(c)}_{\text{prec},1} \|_{2}$} and unpreconditioned residual norm \mbox{$\| \vec{r}^{(c)} \|_{2} / \| \vec{r}^{(c)}_0 \|_{2}$} for different resolutions and polynomial degrees $p$. The red dashed line shows the upper bound on the relative error \mbox{$\| \widehat{\vec{u}}^{(c)} \|_{2} / {\| \vec{u}_{\textrm{ini.}}^{(c)} \|}_{2}$}.
   \label{fig:er_vs_ez}
  }
\end{figure}

Let us consider the linear problem $A \vec{u}^{(c)} = 0$ with the exact solution $\vec{u}^{(c)}_{\text{exact}} = 0$. Now assume that the same problem is solved with the $hp$-multigrid method in \cref{alg:multiplicative_multigrid}, with an initial guess $\vec{u}_{\textrm{ini.}}^{(c)}$ given by the ``two-peak'' reference function \eqref{eq:results:two_peak_exact} evaluated at the nodal points. The iteration is aborted once the solution error satisfies \mbox{$\| \widehat{\vec{u}}^{(c)} \|_{2} / {\| \vec{u}_{\textrm{ini.}}^{(c)} \|}_{2} \le 5 \times 10^{-9}$}. For the approximate solution $\widehat{\vec{u}}^{(c)}$ that is obtained in this way we can now compute the corresponding unpreconditioned residual $\vec{r}^{(c)} = \vec{b}^{(c)}-A\widehat{\vec{u}}^{(c)}$ and the preconditioned residual $\vec{r}^{(c)}_{\text{prec}} :=\widehat{\vec{u}}^{(c)}-\widehat{\vec{u}}^{(c)}_{\text{old}}$, where $\widehat{\vec{u}}^{(c)}_{\text{old}}$ is the numerical solution in the penultimate iteration, see \eqref{eqn:preconditioned_residual}. In \cref{fig:er_vs_ez} the relative norms \mbox{$\| \vec{r}^{(c)} \|_{2} / \| \vec{r}^{(c)}_0 \|_{2}$}  and \mbox{$\| \vec{r}^{(c)}_{\text{prec}} \|_{2} / \| \vec{r}^{(c)}_{\text{prec},1} \|_{2}$} are plotted for different values of the grid spacing $h$ and the polynomial degree $p$. While the value of the relative preconditioned residual norm \mbox{$\| \vec{r}^{(c)}_{\text{prec}} \|_{2} / \| \vec{r}^{(c)}_{\text{prec},1} \|_{2}$} is on the order of $10^{-9}\text{--}10^{-8}$ and in the same ballpark as the relative error \mbox{$\| \widehat{\vec{u}}^{(c)} \|_{2} / {\| \vec{u}_{\textrm{ini.}}^{(c)} \|}_{2}\lesssim 5 \times 10^{-9}$}, the relative unpreconditioned residual norm \mbox{$\| \vec{r}^{(c)} \|_{2} / \| \vec{r}^{(c)}_0 \|_{2}$} increases for larger problem sizes and can be several magnitudes larger than the relative error itself. As a consequence, using the unpreconditioned residual norm in the exit criterion would result in unnecessarily many iterative solver iterations.

\end{document}